\documentclass[12pt]{article}
\usepackage{amsmath,amsthm}
\usepackage{txfonts}
\usepackage{amscd}
\usepackage[english]{babel}
\usepackage[all]{xy}

\usepackage{pdfsync}

\usepackage{hyperref}

\usepackage{fancyhdr}
\usepackage[active]{srcltx}

\usepackage{makeidx}
\newtheoremstyle{cambiado}% name
  {3pt}%      Space above, empty = `usual value'
  {3pt}%      Space below
  {\itshape}% Body font
  {}%         Indent amount (empty = no indent, \parindent = para indent)
  {\bfseries}% Thm head font
  {\!:}%        Punctuation after thm head
  {.5em}%     Space after thm head: " " = normal interword space;
  {\arabic{theorem}. #1 #3}% Thm head spec

\theoremstyle{cambiado}

\newtheorem{theorem}{Teorema}[section]
\newtheorem{teorema}[theorem]{Theorem}
\newtheorem{lemma}[theorem]{Lemma}
\newtheorem{proposicion}[theorem]{Proposition}
\newtheorem{corolario}[theorem]{Corollary}

\newtheoremstyle{naida}% name
  {3pt}%      Space above, empty = `usual value'
  {3pt}%      Space below
  {\itshape}% Body font
  {}%         Indent amount (empty = no indent, \parindent = para indent)
  {\bfseries}% Thm head font
  {:}%        Punctuation after thm head
  {.5em}%     Space after thm head: " " = normal interword space;
  {\arabic{theorem}. #3}% Thm head spec

\theoremstyle{naida}

\newtheorem{naida}[theorem]{}

\newtheoremstyle{vacio}% name
  {3pt}%      Space above, empty = `usual value'
  {3pt}%      Space below
  {}% Body font
  {}%         Indent amount (empty = no indent, \parindent = para indent)
  {\bfseries}% Thm head font
  {}%        Punctuation after thm head
  {.5em}%     Space after thm head: " " = normal interword space;
  {\arabic{theorem}. #3}% Thm head spec

\theoremstyle{vacio}

\newtheorem{vacio}[theorem]{}

\newtheoremstyle{definition2}% name
  {3pt}%      Space above, empty = `usual value'
  {3pt}%      Space below
  {}% Body font
  {}%         Indent amount (empty = no indent, \parindent = para indent)
  {\bfseries}% Thm head font
  {}%        Punctuation after thm head
  {.5em}%     Space after thm head: " " = normal interword space;
  {\arabic{theorem}. #1\,:}% Thm head spec

\theoremstyle{definition2}

\newtheorem{ejemplo}[theorem]{Example}
\newtheorem{ejemplos}[theorem]{Examples}
\newtheorem{definicion}[theorem]{Definition}
\newtheorem{ejercicio}[theorem]{Exercise}
\newtheorem{notacion}[theorem]{Notation}
\newtheorem{Nota}[theorem]{Note}

\newtheoremstyle{definition3}% name
  {3pt}%      Space above, empty = `usual value'
  {3pt}%      Space below
  {}% Body font
  {}%         Indent amount (empty = no indent, \parindent = para indent)
  {\bfseries}% Thm head font
  {}%        Punctuation after thm head
  {.5em}%     Space after thm head: " " = normal interword space;
  {\arabic{theorem}. #1\,:}% Thm head spec

\theoremstyle{definition3}

\newtheorem{observacion}[theorem]{Remark}
\newtheorem{observaciones}[theorem]{Remarks}
\numberwithin{equation}{section}

\newcommand{\teor}{\begin{teorema}}
\newcommand{\eteor}{\end{teorema}}
\newcommand{\prop}{\begin{proposicion}}
\newcommand{\eprop}{\end{proposicion}}
\newcommand{\coro}{\begin{corolario}}
\newcommand{\ecoro}{\end{corolario}}
\newcommand{\lema}{\begin{lemma}}
\newcommand{\elema}{\end{lemma}}
\newcommand{\ejer}{\begin{ejercicio}}
\newcommand{\eejer}{\end{ejercicio}}
\newcommand{\obse}{\begin{observacion}}
\newcommand{\obses}{\begin{observaciones}}
\newcommand{\eobse}{\end{observacion}}
\newcommand{\eobses}{\end{observaciones}}
\newcommand{\defi}{\begin{definicion}}
\newcommand{\edefi}{\end{definicion}}
\newcommand{\demo}{\begin{proof}[Proof]}
\newcommand{\edemo}{\end{proof}}
\newcommand{\ejem}{\begin{ejemplo}}
\newcommand{\eejem}{\end{ejemplo}}
\newcommand{\ejems}{\begin{ejemplos}}
\newcommand{\eejems}{\end{ejemplos}}
\newcommand{\nota}{\begin{notacion}}
\newcommand{\enota}{\end{notacion}}
\newcommand{\note}{\begin{Nota}}
\newcommand{\enote}{\end{Nota}}

\newcommand{\nada}{\begin{naida}}
\newcommand{\enada}{\end{naida}}
\newcommand{\vaci}{\begin{vacio}}
\newcommand{\evaci}{\end{vacio}}

\newcommand{\A}{{\mathcal A}}

\newcommand{\C}{{\mathcal C}}
\newcommand{\V}{{\mathcal V}}

\newcommand{\M}{{\mathcal M}}
\newcommand{\Nc}{{\mathcal N}}
\newcommand{\Ic}{{\mathcal I}}
\newcommand{\Jc}{{\mathcal J}}

\newcommand{\U}{{\mathcal U}}
\newcommand{\Kc}{{\mathcal K}}

\newcommand{\OO}{{\mathcal O}}
\newcommand{\NN}{{\mathbb N}}
\newcommand{\ZZ}{{\mathbb Z}}

\newcommand{\pp}{{\mathfrak p}}
\newcommand{\qq}{{\mathfrak q}}
\newcommand{\rr}{{\mathfrak r}}

\newcommand{\punto}{{\displaystyle \cdot}}
\newcommand{\dosflechasa}[3][]{\xymatrix@1{\ar@<1ex>[r]^-{#2}
\ar@<-1ex>[r]_-{#3} & }}
\newcommand{\dosflechas}{{\dosflechasa[]{}{}}}
\newcommand{\dosflechaso}[3][]{\xymatrix@1{\ar@<0.5ex>[r]^-{#2}
 & \ar@<0.5ex>[l]^-{#3} }}
\newcommand{\dosflechos}{{\dosflechos[]{}{}}}

\newcommand{\tresflechasa}[4][]{\xymatrix@1{
\ar@<-1ex>[r]_-{#4} \ar[r]|-{#3} \ar@<1ex>[r]^-{#2} & }}
\newcommand{\tresflechas}{\tresflechasa[]{}{}{}}
\newcommand{\cuatroflechasa}[5][]{\xymatrix@1{
\ar@<-1.5ex>[r]_-{#5} \ar[r]<-0.5ex>_-{#4} \ar@<0.5ex>[r]^-{#3}
\ar@<1.5ex>[r]^-{#2}& }}

\newcommand{\larga}{\,\,\raise  3.5pt\hbox{$\begin{CD} @>>>
\end{CD} $}\,\,}
\newcommand{\Larga}[1]{\,\,\raise #1\hbox{$\begin{CD} @>>>
\end{CD} $}\,\,}
\newcommand{\pamatrix}[1]{\begin{pmatrix} #1 \end{pmatrix}}
\newcommand{\xymatris}{\def\objectstyle{\scriptstyle}
\def\labelstyle{\scriptstyle}\xymatrix@C=8pt @R=8pt}

\DeclareMathOperator{\limi}{{lim}}
\newcommand{\ilim}[1]{\,\underset{#1}{\underset{\to}{\limi}}\,}
\newcommand{\plim}[1]{\,\underset{#1}{\underset{\leftarrow}{\limi}}\,}
\newcommand{\iny}{\hookrightarrow}

\newcommand{\suma}[2][]{\underset{#2}{\overset{#1}{\sum}}}
\newcommand{\oplusa}[2][]{\underset{#2}{\overset{#1}{\oplus}}}
\newcommand{\capa}[2][]{\underset{#2}{\overset{#1}{\cap}}}
\newcommand{\cupa}[2][]{\underset{#2}{\overset{#1}{\cup}}}
\newcommand{\proda}[2][]{\underset{#2}{\overset{#1}{\prod}}}
\newcommand{\coproda}[2][]{\underset{#2}{\overset{#1}{\coprod}}}

\newcommand{\alineas}[1]{\begin{array}{#1}}
\newcommand{\alinea}{\begin{array}{l}}
\newcommand{\ealinea}{\end{array}}
\newcommand{\ealineas}{\end{array}}
\newcommand{\enumera}{\begin{enumerate}}
\newcommand{\eenumera}{\end{enumerate}}
\newcommand{\wh}{\widehat}
\newcommand{\wt}{\widetilde}

\newcommand{\mroot}[2][]{\text{\ \,\,\font\rad =cmmi5 at
5pt\setbox1=\vbox{\rad 3}\hbox{ \lower 1.3mm\vtop to
0.1222pt{\hsize 2.41pt \moveleft 8.7pt \vbox{\hsize 0.7pt
\hbox{\vsize 2pt  /}\vskip \ht2 minus \ht1\setbox2=\vbox{\rad
#1}\nointerlineskip \vbox{\hsize 5pt\noindent\,\hfill\rad
#1}}}\raise 0.35mm\hbox {\vbox{\hsize 2pt\moveleft
5.5pt\vbox{\hsize 2pt\hbox{\text{$\backslash$}}}}}\vbox{\hsize
2pt\moveleft 1.5pt\vbox{\hsize
1pt\hbox{\underbar{{\,}}}}}}\underbar {\hbox{\raise
2pt\hbox{\text{$#2$\ }}}}}}

\DeclareMathSymbol{\ltimes}{\mathbin}{AMSb}{"6E}
\DeclareMathSymbol{\rtimes}{\mathbin}{AMSb}{"6F}

\DeclareMathSymbol{\nleq}{\mathbin}{AMSb}{"0A}
\DeclareMathSymbol{\functor}{\mathbin}{AMSa}{"20}
\DeclareMathSymbol{\subsetp}{\mathbin}{AMSb}{"24}

\DeclareMathOperator{\Hom}{{Hom}}

\DeclareMathOperator{\Spec}{{Spec}}

\DeclareMathOperator{\Coker}{{Coker}}

\DeclareMathOperator{\Ker}{{Ker}}

\DeclareMathOperator{\Ima}{{Im}}

\DeclareMathOperator{\rad}{{rad}}

 \DeclareMathOperator{\Id}{{Id}}

\DeclareMathOperator{\di}{{d}}

\makeindex

\NoCompileMatrices

\usepackage{multicol}
\usepackage{xcolor}

\begin{document}

\title{Notes on schematic finite spaces}
\author{F. Sancho and P. Sancho}
%\date{May 2020}

\maketitle

\begin{abstract} The schematic finite spaces are those finite ringed spaces where  a theory of quasi-coherent modules  can be developed with minimal natural conditions. We  give various characterizations of these spaces and  their natural morphisms. We show that schematic finite spaces are strongly related to quasi-compact  quasi-separated schemes.

\end{abstract}

\section*{Introduction}

In Topology, Differential Geometry and Algebraic Geometry, it is usual to study their geometric objects considering suitable finite open coverings and studying the associated finite ringed spaces. 
Let us remember how these finite ringed spaces are constructed:

\vaci \label{Parr}
\medskip
Let $S$ be a topological space and let   $\U=\{U_1,\dots,U_n\}$ be a finite open covering of $S$. For each $s\in S$ define $U_s:=\underset{s\in U_i} \cap U_i$.  Observe that the topology generated by $\U$ is equal to the topology generated by $\{U_s\}_{s\in S}$.
We shall say that $\U$ is a minimal open covering of $S$ if $U_i\neq U_j$ if $i\neq j$ and $U_i=U_s$ for some $s\in S$, for every $i$. Define the following equivalence relation on $S$: $s\sim s'$ if the covering $\U$ does not distinguish them, i.e., if $U_s=U_{s'}$. Consider on $S$ the topology generated by the covering $\U$ and let $X:=S/\sim$ be the quotient topological space. $X$ is a finite $T_0$-topological space, then
it is a finite poset as follows: $[s]\leq [s']$ if $U_{s'}\subseteq U_{s}$.
Let $\pi\colon S\to X$, $s\mapsto [s]$ be the quotient morphism and let $U_{[s]}=\{[s']\in S\colon [s']\geq [s]\}$ be the minimal open neighborhood of $[s]$. One has that $\pi^{-1}(U_{[s]})=U_s$. 
Suppose now that $(S,\OO_S)$ is a ringed space (a scheme, a differentiable manifold, an analytic space, etc.).
We have then a sheaf of rings on $X$, namely  $\OO:=\pi_*\OO_S$, so that $\pi\colon (S,\OO_S)\to (X,\OO)$ is a morphism of ringed spaces. We shall  say that $(X,\OO)$ is the  {\it ringed finite space associated with the finite covering $\U$}. Observe that $\OO_{[s]}=\OO(U_{[s]})=\OO_S(U_s)$.

To fix ideas, suppose that $ S $ is a quasi-compact quasi-separated scheme (see \cite{Gr} 1.2.1).  There exists a minimal affine open covering $\U=\{U_{s_1},\ldots,U_{s_n}\}$ of $S$ (see \cite{SanchoHomotopy} 3.13). Consider the associated ringed finite space 
$X$.
It is easy to prove that the functors $\mathcal M \functor \pi_*\M$, $\Nc\functor \pi^*\Nc$ stablish an equivalence of categories between the category of quasi-coherent $\OO_S$-modules
and the category of quasi-coherent $\OO_X$-modules. Besides,
$H^i(S,\M)=H^i(X,\pi_*\M)$ for any quasi-coherent $\OO_S$-module $\M$. Observe that for every $U_{s_j}, U_{{s_{j'}}}\subseteq  U_{s_i}$:
\enumera

\item[-] The restriction morphism $\OO_S(U_{s_i})\to \OO_S(U_{s_j})$ is a flat morphism, since  the morphism $\OO_S(U_{s_i})\to (\OO_S(U_{s_i})\backslash\pp)^{-1}\OO_S(U_{s_i})=\OO_{S,\pp}=(\OO_S(U_{s_j})\backslash\pp)^{-1}\OO_S(U_{s_j})$ is flat, for any  $\pp\in  U_{s_j}=\Spec \OO_S(U_{s_j})$.

\item[-] The natural morphism  $\OO(U_{s_j})\otimes_{\OO_S(U_{s_i})}\OO_S(U_{{s_{j'}}})\to \OO_S(U_{s_j}\cap U_{{s_{j'}}})$ is an isomorphism, because
$U_{s_j}\times_{U_{s_i}} U_{{s_{j'}}}=U_{s_j}\cap U_{{s_{j'}}}$.

\item[-]  The morphism $\OO_S(U_{s_j}\cap U_{{s_{j'}}})\to \prod_{U_{s_k}\subseteq U_{s_j}\cap U_{{s_{j'}}}} \OO_S(U_{s_k})$ is faithfully flat: it is flat because $U_{s_j}\cap U_{{s_{j'}}}$ and $U_k$  are affine and it is faithfully flat because 
$\coprod_{U_{s_k}\subseteq U_{s_j}\cap U_{{s_{j'}}}} U_k\to U_{s_j}\cap U_{{s_{j'}}}$ is a surjective map. 
\eenumera
Therefore, for any points $x_j,x_{j'}\geq x_i$ in $X$:

\enumera
\item[a.] The natural morphism $\OO_{x_i}\to \OO_{x_j}$ is flat.

\item[b.] The natural morphism  $\OO_{x_j}\otimes_{\OO_{x_i}}\OO_{x_{j'}}\to \OO(U_{x_j}\cap U_{x_{j'}})$ is an isomorphism.

\item[c.]  The morphism $\OO(U_{x_j}\cap U_{x_{j'}})\to \prod_{x_k\geq x_j, x_{j'}} \OO_{x_k}$ is faithfully flat.
\eenumera
\evaci

\vaci We shall  say that a ringed finite space is schematic if it satisfies a., b. and c.
In \cite{KS} 4.4, 4.11, it is proved  that a finite ringed space $X$ is schematic iff 
$X$ satisfies a. and $R^n\delta_*\OO_X$ is a quasi-coherent module for any $n$, where $\delta\colon X\to X\times X$ is the diagonal morphism.  In \cite{KG} 4.5, it is proved that  $X$ is schematic iff $X$ satisfies a. and $R^ni_*\OO_{U}$ is quasi-coherent 
for any open subset $U\overset i\subset X$ and for $n\in\NN$. In Algebraic Geometry, it is usual to approach the study of schemes and their morphisms through the category of quasi-coherent modules, for example, the theory of intersection can be studied with the $K$-theory of quasi-coherent modules. We shall  denote by 
${\bf Qc\text{-}Mod}_X$  the category of quasi-coherent $\OO_X$-modules.
We prove that a finite ringed space $X$ is schematic iff ${\bf Qc\text{-}Mod}_X$ satisfies minimal conditions:

\enumera

\item[-]   A ringed finite space $X$ is schematic iff for any morphism $f\colon \M\to \Nc$ of quasi-coherent $\OO_X$-modules $\Ker f$ is  quasi-coherent, and $\delta_*(\mathcal M)\in {\bf Qc\text{-}Mod}_{X\times X}$,  for any $\M\in {\bf Qc\text{-}Mod}_X$, where  $\delta\colon X\to X\times X$, $\delta(x)=(x,x)$ is the diagonal morphism. 

\item[-]  A ringed finite space  $X$ is schematic iff  for any morphism $f\colon \M\to \Nc$ of quasi-coherent $\OO_X$-modules $\Ker f$ is  quasi-coherent, and $i_*(\mathcal M)\in {\bf Qc\text{-}Mod}_X$  for any $\M\in {\bf Qc\text{-}Mod}_U$ and  any open subset  $U\overset i\iny X$. 
\eenumera

Likewise, we study and characterize  affine finite spaces. Let us use the previous notations.
It can be proved that $S$ is an affine scheme iff the morphisms
$\OO(S)\to \prod_{i}\OO_{s_i}$ and $\OO(U_{s_i})\otimes_{\OO(S)}\OO(U_{s_j}) \to \prod_{U_{s_k}\subset U_{s_i}\cap U_{s_j}} \OO(U_{s_k})$ are faithfully flat, for any $i,j$.
We say that a finite ringed space $ X $ is affine if\enumera 
\item[-] The morphism $\OO(X)\to \prod_{x\in X} \OO_x$ is faithfully flat.

\item[-] The morphism $\OO_{x}\otimes_{\OO(X)}\OO_{x'}\to \prod_{z\in U_{x}\cap U_{x'}} \OO_z$ is faithfully flat,  for any $ x,x'\in X$.\eenumera

\noindent Affine finite spaces are schematic and a finite ringed space $X$ is schematic iff $U_x$ is affine, for any $x\in X$.
We  prove that $X$ is affine iff $\OO(X)\to \OO_x$ is flat,  $\OO_X(U_x\cap U_{y})=\OO_x\otimes_{\OO(X)}\OO_y$ and $X$ and $U_x\cap U_{y}$ are acyclic, for any $x, y\in X$.
A schematic space $X$ is affine iff $H^1(X,\M)=0$ for any quasi-coherent module $\M$ (see \cite{KS} 5.11), which is equivalent to saying
that the functor $\M\functor \Gamma(X,\M)$ is exact (see Corollary \ref{corotonto}).\evaci

\vaci Next, we study the morphisms between schematic finite spaces.
Let $f\colon X\to Y$ be a morphism of ringed spaces between schematic finite spaces. We say that $f$
is affine if $f_*\OO_X$ is a quasi-coherent module and $f^{-1}(U)$ is affine, for any affine open subset $U\subset Y$. We prove that $f$ is affine  if and only if $f_*$ preserves quasi-coherence and it is an exact functor.

We say that 
$f$ is schematic if $f_*\OO_X$ is quasi-coherent and the morphism $U_x\to U_{f(x)}$, $x'\mapsto f(x')$ is affine, for any $x\in X$.
 We prove that the following statements are equivalent:

\enumera
\item[-] $f$ is schematic.

\item[-] $f_*$ preserves quasi-coherence.

\item[-] The natural flat morphism $\OO_x\otimes_{\OO_{f(x)}}\OO_y \to \prod_{z\in U_x\cap f^{-1}(U_y)} \OO_z$
 is faithfully flat, for any $x\in X$ and $y\geq f(x)$.

\item[-] $R^n\Gamma_{f*}\OO_X\in {\bf Qc\text{-}Mod}_{X\times Y}$ for any $n$, where $\Gamma_f\colon X\to X\times Y$, $\Gamma_f(x)=(x,f(x))$ is the graph of $f$.

\item[-] $R^n(f\circ i)_*\OO_U\in {\bf Qc\text{-}Mod}_X$, for   any open subset  $U\overset{i}\iny X$ and any $n\geq 0$.
\eenumera\evaci

\vaci Now, we ask ourselves whether schematic finite spaces are determined by the category of their quasi-coherent modules. A   morphism of schemes $F \colon S \to T$ between 
quasi-compact quasi-separated schemes is an isomorphism if and only if 
the functors $${\bf Qc\text{-}Mod}_S\,\, \dosflechaso{F_*}{F^*} \,\,
{\bf Qc\text{-}Mod}_T$$ are mutually inverse.  
Given a schematic morphism $f\colon X\to Y$, we prove that the following statements are equivalent
\enumera
\item[-] The functors ${\bf Qc\text{-}Mod}_X\,\, \dosflechaso{f_*}{f^*} \,\,
{\bf Qc\text{-}Mod}_Y$ are mutually inverse.

\item[-] $\OO_Y=f_*\OO_X$ and $f$ is affine.

\item[-] Essentially, $f$ is the quotient morphism defined on $X$ by a minimal  affine open covering.

\item[-] The cylinder $C(f)=X\coprod Y$ of $f$ is a schematic space and $f$ is a faithfully flat morphism.
\eenumera 
\noindent A morphism that satisfies any of these statements will be called   quasi-isomorphism.\evaci

\vaci \label{Pa2}  Let us talk less accurately.
Given a schematic finite space $X$, consider the ringed space 

 $$\tilde X:=\ilim{x\in X} \Spec \OO_x.$$ 
  $\Spec \OO_x$ is a subspace of $\tilde X$ via the natural morphism $i_x\colon \Spec \OO_x\to \tilde X$ and $\tilde X=\cup_{x\in X} \Spec \OO_x$. 
$\tilde X$ is quasi-compact, the set of its quasi-compact open subsets is a basis of its topology and the intersection of two quasi-compact open subsets is quasi-compact.  Given a quasi-coherent $\OO_X$-module $\M$, let $\tilde \M_x$ be the $\OO_{\Spec \OO_x}$-module of localizations of $\M_x$ and consider the $\OO_{\tilde X}$-module $\tilde \M:=\plim{x\in X} i_{x,*} \tilde\M_x$. We prove that $H^n(X,\M)=H^n(\tilde X,\tilde{\M)}$, for any $n\geq 0$,  and 
 that the category of quasi-coherent modules of $X$ is equivalent to the category of quasi-coherent modules of $\tilde X$. Let $X$ and $Y$ be schematic finite spaces.
We say that a morphism of ringed spaces  $f\colon \tilde X\to \tilde Y$ is schematic if  $f_*$ preserves quasi-coherence.
Let $\C_W$  be the category of schematic finite spaces localized by the quasi-isomorphisms. We  prove that
$$\Hom_{\C_W}(X,Y)=\Hom_{sch}(\tilde X,\tilde Y).$$\evaci

\section{Finite ringed spaces: basic notions}

Let $X$ be a finite set. It is well known ([\ref{Alexandrof}]) that giving a topology on $X$ is equivalent to giving a preorder relation on $X$: 
$$x \leq y  \, \iff \, \bar{x} \subseteq \bar{y}, \ \ \text{ where } \bar{x}, \bar{y} \text{ are the closures of } x \text{ and } y.$$ In addition, the topology is $T_0$ if and only if the preorder is a partial order (i.e. it satisfies antisymmetric property). 

Let $X$ be a finite topological space. For each point $p \in X$, let us denote
$$ U_x  = \text{smallest open subset containing } x,$$
that is, $U_x= \{y \in X : y \geq x \}$.  Then, $ x \leq y  \Leftrightarrow U_y \subseteq U_x$. The family of open subsets $\{U_x\}_{x \in X}$ constitutes a minimal basis of open subsets of $X$ (any other base contains this one).

  A map $f\colon X \to Y$ between finite topological spaces is continuous if and only if it is monotone (i.e. $x \leq y$ implies $f(x) \leq f(y)$).

Let $X$ be a finite topological space and 
let $F$ be a sheaf of abelian groups (resp. rings, etc.) on $X$. The stalk of $F$ at $x\in X$, $F_x$, is an abelian group (resp. ring, etc.) and coincides
with the sections of $F$ on $U_x$. 
For each $x\leq y$, the natural morphism $r_{xy}\colon F_x\to F_y$ is just the restriction morphism  $F(U_x)\to F(U_y)$, which satisfies:
$r_{xx}=\Id$ for any $x$, and $r_{yz}\circ r_{xy}=r_{xz}$ for any $x\leq y\leq z$. 

 Conversely,  consider the following data:

 - An abelian group  (resp. a ring, etc) $F_x$ for each $x\in X$.

 - A morphism of groups (resp. rings, etc) $r_{xy}\colon F_x\to F_y$ for each $x\leq y$, satisfying: $r_{xx}=\Id$ for any $x$, and $r_{yz}\circ r_{xy}=r_{xz}$ for any $x\leq y\leq z$. 
 
\noindent Let $\mathcal F$ be the following presheaf of groups (resp. rings, etc.): For each open subset $U\subset X$,  
$\mathcal F(U):=\plim{x\in U} F_x$.
It is easy to prove that $\mathcal F_x=F_x$ and that $\mathcal F$ is a sheaf.

\defi A ringed space is a pair $(X,\OO)$, where $X$ is a topological space and  $\OO$ is a sheaf of (commutative with unit) rings on  $X$. A morphism of ringed spaces $(X,\OO)\to (X',\OO')$ is a pair $(f,f^\#)$, where $f\colon X\to X'$ is a continuous map and $f^\#\colon \OO'\to f_*\OO$ is a morphism of sheaves of rings (equivalently, a morphism of sheaves of rings $f^{-1}\OO'\to \OO$).
A {\it finite ringed space} is a ringed space $(X,\OO)$ whose underlying topological space $X$ is finite. 
\edefi

A morphism of ringed spaces $(X,\OO)\to (X',\OO')$ between two finite  ringed spaces is equivalent to the following data:

-  a continuous (i.e. monotone) map $f\colon X\to X'$,

-  for each  $x\in X$, a ring homomorphism  $f^\#_x\colon \OO'_{f(x)}\to \OO_x$, such that, for any  $x\leq y$, the diagram 
\[ \xymatrix{ \OO'_{f(x)} \ar[r]^{f^\#_{x}} \ar[d]_{r_{f(x)f(y)}} & \OO_{x}\ar[d]^{r_{xy}}\\ \OO'_{f(y)} \ar[r]^{f^\#_{y}}   & \OO_{y}}\] is commutative. We denote by $\Hom(X,Y)$ the set of morphisms of ringed spaces between two ringed spaces $X$ and $Y$.

\ejem Let $\{*\}$ be the topological space with one element. We denote by $(*,R)$ the finite ringed space whose underlying topological space is $\{*\}$ and the sheaf of rings is a ring $\OO_* =R$. For any ringed space $(X,\OO)$ there is a natural morphism of ringed spaces $(X,\OO) \to (*,\OO(X))$.
\eejem

Let $(X,\OO)$ be a finite ringed space. A sheaf $\M$ of $\OO$-modules (or $\OO$-module) is equivalent to these data: an $\OO_x$-module $\M_x$ for each $x \in X$, and a morphism of $\OO_x$-modules $r_{xy} \colon \M_x \to \M_y$ for each $x \leq y$, such that $r_{xx}=\Id$ and $r_{xz}=r_{yz} \circ r_{xy}$ for any $x \leq y \leq z$. Again, one has that
$$ \M_x = \text{stalk of } \M \text{ at } x = \M(U_x) $$
and $r_{xy}$ is the restriction morphism $\M(U_x) \to \M(U_y)$.

For each $x \leq y$ the morphism $r_{xy}$ induces a morphism of $\OO_y$-modules
$$ \widetilde{r_{xy}} \colon \M_x\otimes_{\OO_x}\OO_y \to   \M_y.$$
An $\OO$-module $\M$ is said to be quasi-coherent if for any $x\in X$ there exist an 
open neighbourhood $U$ of $x$ and an exact sequence of $\OO_{\vert U}$-modules
\[ \OO_{\vert U}^I \to \OO_{\vert U}^J\to\M_{\vert U}\to 0.\]

\teor[(\cite{SanchoHomotopy} 3.6)\,] \label{qc} Let $(X,\OO)$ be a finite ringed space. An $\OO$-module $\M$ is quasi-coherent if and only if for any  $x\leq y$ the morphism
\[\widetilde{r_{xy}} \colon \M_x\otimes_{\OO_x}\OO_y\to\M_y\] is an isomorphism.
\eteor

%{[\ref{SanchoHom},Theorem 2.6]}
\demo $\Rightarrow)$ Let $U$ be an open neighbourhood of $p$ such that  there exists an exact sequence
\[ \OO_{\vert U}^I \to \OO_{\vert U}^J\to\M_{\vert U}\to 0.\]
We can suppose $X=U$. Obviously,
$(\OO^I)_x\otimes_{\OO_x}\OO_y=(\OO_x)^I \otimes_{\OO_x}\OO_y=(\OO_x \otimes_{\OO_x}\OO_y)^I=(\OO_y)^I=(\OO^I)_y$ and 
$(\OO^J)_x\otimes_{\OO_x}\OO_y= (\OO^J)_y$, then
$\M_x\otimes_{\OO_x}\OO_y=\M_y$.

$\Leftarrow)$ Given $x\in X$, consider an exact sequence of $\OO_x$ modules $\OO_{x}^I \to \OO_{x}^J\to\M_{x}\to 0$.
Tensoring by $\otimes_{\OO_x}\OO_y$, for any $x\leq y$, one has the exact sequence $\OO_{y}^I \to \OO_{y}^J\to\M_{y}\to 0$.
Then, one has a sequence of morphisms
$$\OO_{|U_x}^I \to \OO_{|U_x}^J\to\M_{|U_x}\to 0$$
which is exact since it is exact on stalks at $y$, for any $y\in U_x$.
Therefore, $\M$ is quasi-coherent.
\edemo

\medskip

%\defi Let $(X,\OO)$ a ringed space.
 %Let $\M$ be an $\OO$-module (a sheaf of $\OO$-modules). We say that $\M$ is {\it quasi-coherent} if for each $x\in X$ there exist an open neighborhood  $U$ of $x$ and an exact sequence
%\[ \OO_{\vert U}^I \to \OO_{\vert U}^J\to\M_{\vert U}\to 0\] with $I,J$ arbitrary sets of indexes. Briefly speaking, $\M$ is quasi-coherent if %it is locally the  cokernel of a linear morphism between free modules. 
%\edefi

 %Let $\M$ be a sheaf of  $\OO$-modules on a ringed finite space $(X,\OO)$. Thus, for each $p\in X$, $\M_p$ is an $\OO_p$-module and for each  $p\leq q$ one has a morphism of $\OO_p$-modules $\M_p\to\M_q$, hence a morphism of  $\OO_q$-modules
%\[\M_p\otimes_{\OO_p}\OO_q\to\M_q.\]
We shall denote by $\text{\bf Mod}_X$ the category of $\OO$-modules on a ringed space $(X,\OO)$ and by ${\bf Qc\text{-}Mod}_X$ the subcategory of quasi-coherent $\OO$-modules. Also, for any ring $R$, we shall denote  by $\text{\bf Mod}_R$ the category of $R$-modules.

\obses 
\noindent a) If $f\colon X \to Y$ is a morphism of ringed spaces and $\Nc$ is a quasi-coherent $\OO_Y$-module, then $f^* \Nc:=f^{-1}\Nc\otimes_{f^{-1}\OO_Y}\OO_X$ is a quasi-coherent $\OO_X$-module. In particular, this is true for morphisms between finite ringed spaces.

 \noindent b) If $f\colon \mathcal \M \to \Nc$ is a morphism of $\OO_X$-modules where $\M$ and $\Nc$ are quasi-coherent, then $\Coker f$ is quasi-coherent. However, it is not always true that $\Ker f$ is quasi-coherent.
\eobses

\ejem Let $(X,\OO)$ be a finite ringed space and $\pi\colon (X,\OO) \to (*,\OO(X))$ the natural morphism of ring. If $M$ is an $\OO(X)$-module, then $\pi^* M$ is a quasi-coherent $\OO_X$-module, which we denote $\tilde M$. We say that $\tilde M$ is the quasi-coherent module associated with $M$ and we have a functor $\pi^*\colon \text{\bf Mod}_R \to {\bf Qc\text{-}Mod}_X, M \mapsto \tilde M$. Note that $\tilde M_x=M\otimes_{\OO(X)} \OO_x$, for each $x \in X$.
\eejem

\defi A finite ringed space $(X, \OO)$ is a \textit{finite flat-restriction space} (or \textit{finite fr-space}) if the restriction morphisms $r_{xy} \colon \OO_x\to \OO_y$ are flat, for any $x\leq y$.
\edefi

\prop \label{fr} Let $(X,\OO)$ be a finite ringed space. 

$(X,\OO)$ is a finite fr-space $\Leftrightarrow$ For any open subset $U$ (resp. $U_x$) of $X$ and any morphism $f \colon \M \to \Nc$ of quasi-coherent $\OO_U-$modules (resp. $\OO_{U_x}$-modules), $\Ker f$ is quasi-coherent.
\eprop
\demo 
$\Rightarrow)$ Let  $f \colon \M \to \Nc$ be a morphism of quasi-coherent $\OO_U-$modules. We have to prove that, for each $x \leq y \in U$, the morphism 
$$ \tilde r_{xy} \colon (\Ker f)_x \otimes_{\OO_x} \OO_y \to (\Ker f)_y$$
is an isomorphism. This follows from the next conmutative diagram of exact rows
$$\xymatrix{0 \ar[r] & (\Ker f)_x \otimes_{\OO_x} \OO_y  \ar[d]^{\tilde r_{xy}} \ar[r] & \M_x\otimes_{\OO_x} \OO_y  \ar[r]^-{f_x} \ar[d]^{\tilde r_{xy}} &
\Nc_x \otimes_{\OO_x} \OO_y  \ar[d]^-{\tilde r_{xy}}   \\ 0 \ar[r] & (\Ker f)_y \ar[r] & \M_y \ar[r]_-{f_y}  
& \Nc_y}$$ 
in which the first row is exact because $\OO_x \to \OO_y$ is a flat morphism and the second and the third vertical morphisms are isomorphisms because $\M$ and $\Nc$ are quasi-coherent $\OO_U$-modules.

\noindent $\Leftarrow)$ Let $x \leq y \in X$. Let $f_x \colon M_x \hookrightarrow N_x$ be an injective morphism of $\OO_x$-modules. We have to prove that $f_x \otimes 1 \colon M_x \otimes_{\OO_x} \OO_y \to N_x \otimes_{\OO_x} \OO_y $ is still injective. Consider the open subset $U_x$ of $X$ and the functor $\text{\bf Mod}_{\OO_x} \to {\bf Qc\text{-}Mod}_{U_x}, M_x \mapsto \widetilde{M_x}$. Then, the morphism $f_x$ gives us a morphism $\widetilde{f_x} \colon \widetilde{M_x} \to \widetilde{N_x}$ of quasi-coherent $\OO_{U_x}$-modules. Note that $(\widetilde{f_x})_y = f_x \otimes 1$. Since by hypothesis $\Ker \widetilde{f_x}$ is quasi-coherent, we have that:
$$\Ker (f_x \otimes 1)= (\Ker \widetilde{f_x})_y = \Ker f_x \otimes_{\OO_x} \OO_y = 0 \otimes_{\OO_x} \OO_y =0,$$
so we conclude that $f_x \otimes 1$ is injective.
\edemo

\obse Let $(X,\OO)$ be a finite ringed space. The proposition above says that ${\bf Qc\text{-}Mod}_U$ is an abelian category for each open subset $U$ of $X$ if and only if $(X,\OO)$ is a finite flat-restriction space. It is also true that in this case ${\bf Qc\text{-}Mod}_U$ is a Grothendieck category (see [\ref{Enochs}]).
\eobse

\section{Affine finite spaces}

\nota Let $(X,\OO)$ be a finite ringed space. For each $x,y\in X$, let us denote $U_{xy}=U_x\cap U_y$ and $\OO_{xy}=\OO(U_{xy})$. If $\M$ is an $\OO$-module, we denote $\M_{xy}= \M(U_{xy})$.
\enota

\defi \label{9} A finite ringed space $(X,\OO)$ is  called an \textit{affine (schematic) finite space} if it satisfies the following conditions: \enumera
\item $\OO(X)\to \prod_{x\in X} \OO_x$ is faithfully flat, for any $x \in X$.

\item $\OO_x\otimes_{\OO(X)}\OO_y=\OO_{xy}$, for any $x,y\in X$.

\item $\OO_{xy}\to \prod_{z\in U_{xy}} \OO_z$ is faithfully flat, for any $x,y\in X$.\eenumera \edefi

\prop If $(X,\OO)$ is an affine finite space, then it is a finite fr-space.
\eprop

\demo By condition 3. of the definition above, $\OO_{xx}=\OO_x  \to \prod_{z\in U_x} \OO_z$ is faithfully flat. Therefore,  $\OO_x \to \OO_z$ is flat, for any $z\geq x$. 
\edemo

\prop[(\cite{KS} 4.12)\,] \label{carafin} Let $(X,\OO)$ be a ringed finite space. $X$ is affine  iff 
\enumera
\item The morphism $\OO(X)\to \prod_{x\in X} \OO_x$ is faithfully flat.

\item The morphism $\OO_{y}\otimes_{\OO(X)}\OO_{y'}\to \prod_{z\in U_{yy'}} \OO_z$ is faithfully flat,  for any $ y,y'\in X$.\eenumera\eprop

\demo $\Rightarrow)$ It  follows immediately from the definition.
 
$\Leftarrow)$ In first place, note that  for any $x\leq u,u'$, the morphism $\OO_{u}\otimes_{\OO(X)} \OO_{u'}\to \OO_{u}\otimes_{\OO_x}\OO_{u'}$ is an epimorphism and the composite morphism
$$\OO_{u}\otimes_{\OO(X)} \OO_{u'}\to \OO_{u}\otimes_{\OO_x}\OO_{u'}\to \prod_{z\in U_{uu'}} \OO_z$$ is injective, because it is faithfully flat. Therefore, $\OO_{u}\otimes_{\OO(X)} \OO_{u'}=\OO_{u}\otimes_{\OO_x}\OO_{u'}$.

We only have to prove that 
$$ \OO_y \otimes_{\OO(X)} \OO_{y'}= \OO_{yy'}$$
Let us prove it by reduction to absurdity.
Let $y,y'\in X$ be maximal such that the morphism 
$ \OO_y \otimes_{\OO(X)} \OO_{y'}\to  \OO_{yy'}$ is not an isomorphism.

First, if $y \leq y'$, then $U_{yy'}=U_{y'}$ and the epimorphism $ \OO_y \otimes_{\OO(X)} \OO_{y'} \to \OO_{y'}$ is faithfully flat, by 2. Therefore,  $\OO_y \otimes_{\OO(X)} \OO_{y'}=\OO_{y'} =\OO_{yy'}$. So, neither $y \leq y'$, nor $y' \leq y$. The morphism $B:= \OO_{y}\otimes_{\OO(X)}\OO_{y'}\to \prod_{z\in U_{yy'}} \OO_z=:C$ is faithfully flat. Thus, the sequence of morphisms
$$(*)\qquad B\to C\dosflechas C\otimes_BC$$
is exact. By the maximality of $y$ and $y'$, given $z,z'\in U_{yy'}$,  $\OO_z\otimes_{\OO(X)}\OO_{z'}=\OO_{z}\otimes_{\OO_y}\OO_{z'}=\OO_{zz'}$. The natural morphism 
$\OO_z\otimes_{\OO(X)}\OO_{z'}\to \OO_z\otimes_B\OO_{z'}$ is surjective. The composite morphism $\OO_z\otimes_{\OO(X)}\OO_{z'}\to \OO_z\otimes_B\OO_{z'}\to \OO_{zz'}$ is an isomorphism, then $\OO_z\otimes_B\OO_{z'}=\OO_{zz'}$.
Therefore, $ C\otimes_BC=\prod_{z,z'\in U_{yy'}}\OO_{zz'}$. Then, from the diagram $(*)$, $B=\OO_{yy'}$ and we have come to contradiction.
\edemo

\coro \label{minimal affine} A finite ringed space  $(U_x,\OO)$ is  affine  iff the morphism $\OO_{y}\otimes_{\OO_x}\OO_{y'}\to \prod_{z\in U_{yy'}} \OO_z$ is faithfully flat,  for any $ y,y'\geq x$.

\ecoro
\demo It follows easily from the proposition above.
\edemo

\prop \label{afinfp} Let $X$ be an affine finite space and  $U\subset X$ an open set. Then, 
$U$ is affine iff  $\OO(U)\to \prod_{q\in U} \OO_q$ is a faithfully flat morphism.\eprop

\demo $\Rightarrow)$ $\OO(U)\to \prod_{q\in U} \OO_q$ is a faithfully flat morphism by definition of affine finite space. 

$\Leftarrow)$ We have to check that $U$ satisfies the conditions 2. and 3. of Definition \ref{9}. Condition 3. is clear, because $X$ is affine. Now, let us check 2.: for each $x,y \in U$, the morphism $\OO_x\otimes_{\OO(X)}\OO_{y}\to 
\OO_x\otimes_{\OO(U)} \OO_y$ is surjective and the composite morphism $$\OO_x\otimes_{\OO(X)}\OO_{y}\to 
\OO_x\otimes_{\OO(U)} \OO_{y}\to \OO_{xy}$$ is an isomorphism, thus $\OO_x\otimes_{\OO(U)} \OO_{y}\simeq \OO_{xy}$.\edemo

\coro \label{Uxy} If $X$ is an affine finite space,  then $U_{xy}$ is affine for every $x,y\in X$.\ecoro

\demo It follows from  condition 3. of Definition \ref{9} and the proposition above.\edemo

\prop  \label{11}  Let $X$ be an affine finite space and  $\mathcal M$ a quasi-coherent $\OO_X$-module. The natural morphism 
$$\M(V) \otimes_{\OO(X)}\OO(U)\to \M(U\cap V)$$
is an isomorphism, for any open set  $V\subseteq X$ and any affine open set  $U\subseteq X$.

\eprop

\demo  1.  The morphism $\OO(X)\to \prod_{x\in X}\OO_x=:B$ is  faithfully flat. The sequence of morphisms 
$$\OO(X)\to B= \prod_{x\in X}\OO_x\dosflechas \ B\otimes_{\OO(X)} B=\prod_{x,y\in X} \OO_{xy}$$
is a split sequence of morphisms  under a faithfully flat base change ($\OO(X)\to B$), then
this sequence of morphisms is universally exact, i.e., if we tensor the sequence of morphisms by $M\otimes_C-$ (where $C$ is a commutative ring, $M$ is a $C$-module and $\OO(X)$ a $C$-algebra) then we obtain an exact sequence of morphisms.
In particular, $\OO_{xy}\hookrightarrow \prod_{z\in U_{xy}} \OO_z$ is universally injective and the sequence of morphisms 
$$(*)\qquad \OO(X)\to \prod_{x\in X}\OO_x\dosflechas \prod_{x,y\in X, z\in U_{xy}} \OO_z$$
is universally exact. 

2.  Let  $W\subset U_x$ be an affine open set. Consider the universally exact sequence of morphisms
$$\OO(W)\to \prod_{z\in W} \OO_z\dosflechas \prod_{z,z'\in W, z''\in U_{zz'}} \OO_{z''}.$$
Tensoring by $\M_x\otimes_{\OO_x}-$, we obtain the exact sequence of morphisms
$$\M_x\otimes_{\OO_x} \OO(W)\to \prod_{z\in W} \M_z\dosflechas \prod_{z,z'\in W, z''\in U_{zz'}} \M_{z''},$$
which shows that $\M_x\otimes_{\OO_x} \OO(W)=\M(W) $.  Therefore (using Corollary \ref{Uxy}),
$$\mathcal M_x\otimes_{\OO(X)}\OO_y=\M_x\otimes_{\OO_x}\OO_x\otimes_{\OO(X)}\OO_y=
\M_x\otimes_{\OO_x}\OO_{xy} =\M_{xy}.$$

3. Consider the exact sequence of morphisms 
$$\M(V)\to \prod_{y\in V} \M_y\dosflechas \prod_{y,y'\in V, z\in U_{yy'}} \M_z.$$
Tensoring by $\otimes_{\OO(X)}\OO_x,$ we obtain the exact sequence of morphisms
$$\M(V)\otimes_{\OO(X)}\OO_x\to \prod_{y\in V} \M_{xy}\dosflechas \prod_{y,y'\in V, z\in U_{yy'}} \M_{xz}.$$
which shows that  $\mathcal M(V)\otimes_{\OO(X)} \OO_x=\M(V\cap U_x)$.

4. Consider the universally exact sequence $(*)$, where $X=U$.
Tensoring by  $\M(V)\otimes_{\OO(X)}$, we obtain the exact sequence of morphisms
$$\M(V)\otimes_{\OO(X)}\OO(U)\to \prod_{x\in U}\M(V\cap U_x)\dosflechas \prod_{x,y\in U, z\in U_{xy}} \M(V\cap U_z),$$
which shows that $\mathcal M(V)\otimes_{\OO(X)} \OO(U)=\M(V\cap U)$.

\edemo
 
\teor[(\cite{KS} 2.5,\,4.12)\,] Let $(X,\OO)$ be an affine finite space. Consider the canonical morphism $$\pi\colon (X,\OO)\to (*,\OO(X)), \,\,\pi(x)=*,\text{ for any }x\in X.$$
The functors

$$\xymatrix @R=8pt { \text{\bf Qc-Mod}_X \ar[r]^-{\pi_*} & { \text{\bf Mod}}_{\OO(X)}, & \mathcal M \ar@{|->}[r] & \pi_*\M=\M(X)
\\ {\text{\bf Mod}}_{\OO(X)}
 \ar[r]^-{\pi^*} & \text{\bf Qc-Mod}_X, &  M \ar@{|->}[r] & \pi^*\M=\tilde M}
$$
establish an equivalence between the category of  quasi-coherent  $\OO_X$-modules and the category of
$\OO(X)$-modules.
\eteor 

\demo 
The natural morphism $\pi^*\pi_*\M\to\M$ is an isomorphism because this morphism on stalks at $x$ is the morphism
$\M(X)\otimes_{\OO(X)}\OO_x\to \M_x$, which is an isomorphism by Proposition \ref{11}.

The natural morphism $M\to \pi_*\pi^*M=(\pi^*M)(X)$ is an isomorphism: Tensoring the exact sequence of morphisms $(*)$, in the proof of Proposition \ref{11},  by
$M\otimes_{\OO(X)} -$ we obtain the exact sequence of morphisms 
$$M\otimes_{\OO(X)}\OO(X)\to \prod_{x\in X}(\pi^*M)_x\dosflechas \prod_{x,y\in X, z\in U_{xy}} (\pi^*M)_z$$
which shows that $M=M\otimes_{\OO(X)}\OO(X)=(\pi^*M)(X)$.
\edemo
\lema \label{prodaf1} Let $A\to B$ and $A'\to B'$ be flat  (resp. failthfully flat) morphisms of commutative $C$-{a}lgebras. Then, $A\otimes_CA'\to B\otimes_C B'$ is a flat morphism (resp. faithfully flat).\elema

\demo  It follows from the equality $M\otimes_{A\otimes_CA'} (B\otimes_CB')=(M\otimes_AB)\otimes_{A'}B'$.\edemo

\prop \label{inter} The intersection of two affine open sets of an affine finite space is affine.\eprop

\demo Let $U$ and $U'$ be two affine open sets of the affine finite space $X$. Consider the faithfully flat morphisms  $\OO(U)\to \prod_{x\in U}\OO_x$, $\OO(U')\to \prod_{x'\in U'}\OO_{x'}$. The composition of the faithfully flat morphisms  (recall  Lemma \ref{prodaf1})

$$\OO(U'\cap U)\overset{\text{\ref{11}}}=\OO(U)\otimes_{\OO(X)}\OO(U')\to \prod_{(x,x')\in U\times U'}\OO_{xx'}\to \prod_{(x,x')\in U\times U', z\in U_{xx'}} \OO_z,$$ 
is faithfully flat, hence $\OO(U'\cap U)\to \prod_{z\in U\cap U'}\OO_z$ is faithfully flat. By Proposition \ref{afinfp}, $U\cap U'$ is affine.
%
%Si los $U'\cap U_x$ son afines, entonces
%$\OO(U'\cap U_x)\to \prod_{z\in U'\cap U_x} \OO_z$ es fielmente plano, luego  $\OO(U'\cap U)\to \prod_{z\in U\cap U'} \OO_z$ es fielmente plano y $U'\cap U$ es affine. Hemos reducido el problema a demostrar que $U'\cap U_x$ es affine, que se reduce a probar que $U_y\cap U_x$ es affine, y \'{a}stos lo son (por el corolario \ref{Uxy}).

\edemo

Let $R$ be a commutative ring with a unit. A \textit{finite $R$-ringed space} is a finite ringed space $(X,\OO)$ such that $\OO$ is a sheaf of $R$-algebras; that is, for any $x \in X$, $\OO_x$ is an $R$-algebra and for any $x \leq x'$, $r_{xx'} \colon \OO_x \to \OO_{x'}$ is a morphism of $R$-algebras.

Let $X$ and $Y$ be two finite $R$-ringed spaces. The \textit{direct product} $X \times_R Y$ is the finite $R$-ringed space $(X \times Y, \OO_{X \times Y})$, where $(\OO_{X \times Y})_{(x,y)}:= \OO_x \otimes_{R} \OO_y$, for each $(x,y) \in X \times Y$ and the morphisms of restriction are the obvious ones.

\prop[(\cite{KG} 5.27)\,] \label{ProdAf} Let $X$ and $Y$ be affine finite $R$-ringed spaces. Then, $X\times_R Y$ is an affine finite space and $\OO(X\times_R Y)=\OO(X)\otimes_R \OO(Y)$.\eprop

\demo Consider the universally exact sequence $\OO(X)\to \proda{x\in X} \OO_x\dosflechas \proda{x,x';z\in U_{xx'}} \OO_z$. Tensoring by $\otimes_{R}\OO_y$ we obtain the exact sequence
$$\OO(X)\otimes_{R}\OO_y\to \proda{x\in X} \OO(U_x\times_R U_y)\dosflechas \proda{x,x';z\in U_{xx'}} \OO(U_z\times_R U_y).$$
Hence, $\OO(X)\otimes_{R}\OO_y=\OO(X\times_R U_y)$. 
Consider the universally exact sequence $\OO(Y)\to \proda{y\in Y} \OO_y\dosflechas \proda{y,y';z\in U_{yy'}} \OO_z$. Tensoring by $\OO(X)\otimes_{R}$ we obtain the exact sequence
$$\OO(X)\otimes_{R}\OO(Y)\to \proda{y\in Y} \OO(X\times_R U_y)\dosflechas \proda{y,y';z\in U_{yy'}} \OO(X\times_R U_z).$$
Hence, $\OO(X)\otimes_{R}\OO(Y)=\OO(X\times_R Y)$.
In particular, $\OO_{xx'}\otimes_R\OO_{yy'}=\OO_{(x,y)(x',y')}$, for any $x,x'\in X$ and $y,y'\in Y$.
By Lemma \ref{prodaf1}, the morphism
$$\OO(X\times_R Y)=\OO(X)\otimes_R \OO(Y) \to \proda{x\in X}\OO_x\otimes_R \proda{y\in Y} \OO_y = \proda{(x,y)\in X\times_R Y} \OO_{(x,y)}$$ is faithfully flat. By Lemma \ref{prodaf1}, the morphism
$$\OO_{(x,y)(x',y')}=\OO_{xx'}\otimes_R \OO_{yy'} \to \proda{z\in U_{xx'}}\OO_z\otimes_R \proda{z'\in U_{yy'}} \OO_{z'} = \proda{(z,z')\in U_{xx'}\times_R U_{yy'}} \OO_{(z,z')}= \proda{(z,z')\in U_{(x,y)(x',y')}} \OO_{(z,z')}$$ is faithfully flat. Therefore, $X\times_R Y$ is affine.

\edemo

\subsection{Some conmutative algebra results.}

If $X$ is an affine finite space then, for each $x \leq y \in X$, 
the morphism $\OO_x\to \OO_y$ is flat and $\OO_y\otimes_{\OO_x}\OO_y=\OO_{yy}=\OO_y$. In this subsection we study this kind of morphisms. In this paper, we use well-known  properties of flat  morphisms and  faithfully flat morphisms,   that can be found in \cite{Matsumura}.

\nota\label{Notation} Given $\pp\in\Spec R$ and an $R$-module $M$, we denote
$M_\pp:=(R\backslash \pp)^{-1}\cdot M$.\enota

\prop \label{last} Let $f\colon A\to B$ be a morphism of rings and $f^*\colon \Spec B\to \Spec A$ the induced morphism. The following conditions are equivalent:

\enumera \item $A\to B$ is a flat morphism and $B\otimes_A B=B$.

\item $A_{f^*(\pp)}=B_{f^*(\pp)}$, for all $\pp\in\Spec B$.

\item The morphism $f^*\colon \Spec B\to \Spec A$ is injective and $A_{f^*(\pp)}=B_\pp$, for any $\pp\in\Spec B$.
\eenumera

\eprop 

\demo  $1. \Rightarrow \,2.$ The morphism $A_{f^*(\pp)}\to B_{f^*(\pp)}$ is faithful flat, for any $\pp$. Besides, $B_{f^*(\pp)}=A_{f^*(\pp)}\otimes_A B= A_{f^*(\pp)}\otimes_A(B\otimes_A B)=B_{f^*(\pp)}\otimes_{A_{f^*(\pp)}}B_{f^*(\pp)}$, then  $A_{f^*(\pp)}=B_{f^*(\pp)}$. 

$2. \Rightarrow \,3.$  If $A_{f^*(\pp)}=B_{f^*(\pp)}$ then $B_{f^*(\pp)}=B_\pp$ and $f^{*-1}({f^*(\pp)})=\{\pp\}$, then $f^*$ is injective. 

%$3. \Rightarrow \,2.$ Assume $A_{f^*(\pp)}=B_\pp$
%and that $\Spec B\to \Spec A$ is injective. Since the composition of the injective morphisms $\Spec B_\pp\iny \Spec B_{f^*(\pp)}\iny \Spec A_{f^*(\pp)}$ is the identity then $\Spec B_\pp=\Spec B_{f^*(\pp)}$, therefore $B_\pp=B_{f^*(\pp)}$.

$3. \Rightarrow \,1.$ The morphism $A\to B$ is flat: Given an injective morphism $N\hookrightarrow M$ of $A$-modules,
$N_{f^*(\pp)}\to M_{f^*(\pp)}$ is injective, for any $\pp$.  Then,
$N\otimes_A B_\pp\to M\otimes_A B_\pp$ is injective, for any $\pp$ and $N\otimes_A B\to M\otimes_A B$ is injective.
 
$\Spec B_{\pp}\subseteq \Spec B_{f^*(\pp)}\subseteq \Spec A_{f^*(\pp)}$ and $\Spec B_{\pp}=\Spec A_{f^*(\pp)}$.
Hence,  $\Spec B_{\pp}=\Spec B_{f^*(\pp)}$ and $B_{f^*(\pp)}=B_{\pp}$.
Then, 
 $(B\otimes_A B)_\pp=(B\otimes_A B)\otimes_B B_\pp =(B\otimes_AB)\otimes_B B_{f^*(\pp)}=
(B\otimes_AB)\otimes_A A_{f^*(\pp)}= 
 B_{f^*(\pp)} \otimes_{A_{f^*(\pp)}} B_{f^*(\pp)}=B_\pp$, for any $\pp\in\Spec B$. Therefore, $B\otimes_A B=B$.

\edemo

\nota \label{N3.9} Given a morphism $f\colon A\to B$ and an ideal $I\subseteq B$ denote $A\cap I:=f^{-1}(I)$. Denote $(I)_0=\{\pp\in\Spec B\colon I\subseteq \pp\}$.
\enota

\prop \label{sudor} Let  $A\to B$ be a flat morphism of rings such that $B\otimes_A B=B$. Then,
\enumera

\item  
$(I\cap A)\cdot B=I$, for any ideal $I\subseteq B$. 

\item $\Spec B$ is a topological subspace of $\Spec A$, with their Zariski topologies.

\item Let $\qq\in \Spec A$. \enumera \item  If $\qq\notin \Spec B$, then
$\qq\cdot B=B$. \item If $\qq\in\Spec B$, then $\qq\cdot B\subset B$ is a prime ideal and $(\qq\cdot B) \cap A=\qq$.\eenumera

\item $\Spec B=\capa{\Spec B \subseteq\, \text{open set  $U\subseteq \Spec A$}} U$.\eenumera
\eprop

\demo 1.  
Let $\pp\in \Spec B$,  $\qq:=A\cap \pp$ and  $M$ a $B$-module. By  Proposition \ref{last}, $M_\pp=M\otimes_BB_\pp=M\otimes_BB_\qq=M_\qq$. Then,
$$[(I\cap A)\cdot B]_\pp=[(I\cap A)\cdot B]_\qq=(I_\qq\cap A_\qq)\cdot B_\qq=
I_\qq=I_\pp.$$
Hence,  $(I\cap A)\cdot B=I$. 

2. By Proposition \ref{last}, we can think $\Spec B$ as a subset of $\Spec A$. Given an ideal $I\subseteq B$, observe that $(I)_0=((I\cap A)\cdot B)_0=(I\cap A)_0\cap \Spec B$.

3. (a) Suppose that there exists a prime ideal $\pp\subset B$ that contains to
$\qq\cdot B$. Denote $\pp'=\pp\cap A$.
Then, $\qq\in \Spec A_{\pp'}=\Spec B_\pp\subseteq \Spec B.$, which is contradictory. (b) Let $\pp\in \Spec B$ be a prime ideal such that $\pp\cap A=\qq$. Then,
$\pp=(\pp\cap A)\cdot B=\qq\cdot B$.

4. If  $\qq\in \Spec A\backslash \Spec B$, then $(\qq)_0\cap \Spec B=
(\qq\cdot B)_0=(B)_0=
\emptyset$. Then, $\Spec B$ is equal to the intersection of the open sets $U\subseteq \Spec A$ such that $\Spec B\subseteq U$.

\edemo

\section{Schematic finite spaces}

\subsection{Definition, examples and first characterizations}

\defi \label{defesq} We say that a finite ringed space $(X,\OO)$ is a \textit{schematic finite space} if it is locally affine; i.e. if  there exists an open covering $\{U_i\}_{i\in I}$ on $X$, such that $U_i$ is an affine finite space, for each $i \in I$.
\edefi

\prop Let $X$ be a finite ringed space. $X$ is a schematic finite space iff the open subsets $U_x$ are affine finite spaces for all $x \in X$.
\eprop
\demo $\Rightarrow)$ Let $\{U_i\}_{i\in I}$  be an affine open covering of $X$. For each $x \in X$, $U_x$ is an open subset of one of the affine finite spaces $U_i$. So, it follows from Corollary \ref{Uxy}, that $U_x$ is also affine.

$\Leftarrow)$ It is clear, since $\{U_x\}_{x \in X}$ is an open covering of $X$. 
\edemo

\obses
\item 1. All schematic finite spaces are finite fr-spaces.

\item 2. Affine finite spaces are schematic.

\item 3. If $X=U_x$, then $X$ schematic if and only if it is affine.
\eobses
\medskip

The finite ringed space  associated with a minimal affine  finite covering $\U$ of a quasi-compact and quasi-separated scheme is a schematic finite space, by Paragraph \ref{Parr}.

\ejems \label{Eejemplos} Let us give some examples of schematic finite spaces
(below we indicate the ringed space constructed  in  Paragraph \ref{Pa2}):

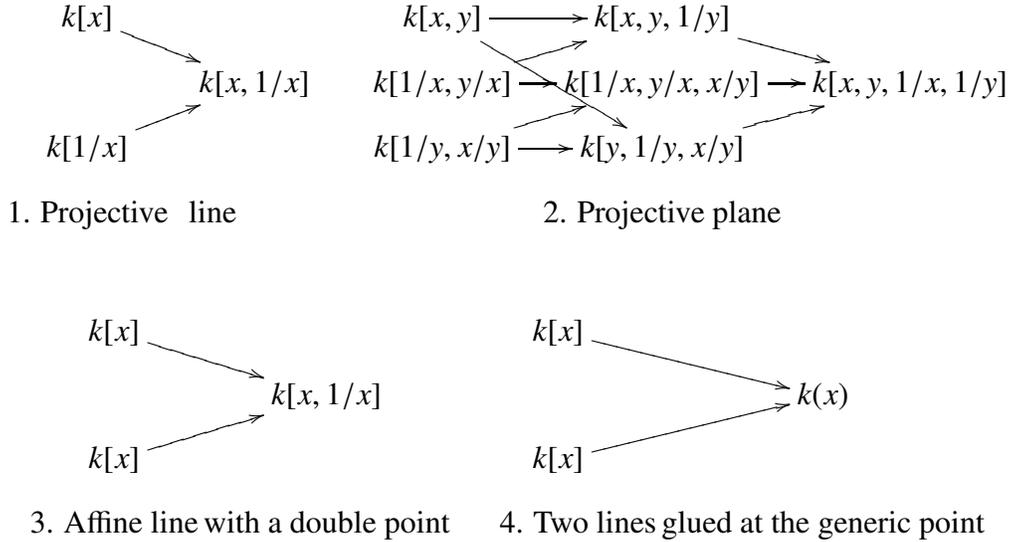
$$\xymatrix@C=6pt @R=8pt{ k[x] \ar[dr] & & & k[x,y] \ar[r]\ar[rdd]  &  k[x,y,1/y] \ar[rd] &   
\\  & k[x,1/x]  & & k[1/x,y/x] \ar[r] \ar[ru] & k[1/x,y/x,x/y] \ar[r] & k[x,y,1/x,1/y]
\\ k[1/x] \ar[ru] & & & k[1/y,x/y] \ar[r] \ar[ru] & k[y,1/y,x/y] \ar[ru] & 
\\\text{1. Projective} & \!\! \!\!\!\!\!\!\!\!\text{line} \quad& &  & \text{2. Projective plane} & }
$$

$$\xymatrix @C=-4pt @R=8pt{ k[x] \ar[dr] & & & k[x] \ar[dr]  & \\ & k[x,1/x]  & && k(x)\\
 k[x] \ar[ru] &&& k[x] \ar[ru]&\\ \text{3. Affine line} & \text{with a double point} & & \quad\text{ 4. Two lines} & \text{glued at the  generic point}}$$ 

It can be proved that the first three examples are finite models of the schemes we indicate, but the fourth it is not the model of any scheme. Also note that none of these examples are affine finite spaces.

%Grupos $H^i(X,\OO_X)$: 1. $H^0=k$, $H^1=0$. 2. $H^0=k[x]$, $H^1=k[1/x]/k$. 3. $H^0=k[x]$, $H^1=k(x)/k[x]$. 4. $H^0=k$, $H^i=0$, para $i>0$. 5. $H^0=k[x,y]$, $H^1=0$ y $H^2=\langle \frac{1}{x^ny^m}\rangle_{n,m>0}$.

\eejems

If $X$ and $Y$ are schematic finite $R$-ringed spaces, then $X\times_R Y$ is an schematic finite space, by Proposition \ref{ProdAf}.

\prop[(\cite{KS} 4.11)\,] \label{first characterization schematic} Let $X$ be a finite ringed space. $X$ is a schematic finite space if and only if it satisfies the following two conditions:
\enumera
 \item The natural morphism  $\OO_{y}\otimes_{\OO_x}\OO_{y'} \to \OO_{yy'}$ is an isomorphism for any $y,y' \geq x$.
 
\item The natural morphism $\OO_{yy'}\to \prod_{z\in U_{yy'}}\OO_z$  is faithfully flat, for any $y,y' \in X$ for which there is an element $x \in X$ such that $y \geq x$ and $y' \geq x$.  
\eenumera
\eprop
\demo It follows easily from the definition of affine finite space that the open subsets $\{U_x\}_{x \in X}$ are affine finite spaces iff the conditions 1. and 2. above are satisfied.
\edemo

\prop[(\cite{KS} 4.11)\,] \label{4pg20} Let $(X,\OO_X)$ be a ringed finite space. $X$ 
is a schematic finite space iff the morphism  $$\OO_{y}\otimes_{\OO_x}\OO_{y'}\to \prod_{z\in U_{yy'}} \OO_z$$ 
is faithfully flat, for any  $x\leq y,y'\in X$.\eprop

\demo  It follows directly from Corollary \ref{minimal affine}.
 \edemo

\subsection{More characterizations of schematic finite spaces}

In this section, we see that schematic finite spaces can be characterized by the good behavior of their quasi-coherent modules.

\prop \label{C4.9} Let $X$ be a schematic finite space, $U\overset{i}\subseteq X$ an open subset and $\Nc$ a quasi-coherent $\OO_U$-module. Then, $i_* \Nc$ is a quasi-coherent $\OO_X$-module. \eprop

\demo Let $x \leq y \in X$. We have to see that the morphism $(i_* \Nc)_x \otimes_{\OO_x} \OO_y \to (i_* \Nc)_y$ is an isomorphism.  This morphism is equal to the morphism $$ \Nc(U \cap U_x) \otimes_{\OO(U_x)} \OO(U_y) \to \Nc(U \cap U_y ),$$
which is an isomorphism by Proposition \ref{11}.
\edemo

\hskip-0.65cm\colorbox{white}{\,\begin{minipage}{15.15cm}\teor \label{T1} Let $X$ be a finite ringed space. $X$ is a schematic finite space if and only if it satisfies the next two conditions:
\enumera
\item $\Ker f$ is quasi-coherent, for any morphism $f \colon \M \to \Nc$ of quasi-coherent $\OO_{X}$-modules.

\item For any open subset $i\colon U_x\iny X$ and any quasi-coherent $\OO_{U_x}$-module $\M$, the $\OO_X$-module $i_*\M$ is quasi-coherent.
\eenumera

\eteor\end{minipage}}

\demo  $\Rightarrow)$ We know that schematic finite spaces are finite fr-spaces. By Proposition \ref{fr}, $\Ker f$ is quasi-coherent. The second condition follows from the proposition above.

$\Leftarrow)$ First, let us prove that $X$ is an fr-space.
Let $i\colon U_x\iny X$ be an open subset and $\M\to \Nc$  a morphism of quasi-coherent $\OO_{U_x}$-modules. 
$\Ker[\M\to \Nc]=\Ker[i_*\M\to i_*\Nc]_{|U_x}$, then it is quasi-coherent. By Proposition \ref{fr}, $X$ is an fr-space.

If $X$ is an fr-space and satisfies condition 2., then it is schematic:

Consider $x \leq x'$, let $j\colon U_{x'}\iny U_x$ be the inclusion morphism  and $\mathcal N$ a quasi-coherent $\OO_{U_{x'}}$-module. Since condition 2. is satisfied, the $\OO_{U_x}$-module $j_*\Nc=i^*((i\circ j)_*\Nc)$ is quasi-coherent. It follows from this result that we can suppose $X=U_x$ (because being finite fr-space and schematic are local conditions).

Now, by Corollary \ref{minimal affine}, we only have to prove that, for each $y,y'\geq x$,
the morphism $$\OO_y \otimes_{\OO_x}\OO_{y'} \to \prod_{z\in U_{yy'}} \OO_{z}$$ is faithfully flat.

Consider the open subset
$i\colon U_y\iny X=U_x$. Since $i_*\OO_{U_y}$ is quasi-coherent,
$$\OO_{yy'}=(i_*\OO_{U_y})(U_{y'})=(i_*\OO_{U_y})(U_{x})\otimes_{\OO_x}\OO_{y'}=\OO_y\otimes_{\OO_x}\OO_{y'}.$$
In particular, $\OO_y=\OO_{yy}=\OO_y\otimes_{\OO_x}\OO_{y}$. The morphism $\OO_y\to \OO_z$ is flat, for any $z\geq y,y'$, then the morphism
$\OO_y\otimes_{\OO_x}\OO_{y'}\to \OO_z\otimes_{\OO_x}\OO_{y'}=\OO_{zy'}=\OO_z$ is flat.

If the morphism $\OO_y\otimes_{\OO_x}\OO_{y'} \to \prod_{z\in U_{yy'}} \OO_{z}$ is not faithfully flat,
there exists an ideal $I\underset\neq\subset \OO_y\otimes_{\OO_x}\OO_{y'}$ such that $I\cdot \prod_{z\in U_{yy'}} \OO_{z}=\prod_{z\in U_{yy'}} \OO_{z}$. Observe that the morphism
$\OO_y\to \OO_y\otimes_{\OO_x}\OO_{y'} $ is flat since $\OO_x\to \OO_{y'}$ is  flat. Besides,
$$(\OO_y\otimes_{\OO_x}\OO_{y'})\otimes_{\OO_y}(\OO_y\otimes_{\OO_x}\OO_{y'})=\OO_{y}\otimes_{\OO_x}(\OO_{y'}\otimes_{\OO_x}\OO_{y'})=\OO_y\otimes_{\OO_x}\OO_{y'}. $$
By Proposition \ref{sudor}, there exists an ideal $J\subset \OO_y$ such that
$J\cdot (\OO_y\otimes_{\OO_x}\OO_{y'})=I$. Let $\mathcal M$ be the quasi-coherent $\OO_{U_y}$-module
 associated with  the $\OO_y$-module $\OO_y/J$. Then, $i_*\mathcal M$ is the quasi-coherent $\OO_X$-module
 associated with the $\OO_x$-module $\OO_y/J$ and 
$$\M(U_{yy'})=(i_*\M)(U_{y'}) =(\OO_y/J)\otimes_{\OO_x}\OO_{y'}=(\OO_y\otimes_{\OO_x}\OO_{y'})/J\cdot (\OO_y\otimes_{\OO_x}\OO_{y'})=(\OO_y\otimes_{\OO_x}\OO_{y'})/I\neq 0.$$
However, ${\mathcal M}_{|U_{yy'}}=0$, since $\mathcal M_z=(\OO_y/J)\otimes_{\OO_y}\OO_z=\OO_z/J\cdot \OO_z=\OO_z/I\cdot \OO_z=0$, for any $z\in U_{yy'}$.
So we have a contradiction; therefore, 
  the morphism $\OO_{y}\otimes_{\OO_x}\OO_{y'} \to \prod_{z\in U_{yy'}} \OO_{z}$ is faithfully flat.

% $(i_*\mathcal M)_{|U_{y'}}$  
%is the quasi-coherent $\OO_{U_{y'}}$-module
% associated with the $\OO_{y'}$-module 
% $$(\OO_y/J)\otimes_{\OO_x}\OO_{y'}=(\OO_y\otimes_{\OO_x}\OO_{y'})/J\cdot (\OO_y\otimes_{\OO_x}\OO_{y'})=(\OO_y\otimes_{\OO_x}\OO_{y'})/I\neq 0.$$
%${\mathcal M}_{|U_{yy'}}=0$, since $\mathcal M_z=(\OO_y/J)\otimes_{\OO_y}\OO_z=\OO_z/J\cdot \OO_z=\OO_z/I\cdot \OO_z=0$, for any $z\in U_{yy'}$. Finally,
%$$0=i'_*({\mathcal M}_{|U_{yy'}})=(i_*\mathcal M)_{|U_{y'}}\neq 0,$$
%where $i'\colon U_{yy'}\to U_{y'}$ is defined by $i'(z)=z$. This is contradictory, then 
%  the morphism $\OO_{y}\otimes_{\OO_x}\OO_{y'} \to \prod_{z\in U_{yy'}} \OO_{z}$ is faithfully flat.
% 
\edemo

\hskip-0.65cm\colorbox{white}{\,\begin{minipage}{15.15cm}
\teor \label{Cdelta} Let $(X,\OO)$ be a finite ringed space.  Let $\delta\colon X\to X\times X$, $\delta(x)=(x,x)$ be the diagonal morphism.
Then, $X$ is schematic iff it satisfies these two conditions:
\enumera
\item $\Ker f$ is quasi-coherent, for any morphism $f \colon \M \to \Nc$ of quasi-coherent $\OO_{X}$-modules.

\item $\delta_{*}\Nc$ is a quasi-coherent $\OO_{X\times X}$-module for any quasi-coherent $\OO_X$-module $\Nc$.
\eenumera 
\eteor\end{minipage}}

\demo $\Rightarrow)$ For any $(x,y)\leq (x',y')$, we have
$$\aligned (\delta_*\Nc)_{(x,y)}\otimes_{\OO_{(x,y)}}\OO_{(x',y')} & = \Nc_{xy} \otimes_{\OO_{x} \otimes \OO_y}(\OO_{x'} \otimes \OO_{y'}) =
\Nc_{xy}\otimes_{\OO_x}\OO_{x'}\otimes_{\OO_{y}}\OO_{y'}
\overset{\text{\ref{11}}}=\Nc_{x'y}\otimes_{\OO_{y}}\OO_{y'}
\overset{\text{\ref{11}}}=\Nc_{x'y'}\\ & =(\delta_*\Nc)_{(x',y')}.\endaligned$$

$\Leftarrow)$ First, note that for any $x \in X$ and any $x'\leq x''$,
$$\OO_{xx'}\otimes_{\OO_{x'}}\OO_{x''}= (\delta_* \OO)_{(x,x')} \otimes_{\OO_x} \OO_x \otimes_{\OO_{x'}} \OO_{x''}= (\delta_* \OO)_{(x,x')} \otimes_{\OO_{(x,x')}} \OO_{(x,x'')}= (\delta_* \OO)_{(x,x'')}= \OO_{xx''}.$$
In consequence, for any open subset $i\colon U_x\iny X$ and any quasi-coherent $\OO_{U_x}$-module $\M$ there exists a quasi-coherent $\OO_X$-module $\Nc$ such that $\Nc_{|U_x}\simeq \M$: define $\Nc_{x'}:=\M_x\otimes_{\OO_x}\OO_{xx'}$, for any $x'\in X$. $\Nc$ is quasi-coherent since for any $x'\leq x''$,
$$\Nc_{x'}\otimes_{\OO_{x'}}\OO_{x''}=\M_x\otimes_{\OO_x}\OO_{xx'}\otimes_{\OO_{x'}}\OO_{x''}=\M_x\otimes_{\OO_x}\OO_{xx''}=\Nc_{x''}.$$
Besides, $\Nc_x=\M_x$, so $\Nc_{|U_x}=\M$.

By Theorem \ref{T1},  
we have to prove that $i_*\M$ is a quasi-coherent $\OO_X$-module.
That is,  we have to prove that
$(i_*\M)_{y'}=(i_*\M)_y \otimes_{\OO_y}\OO_{y'}$, for any
$y\leq y'\in X$:

$$\aligned (i_*\M)_{y'}=\M(U_{y'x})=\Nc(U_{y'x}) & =(\delta_*\Nc)(U_{y'}\times U_{x})=
(\delta_*\Nc)(U_y\times U_{x})\otimes_{\OO_y\otimes \OO_x}\OO_{y'}\otimes \OO_{x}\\ & =\Nc(U_{yx})\otimes_{\OO_y}\OO_{y'}=\M(U_{yx}) \otimes_{\OO_y}\OO_{y'}=(i_*\M)_y \otimes_{\OO_y}\OO_{y'}.\endaligned$$

\edemo 

\obse  The theorem above can be restated by saying that a finite ringed space $X$ is schematic iff it is a finite fr-space and  for any quasi-coherent module $\Nc$, any $x \in X$ and any $x'\leq x'' \in X$, $\Nc_{xx'}\otimes_{\OO_{x'}}\OO_{x''}=\Nc_{xx''}$.
\eobse

\prop \label{K2?} Let $X$ be a finite fr-space and $\Nc$ a quasi-coherent $\OO_X$-module. Then, 
$\Nc_{pq}\otimes_{\OO_p}\OO_{p'}=\Nc_{p'q}$, for any $p\leq p' \in X$ and for any  $q \in X$ iff
$\Nc_{p'}\otimes_{\OO_p}\OO_{p''}=\Nc_{p'p''}$, for any $p\leq p' ,p'' \in X$.

\eprop

\demo $\Rightarrow)$ $\Nc_{p'}\otimes_{\OO_p}\OO_{p''}=\Nc_{pp'}\otimes_{\OO_p}\OO_{p''}=\Nc_{p'p''}$.

$\Leftarrow)$ Let $U\subseteq U_p$ be an open set  and $p\leq p'$. Consider the exact sequence of morphisms
$$\Nc(U)\to \prod_{x\in U} \Nc_x\dosflechas \prod_{z\geq x\in U} \Nc_z.$$
Tensoring by $\otimes_{\OO_p} \OO_{p'}$ we obtain the exact sequence of morphisms
$$\Nc(U)\otimes_{\OO_p} \OO_{p'}\to \prod_{x\in U} \Nc_{p'x}\dosflechas \prod_{z\geq x\in U} \Nc_{p'z},$$
which shows that $\Nc(U)\otimes_{\OO_p} \OO_{p'}=\Nc(U\cap U_{p'})$. In particular,
 $$\Nc_{pq}\otimes_{\OO_p}\OO_{p'}=\Nc_{p'q}.$$\edemo

\coro \label{ultima caracterizacion schem} Let $X$ be a finite ringed space. $X$ is schematic iff  it is a finite fr-space and for any quasi-coherent $\OO_X$-module $\Nc$ and any $x\leq x' ,x'' \in X$, $\Nc_{x'}\otimes_{\OO_x}\OO_{x''}=\Nc_{x'x''}$
\ecoro
\demo It follows directly from Theorem \ref{Cdelta} and Proposition above.
\edemo

\section{Affine morphisms}

\defi  Let $X$ and $Y$ be schematic finite spaces.  A morphism $f\colon X\to Y$ of ringed spaces  is said to be an affine morphism if  $f_*\OO_X$ is a quasi-coherent $\OO_Y$-module and the preimage of any affine open subspace of $Y$ is an affine open subspace of $X$.\edefi

\ejems \label{Kolmogorov} 1. A schematic finite space $X$ is affine iff  $(X,\OO)\to (*,\OO(X))$ is an affine morphism.

2. If $X$ is an affine finite space and   $U\subseteq X$ an affine open subset, the  
inclusion morphism $i\colon U\iny X$ is an affine morphism: $i_*\OO_U$ is quasi-coherent by Proposition \ref{11}, and for any affine open   subset  $V\subseteq Y$, $i^{-1}(V)=V\cap U$ is affine by Proposition \ref{inter}.

3. Let $X$ be a schematic finite space. Given $x.x'\in X$, we shall say that 
$x\sim x'$ if $x\leq x'$ and $x'\leq x$. Let $\bar X:=X/\sim$ be the Kolmogorov quotient of $X$
and define $\OO_{[x]}:=\OO_x$, for any $[x]\in\bar X$. Then, $\bar X$ is a schematic space, the quotient morphism $\pi\colon \bar X\to X$, $\pi(x):=[x]$ is 
affine and  $\pi_*\OO_X=\OO_{\bar X}$.

\eejems

%\ejem \label{xX} Let $X$ be an affine finite space. Let $X':=\{x\}\coprod X$ be the ringed finite space defined as follows: $x<x'$ for any $x'\in X$, and the orden relation on $X$ is the pre-established order relation; $\OO_{X',x}:=\OO_{X}(X)$ and $\OO_{X',x'}:=\OO_{X,x'}$ and the morphisms between the stalks are the obvious morphisms. It is easy to check that $X$ is an affine finite morphism and the inclusion $X\iny X'$ morphism is an affine morphism, by the previous example.\eejem

\prop Let $X$ and  $Y$ be affine finite spaces and $f\colon X\to Y$ an affine morphism. Let
$M$ be an $\OO(X)$-module (therefore, an $\OO(Y)$-module). Then,
$$f_*\tilde M=\tilde M.$$

\eprop

\demo For any open set $U_y$, $$\aligned (f_*\tilde M)(U_y) & =
\tilde M(f^{-1}(U_y))\overset{\text{\ref{11}}}=
\tilde M(X)\otimes_{\OO(X)} \OO_X(f^{-1}(U_y)) =
M\otimes_{\OO(X)} \OO_X(f^{-1}(U_y))\\ & =
M\otimes_{\OO(X)} (f_*\OO_X)(U_y)\overset{\text{\ref{11}}}  =
M\otimes_{\OO(X)}  {\OO(X)} \otimes_{\OO(Y)}  \OO_Y(U_y)=M\otimes_{\OO(Y)}  \OO_Y(U_y)\\ & =\tilde M(Y)\otimes_{\OO(Y)}  \OO_Y(U_y)
\overset{\text{\ref{11}}}  =\tilde M(U_y).\endaligned$$
\edemo

\prop \label{qctoqc} Let $f\colon X\to Y$ be an affine morphism and  $\mathcal M$ a quasi-coherent $\OO_X$-module. Then,
$f_*\M$ is a quasi-coherent $\OO_Y$-module.\eprop

\demo   Being  $f_*\M$ a quasi-coherent $\OO_Y$-module is a local property. We can suppose that $Y$ is affine. Then $X$ is affine. The proof is completed by the previous proposition.
\edemo
%\demo Sea $y\leq y'\in Y$. Entonces,
%$$\aligned & (f_*\M)(U_y') =\M(f^{-1}(U_y'))\overset{\text{\ref{12}}}=\M(f^{-1}(U_y))\otimes_{\OO(f^{-1}(U_y))} \OO(f^{-1}(U_y')) \\ &=
%\M(f^{-1}(U_y))\otimes_{\OO(f^{-1}(U_y))}\OO(f^{-1}(U_y))\otimes_{\OO(U_y)} \OO(U_y'))=
%M(f^{-1}(U_y))\otimes_{\OO(U_y)} \OO(U_y'))\endaligned$$
%
%\edemo 

\prop  \label{PObv} The composition of affine morphisms is affine\eprop

\demo It is obvious.\edemo

\prop \label{hmm} Let $X$ and $Y$ be schematic finite spaces.  A morphism of ringed spaces $f\colon X\to Y$ is affine iff  $f_*\OO_X$ is quasi-coherent and 
$f^{-1}(U_y)$ is affine for any $y\in Y$.\eprop

\demo $\Rightarrow)$ It is obvious.

$\Leftarrow)$ Let us proceed by induction on $\# Y$. 
If $\# Y=1$, it is obvious. We can suppose that $Y$ is affine and we only have to prove that $X$ is affine.
The morphism $\OO(Y)\to \prod_{y\in Y}\OO_{y}$ is faithfully flat, then the morphism 

$$\OO(X)\to \prod_{y\in Y} \OO(X)\otimes_{\OO(Y)} \OO_{y}=\prod_{y\in Y} (f_*\OO_X)(Y)\otimes_{\OO(Y)} \OO_{y}
\overset{\text{\ref{11}}}=\prod_{y\in Y}  \OO(f^{-1}(U_y))$$
is faithfully flat. Since $f^{-1}(U_{y})$ is affine, the morphism $\OO(f^{-1}(U_y))\to \prod_{x\in f^{-1}(U_y)} \OO_x$ is faithfully flat. The composition of  faithfully flat morphisms is faithfully flat, then $\OO(X)\to \prod_{y\in Y, x\in f^{-1}(U_y)}  \OO_x$ is faithfully flat. Therefore, the morphism 
$$\OO(X)\to \prod_{x\in X}  \OO_x$$
is faithfully flat. 

Let $x,x'\in X$. Given an open set $V\subseteq Y$ denote $\bar V:=f^{-1}(V)$. 
$\bar U_{f(x)f(x')}$ is an affine open subset of $\bar U_{f(x)}$ (by induction hypothesis), then 
$\bar U_{f(x)f(x')}\cap U_x$ is  affine, and it is included in $\bar U_{f(x')}$, hence 
$\bar U_{f(x)f(x')}\cap U_x\cap U_{x'}$ is affine. 
Then,
$$U_{xx'}=\bar U_{f(x)f(x')}\cap U_x\cap U_{x'}$$
is affine and the morphism $\OO_{xx'}\to \prod_{z\in U_{xx'}} \OO_z$ is faithfully flat.
Since $f_*\OO_X$ is quasi-coherent (and Proposition  \ref{11}), 
$$\aligned \OO(\bar U_{f(x)}) & =f_*\OO_X(U_{f(x)})  =\OO(X)\otimes_{\OO(Y)}\OO_{f(x)}, \quad \OO(\bar U_{f(x')})=\OO(X)\otimes_{\OO(Y)}\OO_{f(x')}. \\
 \OO(\bar U_{f(x)f(x')}) & =f_*\OO_X(U_{f(x)f(x')}) =\OO(X)\otimes_{\OO(Y)}\OO_{f(x)f(x')}=\OO(X)\otimes_{\OO(Y)} \OO_{f(x)}\otimes_{\OO(Y)}\OO_{f(x')}\\ & =(\OO(X)\otimes_{\OO(Y)} \OO_{f(x)})\otimes_{\OO(X)}(\OO(X)\otimes_{\OO(Y)}\OO_{f(x')})=\OO(\bar U_{f(x)})\otimes_{\OO(X)}\OO(\bar U_{f(x')}.)\endaligned \qquad (*)$$
Now it is easy to prove that $$\OO_{xx'}=\OO(\bar U_{f(x)f(x')}\cap U_x\cap U_{x'})\overset{\text{\ref{11}}}=(\OO(\bar U_{f(x)f(x')})\otimes_{\OO(\bar U_{f(x)})} \OO_x)\otimes_{\OO(\bar U_{f(x')})} \OO_{x'}\overset{(*)}=\OO_x\otimes_{\OO(X)}\OO_{x'}.$$
Therefore, $X$ is affine.
\edemo 

\coro \label{afloc} Let $X$ and $Y$ be schematic finite spaces and let $f\colon X\to Y$ be a morphism of ringed spaces. Then,  being $f$ affine  is a local property on $Y$. 
\ecoro 

\ejem \label{E5.8} Let $(X,\OO)$ be a  schematic finite space and $\OO\to \OO'$
 a morphism of sheaves of rings, such that $\OO'$  is a quasi-coherent $\OO$-module. $(X,\OO')$ is a schematic finite space:
Given $x\leq y,y'$, the morphism
$\OO_y\otimes_{\OO_x} \OO_{y'}\to \prod_{z\in U_{yy'}} \OO_z$ is faithfully flat, by Proposition \ref{4pg20}. Tensoring by $\otimes_{\OO_x}\OO'_{x}$ we obtain the faithfully flat morphism
$\OO'_y\otimes_{\OO'_x} \OO'_{y'}\to \prod_{z\in U_{yy'}} \OO'_z$. Hence, $(X,\OO')$ is a schematic finite space by Proposition
\ref{4pg20}. The obvious morphism $\Id\colon (X,\OO')\to (X,\OO)$ is affine. \eejem

\section{Schematic morphisms}

\defi Let $X,Y$ be schematic finite spaces.  A morphism of ringed spaces  $f\colon X\to Y$ 
is said to be a schematic morphism if for any $x\in X$ the morphism $f_x\colon U_x\to U_{f(x)}$, $f_x(x'):=f(x')$ is affine.\edefi

\ejem If $X$ is a schematic finite space, $X\to (*,\OO(X))$ is a schematic morphism.\eejem

\ejem If $U$ is an open subspace of a schematic finite space $X$, then the inclusion morphism $U\iny X$ is schematic.\eejem

\obse Let $X$ and $Y$ be schematic finite spaces and $f\colon X\to Y$ be a morphism of ringed spaces.
Then,  being $f$ schematic  is a local property on $Y$ and on $X$.
\eobse

\prop The composition of schematic morphisms is schematic.\eprop

\demo It is a consequence of Proposition \ref{PObv}.\edemo

\prop \label{afinesqu} Affine morphisms between schematic finite spaces are schematic morphisms.
\eprop

\demo Let $f\colon X\to Y$ be an affine morphism. 
Then, $f^{-1}(U_{f(x)})$ is affine,  $U_x\iny f^{-1}(U_{f(x)})$ is an affine morphism and $f^{-1}(U_{f(x)})\to U_{f(x)}$ is affine, by Corollary \ref{afloc}. The composition $U_x\iny f^{-1}(U_{f(x)})\to 
U_{f(x)}$  is affine, by Proposition \ref{PObv}.
Hence, $f$ is a schematic morphism.

\edemo

\prop \label{K40} Let $f\colon X\to Y$ be a schematic morphism  and $\M$  a quasi-coherent $\OO_X$-module. 
Then,  $f_*\M$ is a quasi-coherent  $\OO_Y$-module. \eprop

\demo We can suppose that $Y$ is affine. Consider an open set $U_x\overset i\subseteq X$ and denote $\M_{U_x}=i_*{\M}_{|U_x}$. Observe that $f_*\M_{U_x}=(f\circ i)_*{\M}_{|U_x}$ is a  quasi-coherent  $\OO_Y$-module,  because the composite morphism  $f\circ i\colon U_x\to U_{f(x)}\iny Y$ is  affine and by  Proposition \ref{qctoqc}.

Let  $\{U_{x_1},\ldots, U_{x_n}\}$ be an open covering of $X$ and  $\{U_{x_{ijk}}\}_k$ an open covering of  $U_{x_i}\cap U_{x_j}$, for each $i,j$. Consider the exact sequence of morphisms 
$$\M\to \prod_i \M_{U_{x_i}} \dosflechas \prod_{i,j,k} \M_{U_{x_{ijk}}}$$
Taking $f_*$, we obtain an exact sequence of morphisms, then $f_*\M$ is a quasi-coherent $\OO_Y$-module.

\edemo

\coro \label{C4.10} Let $X$ be a schematic finite space and $U \overset i \iny X$ an open subset. Given a quasi-coherent $\OO_U$-module $\Nc$, there exists a quasi-coherent $\OO_X$-module $\M$, such that $\M_{|U}\simeq \Nc$.
\ecoro

\demo Define $\M:=i_*\Nc$.

\edemo

\lema Let $X$ be an  affine finite space and  $U\subset X$ an open set. Then, 
$U$ is affine iff $U\cap U_x$ is affine, for any $x\in X$.\elema

\demo If $U$ is affine, then $U\cap U_x$ is affine, for any  $x\in X$, by Proposition \ref{inter}.
Let us prove the converse implication. The inclusion morphism  $i\colon U\iny X$ is an affine morphism, by Proposition \ref{hmm}. Hence, $U=i^{-1}(X)$ is affine.
\edemo

\prop A morphism of ringed spaces  $f\colon X\to Y$ between affine finite spaces is affine iff it is a schematic morphism.\eprop

\demo  $\Rightarrow)$ It is known (see Proposition \ref{afinesqu}).

$\Leftarrow)$ By Proposition \ref{K40}, we only have to prove that  $f^{-1}(U)$ is affine, for any affine open subset  $U\subseteq Y$.
By the previous lemma, we  only have to prove that $f^{-1}(U)\cap U_x$ is affine. The composition of affine morphisms is affine, then 
$U_x\to U_{f(x)}\iny Y$ is affine. Hence, $f^{-1}(U)\cap U_x$ is affine.

\edemo 

\coro \label{C9.10} Let $f\colon X\to Y$ be a schematic morphism. Then, $f$ is affine iff there exists an affine open covering of $Y$, $\{U_i\}$, such that $f^{-1}(U_i)$ is affine, for any $i$.\ecoro

\demo Recall that  being  $f$ affine is a local property on $Y$.

\edemo

In \cite{KG} 5.6, it is proved that a morphism of ringed spaces $f\colon X\to Y$ is schematic iff $R^if_*\M$ is quasi-coherent for any quasi-coherent module $\M$ and any $i$.

\hskip-0.65cm\colorbox{white}{\,\begin{minipage}{15.15cm}
\teor \label{Tcohsch} Let $X$ and $Y$ be schematic finite spaces and $f\colon X\to Y$ a morphism of ringed spaces. Then, $f$ is a schematic morphism iff  $f_*\mathcal M$ is quasi-coherent, for any quasi-coherent $\OO_X$-module $\mathcal M$. \eteor\end{minipage}}

\demo  $\Rightarrow)$ Recall Proposition \ref{K40}.

$\Leftarrow)$ We can suppose that $X$ and $Y$ are affine. We only have to prove that 
$f^{-1}(U_y)$ is affine, for any $y\in Y$, by Proposition \ref{hmm}. By Proposition \ref{afinfp}, we only have to prove that
the morphism $\OO_X(f^{-1}(U_y))\to \prod_{x\in f^{-1}(U_y)} \OO_{x}$ is faithfully flat.
Observe that
$$\OO_X(f^{-1}(U_y))=(f_*\OO_X)(U_y)=(f_*\OO_X)(Y)\otimes_{\OO_Y(Y)}\OO_y=\OO_X(X)\otimes_{\OO_Y(Y)}\OO_y$$
and  for any $x\in f^{-1}(U_y)$
$$\OO_{x}\otimes_{\OO_Y(Y)} \OO_y=\OO_{x}\otimes_{\OO_y} \OO_y\otimes_{\OO_Y(Y)} \OO_y=\OO_{x}\otimes_{\OO_y} \OO_y=\OO_x.$$
The morphism $\OO_X(X)\to \OO_x$ is flat, then tensoring by $\otimes_{\OO_Y(Y)} \OO_y$
the morphism $\OO_X(f^{-1}(U_y))\to \OO_x$ is flat, for any $x\in f^{-1}(U_y)$.

If the morphism $\OO_X(f^{-1}(U_y))\to \prod_{x\in f^{-1}(U_y)} \OO_{x}$ is not faithfully flat,
there exists an ideal $I\underset\neq\subset \OO_X(f^{-1}(U_y))$ such that $I\cdot \prod_{x\in f^{-1}(U_y)} \OO_{x}=\prod_{x\in f^{-1}(U_y)} \OO_{x}$. Observe that the morphism
$\OO_X(X)\to \OO_X(X)\otimes_{\OO_Y(Y)}\OO_y=\OO_X(f^{-1}(U_y))$ is flat since $\OO_Y(Y)\to \OO_y$ is flat. Besides,
$$\aligned \OO_X(f^{-1}(U_y)) & \otimes_{\OO_X(X)} \OO_X(f^{-1}(U_y))  =
\OO_X(X)\otimes_{\OO_Y(Y)}\OO_y \otimes_{\OO_X(X)}   \OO_X(X)\otimes_{\OO_Y(Y)}\OO_y\\ & =\OO_X(X)\otimes_{\OO_Y(Y)} \OO_y\otimes_{\OO_Y(Y)}\OO_y=\OO_X(X)\otimes_{\OO_Y(Y)} \OO_y=\OO_X(f^{-1}(U_y)).\endaligned$$
By Proposition \ref{sudor}, there exists an ideal $J\subset \OO_X(X)$ such that
$J\cdot \OO_X(f^{-1}(U_y))=I$. Let $\mathcal M$ be the quasi-coherent $\OO_X$-module
 associated with  the $\OO_X(X)$-module $\OO_X(X)/J$. Then, $f_*\mathcal M$ is the quasi-coherent $\OO_Y$-module
 associated with the $\OO_Y(Y)$-module $\OO_X(X)/J$ and 
 
$$\aligned \mathcal M(f^{-1}(U_y)) & =f_*\M(U_y)=f_*\M(Y)\otimes_{\OO_Y(Y)}\OO_y=(\OO_X(X)/J)\otimes_{\OO_Y(Y)}\OO_y\\ & =(\OO_X(X)\otimes_{\OO_Y(Y)}\OO_y)/J\cdot (\OO_X(X)\otimes_{\OO_Y(Y)}\OO_y)=\OO_X(f^{-1}(U_y))/I\neq 0.\endaligned$$
However, ${\mathcal M}_{|f^{-1}(U_y)}=0$ since $\mathcal M_x=\OO_x/J\cdot \OO_x=\OO_x/I\cdot \OO_x=0$, for any $x\in f^{-1}(U_y)$. This is contradictory, then 
  the morphism $\OO_X(f^{-1}(U_y))\to \prod_{x\in f^{-1}(U_y)} \OO_{x}$ is faithfully flat.

% 
% 
% 
% $(f_*\mathcal M)_{|U_y}$  
%is the quasi-coherent $\OO_{U_y}$-module
% associated with the $\OO_y$-module 
% $$(\OO_X(X)/J)\otimes_{\OO_Y(Y)}\OO_y=(\OO_X(X)\otimes_{\OO_Y(Y)}\OO_y)/J\cdot (\OO_X(X)\otimes_{\OO_Y(Y)}\OO_y)=\OO_X(f^{-1}(U_y)/I\neq 0.$$
%${\mathcal M}_{|f^{-1}(U_y)}=0$, since $\mathcal M_x=\OO_x/J\cdot \OO_x=\OO_x/I\cdot \OO_x=0$, for any $x\in f^{-1}(U_y)$. Finally,
%$$0=f'_*({\mathcal M}_{|f^{-1}(U_y)})=(f_*\mathcal M)_{|U_y}\neq 0$$
%where $f'\colon f^{-1}(U_y)\to U_y$ is defined by $f'(x)=f(x)$. This is contradictory, then 
%  the morphism $\OO_X(f^{-1}(U_y))\to \prod_{x\in f^{-1}(U_y)} \OO_{x}$ is faithfully flat
% 
\edemo

\nota Let $f\colon X\to Y$  be a morphism of ringed spaces, between ringed finite spaces,  $x\in X$ and  $y\in Y$. We shall denote
$U_{xy}:=U_x\cap f^{-1}(U_y)$ and $\OO_{xy}:=\OO(U_x\cap f^{-1}(U_y))$.\enota

\hskip-0.65cm\colorbox{white}{\,\begin{minipage}{15.15cm}
\prop  \label{K6} A morphism of ringed spaces $f\colon X\to Y$ between schematic finite spaces is schematic iff for any  $x\in X$ and  $y\geq f(x)$ 
\enumera \item $U_{xy}$ is affine.

\item $\OO_{xy}=\OO_{x}\otimes_{\OO_{f(x)}}\OO_{y}$.\eenumera
\eprop\end{minipage}}

\demo Consider the morphism $f_x\colon U_x\to U_{f(x)}$. Then,
$f_{x*}\OO_{U_x}$ is a quasi-coherent module iff condition 2. is satisfied. By Proposition \ref{hmm}, $f_x$ is affine iff the conditions 1. and 2. are satisfied. Then, $f$ is schematic iff the conditions 1. and 2. are satisfied. 
\edemo

\hskip-0.65cm\colorbox{white}{\,\begin{minipage}{15.15cm}
\teor \label{K7} Let $f\colon X\to Y$ be a morphism of ringed spaces between schematic finite spaces.
Then, $f$ is schematic iff the  induced morphism on spectra by  the morphism of rings $\OO_x\otimes_{\OO_{f(x)}} \OO_y\to \prod_{z\in U_{xy}} \OO_{z}$ is surjective, for any $x$ and $y\geq f(x)$.\eteor\end{minipage}}

\demo $\Rightarrow)$  By Proposition \ref{K6},  $U_{xy}$ is affine
and 
$\OO_x\otimes_{\OO_{f(x)}} \OO_y=
\OO_{xy}$.  Therefore,  the morphism
$$\OO_x\otimes_{\OO_{f(x)}} \OO_y=\OO_{xy}\to \prod_{z\in U_{xy}} \OO_{z}$$
is faithfully flat, and the induced morphism on spectra is surjective.

$\Leftarrow)$ Let $z\in U_{xy}$. Since $\OO_x\to \OO_z$ is a flat morphism, the morphism
$$\OO_x\otimes_{\OO_{f(x)}}\OO_y\to \OO_z\otimes_{\OO_{f(x)}}\OO_y=\OO_z\otimes_{\OO_{f(z)}}\OO_{f(z)}\otimes_{\OO_{f(x)}}\OO_y=\OO_z\otimes_{\OO_{f(z)}} \OO_{f(z)}=\OO_z$$
is flat. Denote $B=\OO_x\otimes_{\OO_{f(x)}} \OO_y$ and $C=\prod_{z\in U_{xy}} \OO_{z}$. 
The morphism $B\to C$ is faithfully flat. 

Let $z,z'\in U_{xy}$.  The morphism $\OO_z\otimes_{\OO_x}\OO_{z'}\to \OO_z\otimes_{B}\OO_{z'}$ is surjective and the composite morphism $$\OO_z\otimes_{\OO_x}\OO_{z'}\to \OO_z\otimes_{B}\OO_{z'}\to \OO_{zz'}$$ is an isomorphism. Therefore, $\OO_z\otimes_{B}\OO_{z'}= \OO_{zz'}$. The exact sequence of morphisms 
$$B\to C =\prod_{z\in U_{xy}} \OO_{z} \dosflechas  C\otimes_BC=\prod_{z,z'\in U_{xy}} \OO_{zz'}$$
shows that  $B=\OO_{xy}$. 
Therefore the morphism $\OO_{xy}\to \prod_{z\in U_{xy}} \OO_{z}$ is faithfully flat. By Proposition \ref{afinfp},
$U_{xy}$ is affine.  
By Proposition \ref{K6}, $f$ is schematic.
\edemo

\section{Removable points. Minimal schematic space}

\prop \label{22} Let $X$ be a schematic finite space. Let $p\in X$ be a point such that the morphism $\OO_p\to \prod_{q>p}\OO_q$ is faithfully flat\footnote{If $I:=\{q\in X\colon q>p\}=\emptyset$, define $\prod_{q\in I}\OO_q:=\{0\}$.}. Consider the ringed subspace $Y=X-\{p\}$ of $X$ ($\OO_{Y,y}:=\OO_{X,y}$, for any $y\in Y$). Then, $Y$ is a schematic finite space,
 the inclusion map  $i\colon Y\iny X$ is an affine morphism and  $i_*\OO_Y=\OO_X$.\eprop

\demo $Y$ is a schematic finite space by Proposition \ref{4pg20}, because  $X$ is a schematic finite space. Let us prove that $i\colon Y\iny X$ is affine and  $i_*\OO_Y=\OO_X$:
Consider $U_x\subset X$. If $x\neq p$,  then $i^{-1}(U_x)=U_x$ and $(i_*\OO_Y)(U_x)=\OO_{Y,x}=\OO_{X,x}=\OO_X(U_x)$. If  $x=p$ denote $U:=U_p-\{p\}$. Observe that ${\OO_X}_{|U}={\OO_Y}_{|U}$.
The morphism $A=\OO_p\to
 \prod_{x\in U}\OO_x=B$ is faithfully flat. The exact sequence of morphisms
 $$A\to B\dosflechas B\otimes_AB$$
 and the equality $B\otimes_A B=\prod_{x,y\in U} \OO_{xy}$ show that $A=\OO_X(U)=\OO_Y(U)$. By Proposition \ref{afinfp}, $U$ is affine. Then, $i^{-1}(U_p)=U$ is affine and $(i_*\OO_Y)(U_p)=\OO_Y(U)=A=\OO_X(U_p)$.
 
  \edemo

\obse  If $U$ is affine then the morphism $\OO_U\to \prod_{q\in U} \OO_q$ is faithfully flat.
We have proved that  $\OO_p\to \prod_{q>p}\OO_q$ is faithfully flat iff $U:=U_p-\{p\}$ is affine and $\OO_p=\OO(U)$. \eobse

\lema \label{flat} Let $X$ be an affine finite space and $U\subseteq X$ an affine open subset. Then, the restriction morphism $\OO(X)\to \OO(U)$ is flat.

\elema

\demo The morphism $\OO(U)\to \prod_{x\in U} \OO_x$ is faithfully flat and the composite morphism
$\OO(X)\to \OO(U)\to  \prod_{x\in U} \OO_x$ is flat. Then, the morphism $\OO(X)\to \OO(U)$ is flat.

\edemo

\prop \label{bla}Let $X$ be a schematic finite space. Let $p\in X$ be a point such that the morphism $\OO_p\to \prod_{q>p}\OO_q$ is faithfully flat and let $Y:=X-\{p\}$. An open set $V\subseteq X$ is affine iff $V\cap Y$ is affine\eprop

\demo
We only have to prove the converse implication. We can suppose that $V=Y$ and we have to prove that $X$ is affine.

By Lemma \ref{flat}, the morphism $\OO_X(X)=\OO_Y(Y)
\to \OO_Y(U_x\cap Y)=\OO_X(U_x)=\OO_x$ is flat, for any $x\in X$. 
Then, the morphism $\OO(X)=\OO(Y)\to \prod_{y\in Y} \OO_{Y,y}=\prod_{y\in Y} \OO_{X,y}$ is faithfully flat. Therefore the morphism
$\OO(X)\to \prod_{x\in X} \OO_{X,x}$ is faithfully flat.
Likewise, given  $U_x,U_{x'}\subseteq X$, the morphism
$\OO_X(U_x\cap U_{x'})\to \prod_{x''\in U_x\cap U_{x'}} \OO_{X,x''}$
is faithfully flat. Besides,
$$\OO_x\otimes_{\OO(X)}\OO_{x'}=\OO_Y(U_x\cap Y)\otimes_{\OO(Y)}\OO_Y(U_{x'}\cap Y)\overset{\text{\ref{11}}}=\OO_Y(
U_x\cap Y\cap U_{x'}\cap Y)=\OO_X(U_x\cap U_{x'})$$
Therefore, $X$ is affine.\edemo

\defi Let $X$ be a schematic finite space. We shall say that $x\in X$ is removable if  $\OO_x\to \prod_{x'>x} \OO_{x'}$ is faithfully flat.\edefi

If $\OO_x=0$ obviously $x$ is a removable point.

\prop  Let $X$ be a schematic finite space and let $p,p'\in X$ be two  points. 
Then, $p,p'$ are removable points of $X$ iff $p$ is a removable  point of $X$ and $p'$ is a removable point of $X-p$.\eprop

\demo It is immediate.\edemo

\prop Let $U'\subseteq U$ be affine open subsets of  a schematic finite space $X$ and suppose that the morphism $\OO(U)\to \OO(U')$ is {faithfully} flat.  Then, $U-U'$ is a set of removable points of $X$.\eprop

\demo $\OO(U)=\OO(U')$ because  the morphism $\OO(U)\to \OO(U')$ is faithfully flat and 
$$\OO(U)\otimes_{\OO(U)}\OO(U')=\OO(U')
\overset{\text{\ref{11}}}=\OO(U')\otimes_{\OO(U)}\OO(U').$$ 
Let  $x\in U-U'$,  then $\OO_x=\OO(U)\otimes_{\OO(U)}\OO_x=\OO(U')\otimes_{\OO(U)}\OO_x
\overset{\text{\ref{11}}}=\OO(U'\cap U_x)$. By Proposition \ref{inter}, $U'\cap U_x$ is affine,  then the morphism $\OO_x=\OO(U'\cap U_x)\to \prod_{x'\in U'\cap U_x}\OO_{x'}$ is faithfully flat.
Hence,  $\OO_x\to \prod_{y>x} \OO_y$ is faithfully flat, and  $x$ is a removable point.\edemo

\obse In addition, we have proved that $\OO(U)=\OO(U')$.\eobse

\defi A schematic finite space $X$ is said to be minimal if there are no removable points in $X$ and it is $T_0$.  
Let $\tilde X$ the Kolmogorov space of $X$ and 
$P$ be the set of all the removable points of $\tilde X$, we shall denote
$X_M:=\tilde X\backslash P$. 

\edefi

By Proposition \ref{22}, $X_M$ is a schematic finite space, the natural morphism $X_M\overset{i}\iny \tilde X$ is affine and $\OO_{\tilde X}=i_*\OO_{X_M}$.

\prop Let  $f\colon X\to Y$ be a schematic morphism. If  $x\in X$ is not a removable point, then $f(x)\in Y$ is not a removable point. Then, we have the commutative diagram $$\xymatrix{ X \ar[d]^-f \ar[r] & \tilde X \ar[d]^-{\tilde f}& \,\, X_M \ar@{_{(}->}[l] \ar[d]^-{f_{M}}  \\  Y  \ar[r] & \tilde Y & \,\,  Y_M \ar@{_{(}->}[l]  }$$
where $\tilde X$ and $\tilde Y$ are the Kolmogorov spaces of $X$ and $Y$ respectively, and $\tilde f$ and $f_M$ are the induced morphisms. 
\eprop

\demo Consider the affine morphism $f_x\colon U_x\to U_{f(x)},\, f_x(x'):=f(x')$.  Since $f_{x*}\OO_{U_x}$ is a quasi-cohe\-rent   
$\OO_{U_{f(x)}}$-module, then  $\OO_x\otimes_{\OO_{f(x)} }\OO_y=
\OO(f_x^{-1}(U_y))$, for any $y\in U_{f(x)}$. If $f(x)$ is a removable point, then  $\OO_{f(x)}\to \prod_{y>f(x)} \OO_y$ is faithfully flat. Tensoring by  $\OO_x\otimes_{\OO_{f(x)}}$, one has the faithfully flat morphism
$\OO_x\to \prod_{y>f(x)} \OO(f_x^{-1}(U_y))$. The open sets $f_x^{-1}(U_y)$ are affine, because $f_x\colon U_x\to U_{f(x)}$ is affine. Then, the morphisms 
$\OO(f_x^{-1}(U_y))\to \prod_{x'\in f_x^{-1}(U_y)} \OO_{x'}$ are faithfully flat. Hence, the morphism $\OO_x\to \prod_{x'\in f_x^{-1}(U_y), y>f(x)}\OO_{x'}$ is faithfully flat. Therefore, $\OO_x\to \prod_{x'>x}\OO_{x'}$ is faithfully flat and   $x$ is a removable point.

\edemo

\prop \label{P6.11} Let $p\in Y$ be a removable point, $i\colon Y-\{p\}\to Y$ be the inclusion morphism and $f\colon X\to Y$ a schematic morphism. If $f(X)\subseteq Y-p$ and  $g\colon X\to Y-\{p\}$ is the morphism of ringed spaces such that $f=i\circ g$, then $g$ is a schematic morphism. 

Therefore, $f_{_M}\colon  X_M\to  Y_M$ is a schematic morphism.

\eprop

\demo It is an immediate consequence of Theorem \ref{K7}.
%We have to prove that $g_x\colon U_x\to U_{g(x)}$ is affine. We can %suppose
%$X=U_x$, $Y=U_{f(x)}$ and $g_x=g$. 
% Then,
%$$(g_*\OO_X)_y\otimes_{\OO_{y}}\OO_{y'}=(f_*\OO_X)_{i(y)}\otimes_{\OO_{i(y)}}\OO_{i(y')}=
%(f_*\OO_X)_{i(y')}=(g_*\OO_X)_{y'},$$
%for any $y'\geq y\in Y-p$. Then, $g_*\OO_X$ is quasi-coherent.
%Given an affine open subset $U\subseteq Y-\{p\}$, there exists an affine open subset $V\subseteq Y$ such that $V\cap (Y-\{p\})=U$. Then, $g^{-1}(U)=f^{-1}(V)$ is affine. Hence, $g$ is affine.

\edemo

\prop Let $X$ be a schematic finite space and  $U\subset X$ an affine open subspace.  Let  $X':=X\coprod \{u\}$ be the ringed finite space defined by

\enumera \item The 
preorder   on  $X\subset X'$ is the  pre-established preorder.  Given $x\in X$, then 
$u< x$ iff  $x\in U$, and  $x<u$ iff  $x\leq x'$, for any  $x'\in U$.

\item  $\OO_{X',x}:=\OO_{X,x}$ for any $x\in X$, and  $\OO_{X',u}:=\OO_X(U)$. The restriction morphisms are the obvious morphisms.
\eenumera

Then,
$X'$ is a schematic finite space and $u$ is a removable point of $X'$.\eprop

\demo Let us denote $U_u:=U\subset X$ and $\tilde U_{x'}:=\{y\in X'\colon y\geq x'\}$, for any $x'\in X'$. 
 By Proposition \ref{inter}, the morphism
$$\OO_{X',y}\otimes_{\OO_{X',x'}}\OO_{X',y'}\overset{\text{\ref{11}}}
= \OO_X(U_y\cap U_{y'})\to \prod_{x\in U_y\cap U_{y'}} \OO_{X,x}=
\prod_{x\in U_y\cap U_{y'}} \OO_{X',x}$$
is faithfully flat, for any $y,y'\geq x'$. If $U\subseteq U_y\cap U_{y'}$, the morphism $\OO_X(U_y\cap U_{y'})\to \OO_X(U)$ is flat, by Lemma \ref{flat}.
Hence, the morphism
$\OO_{X',y}\otimes_{\OO_{X',x'}}\OO_{X',y'}\to
\prod_{x\in \tilde U_y\cap \tilde U_{y'}} \OO_{X',x}$ is faithfully flat. By Proposition \ref{4pg20},
$X'$ is a schematic finite space.

The morphism $$\OO_{X',u}=\OO_X(U)\to \prod_{x\in U} \OO_{X,x}
=\prod_{x'>u} \OO_{X',x'}$$ is faithfully flat, because $U$ is affine. Hence, $u$ is removable.

\edemo

%\prop \label{preqc} Sea $Y$ un espacio finito af\'{i}n y $X$ un espacio finito esquem\'{a}tico minimal. Se cumple que la aplicaci\'{o}n
%$$\Hom_{sch}(X,Y)\to \Hom_{anillos}(\OO(Y),\OO(X)),\,\, f\mapsto f_Y$$
%es inyectiva.\eprop
%
%\demo Sea $f,g$ dos morfismos tales que $f_Y=g_Y$. Para probar que
%$f=g$ basta verlo localmente en $X$, luego podemos suponer que $X$ es af\'{i}n. Por tanto, $f$  y $g$ son un morfismos afines.
%Dados $y,y'\in Y$, tenemos que $\OO(f^{-1}(U_y))=\OO(X)\otimes_{\OO(Y)}\OO_y$ e igualmente
%$\OO(g^{-1}(U_y))=\OO(X)\otimes_{\OO(Y)}\OO_y$ y
%$$\aligned \OO(f^{-1}(U_y) & \cap g^{-1}(U_y))  =\OO(f^{-1}(U_y))\otimes_{\OO(X)}
%\OO(g^{-1}(U_y)) \\ & =(\OO(X)\otimes_{\OO(Y)}\OO_y)\otimes_{\OO(X)}
%(\OO(X)\otimes_{\OO(Y)}\OO_y)
% =\OO(X)\otimes_{\OO(Y)}\OO_y\otimes_{\OO(Y)}\OO_y
%\\ &  =\OO(X)\otimes_{\OO(Y)}\OO_y=\OO(f^{-1}(U_y))\endaligned$$
%Como $X$ es minimal, esto implica que $f^{-1}(U_y)\cap g^{-1}(U_y)=f^{-1}(U_y)$. Ahora es f\'{a}cil concluir que $f^{-1}(U_y)=g^{-1}(U_y)$, para todo $y$. Es f\'{a}cil concluir que las fibras de $f$ coinciden con las de $g$, es decir, que $f=g$ como aplicaciones de conjuntos. Es f\'{a}cil concluir que $f=g$, pues los morfismos
%$\OO_{f(x)}\to \OO_x$ est\'{a}n determinados por el morfismo $\OO(Y)\to \OO(X)\to \OO_x$.
%
%\edemo

\section{Serre Theorem}

Let  $X$ be a finite topological space and  $F$ a sheaf of abelian groups on  $X$.

\prop\label{aciclicity} If $X$ is a finite topological space with a minimum, then $H^i(X,F)=0$ for any sheaf $F$ and any $i>0$. In particular, for any finite topological space one has
\[ H^i(U_p,F)=0\]
for any $p\in X$, any sheaf $F$ and any $i>0$.
\eprop

\demo Let $p$ be the minimum of $X$. Then $U_p=X$ and, for any sheaf $F$, one has $\Gamma(X,F)=F_p$; thus, taking global sections is the same as taking the stalk at $p$, which is an exact functor.
\edemo

Let $f\colon X\to Y$ a continuous map between finite topological spaces and $F$ a sheaf on $X$. The i-th higher direct image $R^if_*F$ is the sheaf on $Y$ given by:
$$ [R^if_*F]_y=H^i(f^{-1}(U_y),F).$$

%\begin{rem} Let $X,Y$ be two finite topological spaces and $\pi\colon X\times Y\to Y$ the natural projection. If $X$ has a minimum ($X=U_x$), then, for any sheaf $F$ on $X\times Y$, $R^i\pi_*F=0$ for $i>0$, since $(R^i\pi_*F)_y=H^i(U_x\times U_y,F)=0$ by Proposition \ref{aciclicity}. In particular, $H^i(X\times Y, F)=H^i(Y,\pi_*F)$.
%\end{rem}

Let $F$ be a sheaf on a finite topological space $X$. We define $C^nF$ as the sheaf on $X$  whose sections on an open subset $U$ are
$$ (C^nF)(U)=\proda{U \ni x_0<\cdots <x_n } F_{x_n}$$and whose  restriction morphisms $(C^nF)(U)\to (C^nF)(V)$ for any $V\subseteq U$ are the natural projections.

One has morphisms $d\colon C^nF \to C^{n+1}F$, $a=(a_{x_0<\cdots < x_{n}})\mapsto d(a)=(d(a)_{x_0<\cdots < x_{n+1}})$ defined in each open subset $U$ by the formula
$$ (\di a) _{x_0<\cdots < x_{n+1}}= \suma{0\leq i\leq n} (-1)^i a_{x_0<\cdots \wh{x_i}\cdots <x_{n+1}} + (-1)^{n+1} \bar a _{x_0<\cdots <x_n}   $$ where $\bar a _{x_0<\cdots <x_n}$ denotes the image of  $ a _{x_0<\cdots <x_n}$ under the morphism $F_{x_n}\to F_{x_{n+1}}$. There is also a natural morphism   $\di\colon F\to C^0F$. One easily checks that $\di^2=0$.

\teor[(\cite{KG} 2.15)\,]  $C^\punto F$ is a finite and flasque resolution of  $F$ (in fact, it is the Godement resolution of $F$).
\eteor

\demo By definition, $C^nF=0$ for $n>\dim X$. It is also clear that  $C^nF$ are flasque. Let us see that
\[ 0\to F\to C^0F \to \cdots\to C^{\dim X}F\to 0\] is an exact sequence.  We have to prove that  $(C^\punto F)(U_p)$ is a resolution of $F(U_p)$. One has a decomposition
\[ (C^nF)(U_p)= \proda{p=x_0<\cdots <x_n } F_{x_n}\times \proda{p<x_0<\cdots <x_n } F_{x_n} = (C^{n-1}F)(U^*_p)\times (C^nF)(U^*_p)\] with $U_p^*:=U_p-\{ p\}$; via this decomposition, the differential  $\di \colon (C^nF)(U_p) \to (C^{n+1}F)(U_p)$ becomes:
\[ \di(a,b)=(b-\di^*a,\di^*b)\] with $\di^*$ the differential of   $(C^\punto F)(U_p^*)$. 
If $d(a,b)=0$, then $b=d^*a$ and $d(0,a)=(a,b)$.

It is immediate now that every cycle is a boundary.
\edemo

This theorem, together with De Rham's theorem (\cite{Godement}, Thm. 4.7.1), yields that the cohomology groups of a sheaf can be computed with the standard resolution, i.e., $H^i(U,F)=H^i\Gamma(U,C^\punto F)$, for any open subset $U$ of $X$ and any sheaf $F$ of abelian groups on $X$.

%\coro For any finite topological space $X$, any sheaf $F$ of abelian groups on $X$ and any family of supports $\phi$, one has
%\[ H^n_\phi(X,F)=0,\quad \text{for any } n>\text{\rm dim} X.\] Moreover, if $F_p$ is a finitely generated $\ZZ$-module for any $p\in X$, then
%$H^i_\phi(X,F)$ is a finitely generated $\ZZ$-module for any $i\geq 0$.
%\ecoro

\teor[(\cite{KS} 4.3,\,4.12)\,] \label{afin} Every quasi-coherent module of an  affine finite space is acyclic.
\eteor

\demo Let $X$ be an affine finite space. 

Proceed by induction over the order of $X$.  The open sets  $U_{xy}$ are affine by  Corollary \ref{Uxy}. 
If $X=U_{xy}$, then $X=U_x$ and every sheaf is acyclic. 
We can suppose that every quasi-coherent module on $U_{xy}$ is acyclic, by induction hypothesis.
Let us prove that any quasi-coherent $\OO_X$-module $\M$ is acyclic. We have to prove that the sequence of morphisms 
$$\M(X)\to \prod_{x_1\in X} \M_{x_1}\to \prod_{x_1<x_2} \M_{x_2}\to \cdots$$
is exact. It is sufficient to check that tensoring  the previous sequence by 
 $\otimes_{\OO(X)}\OO_z$, for any $z\in X$, the sequence of morphisms  
$$\xymatrix @C18pt{ \M_z\ar[r] &   \proda{z\leq x_1} \M_{x_1} \times 
\proda{z\not\leq x_1} \M_{x_1z} \ar[r] &  \proda{z\leq x_1<x_2} \M_{x_2}\,\, \times\!
\proda{\scriptsize \alinea \,\,\, \,\, \, \, x_1<x_2 \\z\not\leq x_1, z\leq x_2\ealinea} \M_{x_2}\,\, \times\! \proda{\scriptsize \alinea x_1<x_2\\ \, \,\, z\not\leq x_2\ealinea} \M_{x_2z} \, \ar[r] &\,\,  \cdots}$$
is exact. That is, we have to prove that the sequence of morphisms  $(S)$

$$\xymatrix  @C18pt @R8pt{
\M(U_z) \, \ar[r] & \, 
 C^0(U_z,\M) \times \proda{z\not\leq x_1} \M(U_{x_1z})    \ar[r] & \, 
C^1(U_z,\M)\times \proda{z\not\leq x_1} C^0(U_{x_1z},\M)\,\,\times \proda{\scriptsize \alinea x_1<x_2\\ \,\,\, z\not\leq x_2 \ealinea} \! \M(U_{x_2z})
\\ \qquad \,\, \ar[r] & \quad 
\cdots  \qquad \qquad \qquad \qquad  & 
 }$$ 
is exact. Let
$D^\punto_r:=\oplus_{x_1<\cdots<x_r,z\not\leq x_r}( \M(U_{x_rz})\oplus C^\punto  ({U_{x_rz}},\M))[-r]$ and let $d_r$ be the differential such that over  each  direct summand $ (\M(U_{x_rz})\oplus C^\punto  ({U_{x_rz}},\M))[-r]$  is the known differential of $\M(U_{x_rz})\oplus C^\punto  ({U_{x_rz}},\M)$ multiplied by $(-1)^r$. $H^i(D^\punto_r)=0$ for any $i\geq 0$.
The sequence of morphisms $(S)$ is equal to the differential complex $D^\punto:=D^\punto_0\oplus D^\punto_1\oplus\cdots \oplus D^\punto_n$ with the diffrential 

$$d=\pamatrix{d_0 & 0 & 0 &  \cdots & 0 & 0\\ - & d_1 & 0 & \cdots & 0& 0\\ \vdots & \vdots & \ddots & & \vdots  & \vdots \\ \vdots & \vdots & \vdots & \ddots & \vdots  & \vdots \\ - & - & - & \cdots & d_{n-1} & 0 \\
- & - & - & \cdots & - & d_n}$$
Let  $D_{>0}^\punto =\oplus_{i>0} D^\punto_i$.
Consider the exact sequence of morphisms  of complexes $$0\to D_{>0}^\punto\to D^\punto \to D^\punto_0\to 0.$$ Then,  $H^i(D^\punto)=H^i(D_{>0}^\punto)$, for any $i\geq 0$.  Let $D_{>1}^\punto =\oplus_{i>1} D^\punto_i$.
Consider the exact sequence of morphisms  $0\to D_{>1}^\punto \to D_{>0}^\punto \to D^\punto_1\to 0$. Then, $H^i(D^\punto)=H^i(D_{>0}^\punto)=H^i(D_{>1}^\punto)$.
Recursively  $H^i(D^\punto)=H^i(D_{>n}^\punto)=0$ for any $i\geq 0$, and the sequence of morphisms  $(S)$ is exact.

\edemo

Let $R$ and $R'$ be commutative rings and $R\to R'$ a flat morphism of rings.
Let $(X,\OO)$ be an $R$-ringed finite space. Let $\OO\otimes_R R'$ be the sheaf of rings on $X$ defined by
$(\OO\otimes_R R')(U):=\OO(U)\otimes_R R'$. Consider the obvious morphism $\pi\colon (X,\OO\otimes_R R')\to (X,\OO)$, $\pi(x)=x$. Let $\M$ be a sheaf of $\OO$-modules. Then,
$$H^i(X,\pi^*\M)=H^i(X,\M)\otimes_R R'.$$
If $\Nc$ is a quasi-coherent $\OO\otimes_R R'$-module, then $\pi_*\Nc=\Nc$ is a quasi-coherent $\OO$-module. 

Let $S\subset R$ be a multiplicative system, $R'=S^{-1}\cdot R$ and   $\Nc$ a quasi-coherent $\OO\otimes_R R'$-module. Then,
$\pi_*\Nc$ is a quasi-coherent $\OO$-module,  $\Nc=\pi^*\pi_*\Nc$ and $H^i(X,\Nc)=H^i(X,\pi_*\Nc).$

\nada[Serre Theorem (\cite{KS} 5.11)\,] \label{Serre}  Let $X$ be a schematic finite space. $X$ is affine  iff every quasi-coherent $\OO_X$-module $\M$ is acyclic (or $H^1(X,\M)=0$).\enada

\demo $\Leftarrow)$ Let $R:=\OO(X)$.  Recall Notation \ref{Notation}. Given $\pp\in \Spec R$, consider the sheaf of rings on $X$, $\OO\otimes_RR_\pp$. Obviously,  $(X,\OO\otimes_RR_\pp)$ is a schematic finite space and $(X,\OO)$ is  affine iff  $(X,\OO\otimes_RR_\pp)$ is affine for any $\pp$.
Hence, we can suppose $R$ is a local ring. We can suppose that $X$ is minimal. Let $X'$ be the set of the closed points of $X$. Let $x'\in X'$. The morphism $\OO_{x'}\to \prod_{x> x'} \OO_{x}$ is flat but it is not faithfully flat, then there exists a prime ideal $I_{x'}\subset \OO_{x'}$ such that
$I_{x'}\cdot \prod_{x>x'} \OO_{x}=  \prod_{x>x'} \OO_{x}$. Let $\pp$ be the quasicoherent ideal defined by
$\pp_{x'}:=I_{x'}$ if $x'\in X'$ and $\pp_x:=\OO_x$ if $x\notin X'$.
Observe that $(\OO_X/\pp)_x=0$, for any $x\in X\backslash X'$, then
$(\OO_{X}/\pp)(X)=\prod_{x'\in X'} \OO_{X}/I_{x'}$.
 Consider the exact sequence of morphisms
$$0\to \pp \to \OO\to \OO_{X}/\pp\to 0$$
The morphism $R=\OO(X)\to (\OO_{X}/\pp)(X)=\prod_{x'\in X'} \OO_{X}/I_{x'}$ is surjective, because $H^1(X,\pp)=0$. $R$ is a local ring,  then $\prod_{x'\in X'} \OO_X/I_{x'}$  is a local ring, hence 
$X'=\{x'\}$. Therefore, $X=U_{x'}$, which is affine.

\edemo

This theorem yields the usual Serre's criterion on algebraic varieties (see \cite{KG} 4.13 and \cite{Serree}).

\coro  \label{corotonto} A schematic finite space  $X$ is affine  iff
the functor $$\Gamma\colon {\bf Qc\text{-}Mod}_X\to {\bf Mod}_{\OO(X)},\,\,
\M\mapsto \Gamma(X,\M)$$ is exact.\ecoro

\demo $\Rightarrow)$ By  the Serre Theorem $H^1(X,\M)=0$, for any quasi-coherent module $\M$. Hence $\Gamma$ is exact.

$\Leftarrow$) It has been proved in the proof of Serre Theorem.

\edemo

\coro A schematic finite space $X$ is affine iff $H^1(X,\Ic)=0$ for any quasi-coherent ideal $\Ic\subseteq \OO$.\ecoro 

\demo $\Leftarrow)$  Let $R=\OO(X)$. We have just proved this implication when $R$ is a local ring, in the  proof of the Serre Theorem. Let $\pp\in\Spec R$
and let $\Ic^\pp\subset \OO\otimes_RR_\pp$ be a quasi-coherent ideal. Consider the obvious morphism 
$\OO\to  \OO\otimes_RR_\pp$.  $\Jc:=\OO\times_{\OO\otimes_R R_\pp}\Ic^\pp$ is a quasi-coherent ideal of $\OO$ and
$\pi^*\Jc=\Ic^\pp$ (where $\pi\colon (X,\OO\otimes_R R_\pp)\to (X,\OO)$ is defined by $\pi(x):=x$).  Then, $H^1(X,\Ic^\pp)=H^1(X,\Jc)\otimes_R R_\pp=0$. Then, $(X,\OO\otimes_RR_\pp)$ is affine, for any $\pp$. Therefore, $X$ is affine.

\edemo

\coro  Let  $X$ be a minimal  affine finite space   and suppose that $\OO(X) $ is a local ring. Then, there exists a  point $p\in X$ such that 
$X=U_p$.\ecoro

\teor \label{R1} Let  $X$ and $Y$ be schematic finite spaces. A ringed space morphism $f\colon X\to Y$ is affine iff  $f_*\M$ is quasi-coherent and   $R^1f_*\M=0$, for any quasi-coherent $\OO_X$-module $\M$.\eteor

\demo $\Rightarrow)$ It is obvious.

$\Leftarrow)$ 
Let $U\subseteq Y$ be an affine open subspace. 
By Corollary \ref{C4.10}, any quasi-coherent module $\M$ on $f^{-1}(U)$ is the restriction of a quasi-coherent module on $X$.
$H^1(f^{-1}(U),\M)=H^1(U,f_*\M)=0$, for any quasi-coherent  $\OO_{X}$-module $\M$. By Serre Theorem \ref{Serre},
$f^{-1}(U)$ is affine. Hence, $f$ is affine.
\edemo

\hskip-0.65cm\colorbox{white}{\,\begin{minipage}{15.15cm}
\teor \label{Tafex} Let $f\colon X\to Y$ be a schematic morphism. The functor $$f_*\colon \,{\bf Qc\text{-}Mod}_X\, \to \,{\bf Qc\text{-}Mod}_Y,\quad \mathcal M\functor f_*\mathcal M$$ is 
%faithfully 
exact iff $f$ is affine.\eteor\end{minipage}}

\demo $\Rightarrow)$ 1. Let  $U\subseteq Y$ be an open subset, $V=f^{-1}(U)$ and $f_{|V}\colon V\to U$, $f_{|V}(x):=f(x)$. Then $f_{|V*}$ is 
%faithfully 
exact: Any short exact sequence of quasi-coherent modules $\mathcal N^\bullet$ on $V$ is a restriction of a short exact sequence of quasi-coherent modules $\mathcal M^\bullet$ on $X$; and $f_{V*}\mathcal N^\bullet=(f_*\mathcal M^\bullet)_{|U}$.
%Hence, $f_{|V*}$ is exact. 
%Finally, 
%let $\mathcal N\neq 0$ be a quasi-coherent $\OO_V$-module. Observe that $i_*f_{|V*}\mathcal N=f_*j_*\mathcal N\neq 0$, where
%$i\colon U\iny X$ and $j\colon V\iny Y$ are the obvious inclusions. Hence, $f_{|V*}\mathcal N\neq 0$.

2. We can suppose that $Y$ is affine. We can suppose that $Y=(*,A)$.  We can suppose that $A=\OO_X(X)$.

3. We have to prove that $X$ is affine. The functor, 
$$\,{\bf Qc\text{-}Mod}_X\, \to \,{\bf Mod}_{\OO(X)},\quad \mathcal M\functor f_*\mathcal M=\Gamma(X,\M)$$ is 
exact. By Corollary \ref{corotonto}, $X$ is affine.

%3. Let $\pp\in\Spec A$ be a prime ideal. The morphism $f_{\pp*}$ induced by the obvious morphism $f_\pp\colon (X, \OO_X\otimes_A A_\pp)\to
%(*,A_\pp)$ is 
%exact: Observe that any quasi-coherent $\OO_X\otimes_A A_\pp$-module
%can be considered as a quasi-coherent $\OO_X$-module and $f_{\pp,*}=f_*$.
%
%4. We have to prove that $X$ is affine, which is equivalent to saying that $(X, \OO_X\otimes_A A_\pp)$ is affine, for any $\pp\in \Spec A$. We can suppose that $A$ is a local ring and that $X$ is a minimal schematic finite space.
%
%5. By the hypotheses of the theorem, if a morphism of quasi-coherent $\OO_X$-modules  
%$\M_1\to \M_2$ is an epimorphism then
%$\M_1(X)\to \M_2(X)$ is surjective. 
%Let  $X'=\{x_1,\ldots,x_n\}$ be the set of the closed points of  $X$.  The morphism
%$\OO_{x_i}\to \prod_{z>x_i}\OO_z$ is flat and it is not faithfully flat. Let $I_i\underset\neq \subset \OO_{x_i}$ an ideal such that $I_i\cdot  \prod_{z>x_i}\OO_z=  \prod_{z>x_i}\OO_z$. 
%Let  $\OO'$ be the sheaf of rings and the quasi-coherent $\OO_X$-module defined by $\OO'_{x}:=\OO_{x_i}/I_i$, if  $x=x_i$,  for some  $i$,
%and  $\OO'_{x}:=0$ in other case. Consider the obvious epimorphism  $\OO_X\to \OO'$.
%Then,  $A=\OO_X(X)\to \OO'(X)$
%is surjective.   $\OO'(X)=\prod_{i=1}^n \OO_{x_i}/I_i $ is a local ring because $A$ is local.
%Then, $n=1$ and $X=U_{x_1}$, which is affine.

$\Leftarrow)$ By Theorem \ref{R1}, $R^1f_*\M=0$, for any quasi-coherent module $\M$. Then, $f_*$ is exact.
%Let $\mathcal M^\bullet$ be a short exact sequence of quasi-coherent $\OO_X$-modules.
%We have to prove that $f_*\mathcal M^\bullet$ is  a short exact sequence of quasi-coherent $\OO_Y$-modules.
%It is a local property on $Y$, then we can suppose that $X$ and $Y$ are affine.
%Then,  $\mathcal M^\bullet(X)$ is  a short exact sequence of $\OO(X)$-modules, then
%$(f_*\mathcal M^\bullet)(Y)=\mathcal M^\bullet(X)$ is  a short exact sequence of $\OO(Y)$-modules and $f_*\mathcal M^\bullet$ is  a short exact sequence of quasi-coherent $\OO_Y$-modules.
\edemo

\section{Cohom. characterization of schematic finite spaces}

We say that an $R$-ringed space $(X,\OO_X)$ is a flat $R$-ringed space if the morphism $R\to \OO_x$ is flat, for any $x\in X$.

\teor \label{2} Let  $(X,\OO)$ be a flat $R$-ringed finite space and $M$ an $R$-module.
If  $H^i(X,\OO)$ is a flat $R$-module, for any  $i>0$,  then $\OO(X)$ is a flat $R$-module and  $$H^i(X,\tilde M)=H^i(X,\OO)\otimes_{R} M,\quad \forall i.$$
\eteor

\demo Let $C^i:=\Ker[C^i(X,\OO)\to C^{i+1}(X,\OO)]$ and $B^i:=\Ima[C^{i-1}(X,\OO)\to C^{i}(X,\OO)]$.  Let $n=\dim X$. $C^n=C^n(X,\OO)$ is $R$-flat.
The sequence $0\to B^n\to C^n\to H^n(X,\OO)\to 0$ is exact, then 
$B^n$ is flat. The sequence $0\to C^{n-1}\to C^{n-1}(X,\OO)\to B^n\to 0$ is exact then $C^{n-1}$ is $R$-flat. Recursively,  $C^i$ and $B^i$ are $R$-flat for any $i$.  Hence, $\OO(X)=C^0$ is $R$-flat and
$$H^i(X,\tilde M)=H^i(X,C^\punto\OO\otimes_R M)=H^i(X,C^\punto\OO)\otimes_R M=H^i(X,\OO_X)\otimes_R M.$$

\edemo

\coro \label{coro2} Let $X$ be a flat  $\OO(X)$-ringed finite space.  Assume $H^i(X,\OO)$ is a flat $\OO(X)$-module, for any $i>0$. Then, the morphism  $\OO(X)\to \prod_{x\in X} \OO_x$ is faithfully flat.\ecoro

\demo Let  $M$ be an $\OO(X)$-module. By Theorem \ref{2},  $\tilde M(X)=M$. Then, the morphism
$$M=\tilde M(X)\iny \prod_{x\in X} \tilde  M_x=  M\otimes_{\OO(X)} \prod_{x\in X} \OO_x$$ is injective.
Therefore,  the flat morphism $\OO(X)\to \prod_{x\in X} \OO_x$ is faithfully flat.\edemo

\teor \label{afin'} Let $X$ be a flat $\OO(X)$-ringed finite space.  Then,  $X$ is affine iff \enumera
\item[1'.] 
$X$ is acyclic.

\item[2'.]  $\OO_x\otimes_{\OO(X)} \OO_y=\OO_{xy}$, for any $x,y$,

\item[3'.] $U_{xy}$ is acyclic, for any $x,y$.\eenumera 

\eteor
  
  \demo $\Rightarrow)$ $U_{xy}$ is affine by  Corollary \ref{Uxy}.  By Theorem \ref{afin},  $X$ and $U_{xy}$ are acyclic. 
  
$\Leftarrow)$  Let $z\in U_{xy}$. The morphism $$\OO_{xy}\to \OO_z\otimes_{\OO(X)} \OO_{xy}=
\OO_z\otimes_{\OO(X)} \OO_x\otimes_{\OO(X)}\OO_y=\OO_{zx}\otimes_{\OO(X)}\OO_y=\OO_z\otimes_{\OO(X)}\OO_y=\OO_{zy}=\OO_z$$
is flat, since the morphism $\OO(X)\to \OO_z$ is flat.

By Corollary \ref{coro2}, the morphisms  $\OO(X)\to \prod_{x\in X}\OO_x$ and 
  $\OO(U_{xy})\to \prod_{z\in U_{xy}}\OO_z$ are faithfully flat. 
  Hence, $X$ is affine.
  
\edemo

\teor \label{T5.9} A finite fr-space $X$ is schematic  iff  for any $x\leq y,y'$, 

\enumera 

 \item $\OO_{y}\otimes_{\OO_x}\OO_{y'}=\OO_{yy'}$.

\item $U_{yy'}$ is acyclic.\eenumera\eteor

\demo $\Rightarrow)$  $U_x$ is an affine finite space. By Theorem \ref{afin'}, we are done.

$\Leftarrow)$  $U_x$ is an affine finite space by Theorem \ref{afin'}, then $X$ is schematic.\edemo

\prop \label{12} Let  $X$ be an affine finite space. An open subset  $U\subseteq X$ is affine iff it is  acyclic.\eprop

\demo $\Leftarrow)$ $U$ satisfies  1'. and  3' of Theorem \ref{afin'}.
The composite morphism of the epimorphism   $\OO_x\otimes_{\OO(X)}\OO_{y}\to 
\OO_x\otimes_{\OO(U)} \OO_y$ and the morphism  $ 
\OO_x\otimes_{\OO(U)} \OO_{y}\to \OO_{xy}$ is an isomorphism, then $\OO_x\otimes_{\OO(U)} \OO_{y}\to \OO_{xy}$ is an isomorphism.

Besides,  the morphism $\OO(U)=\OO(X)\otimes_{\OO(X)}\OO(U)\to \OO_x\otimes_{\OO(X)}\OO(U)\overset{\text{\ref{11}}}=\OO_x$ is flat.

\edemo

\prop \label{K4} Let $f\colon X\to Y$ be a schematic morphism  and $\M$  a quasi-coherent $\OO_X$-module. 
Then,  $R^if_*\M$ is a quasi-coherent  $\OO_Y$-module,  for any  $i\geq 0$. If $Y$ is affine, $R^if_*\M=\wt{H^i(X,\M)}$. \eprop

\demo We can suppose that $Y$ is affine. 
Given $y\in Y$, the inclusion morphism $f^{-1}(U_y)\overset{j}\iny X$ is affine: 
Let $U\subset X$ be an affine subset and $i$ be the composite morphism $U\iny X\overset f\to Y$, which is an affine morphism. Then, $j^{-1}(U)=f^{-1}(U_y)\cap U=i^{-1}(U_y)$ is affine.

Observe that $R^nf_*(j_*\Nc)=R^n(f\circ j)_*\Nc=0$ for any $n>0$ and any quasi-coherent module $\Nc$, since $j$ and $f\circ j$ are affine morphisms. Likewise, $H^n(X,j_*\Nc)=H^n(f^{-1}(U_y),\Nc)=H^n(Y, (f\circ j)_*\Nc)=0$, for any $n>0$.

Denote $\M_{f^{-1}(U_y)}=j_*\M_{|f^{-1}(U_y)}$ and consider the obvious exact sequence of morphisms
$$0\to \M\to \oplus_{y\in Y} \M_{f^{-1}(U_y)}\overset\pi\to \M'\to 0.$$
Then, $H^1(X,\M)=\Coker \pi_X$ and $H^{n-1}(X,\M')=H^n(X,\M)$
for any $n>1$. Besides,  $R^1f_*\M=\Coker \pi$, which is quasi-coherent,  and $R^nf_*\M=R^{n-1}f_*\M'$ for any $n>1$. 
Therefore, $R^1f_*\M=\widetilde{H^1(X,\M)}$ since $Y$ is affine. 
Hence, $R^1f_*\M'=\widetilde{H^1(X,\M')}$, since $\M'$ is quasi-coherent. By induction on $n$,
$$R^nf_*\M=R^{n-1}f_*\M'=\widetilde{H^{n-1}(X,\M')}=
\widetilde{H^{n}(X,\M)}$$
for any $n>1$.

%
%
%Given an open set $U\overset i\subseteq X$, 
%denote $\M_{U}=i_*{\M}_{|U}$, which is quasi-coherent by Proposition \ref{K40}. 
%$$\M(f^{-1}(U_y)\cap U)=\M((f\circ i)^{-1}(U_y))\overset{\text{\ref{K40}}}=\M((f\circ i)^{-1}(Y))\otimes_{\OO(Y)} \OO_y=
%\M(U)\otimes_{\OO(Y)} \OO_y,$$ that is, $\M_{f^{-1}(U_y)}(U)=\M(U)\otimes_{\OO(Y)}\OO_y$.  Consider the standard chain complex $\M(X)\to C^\punto\M(X)$. Tensoring by $\otimes_{\OO(Y)}\OO_y$ we obtain the chain complex $\M_{f^{-1}(U_y)}(X)\to  C^\punto\M_{f^{-1}(U_y)}(X)$ and  $$(a)\qquad H^i(X,\M)\otimes_{\OO(Y)}\OO_y=H^i(X,\M_{f^{-1}(U_y)}),$$ because $\OO(Y)\to \OO_y$ is a flat morphism.
%
%Let $x\in X$. The composite morphism  $U_x\to U_{f(x)}\iny Y$ is an affine morphism. Then, $U_x\cap f^{-1}(U_y)$ is affine, for any $x\in X$ and $y\in Y$. Therefore, $f^{-1}(U_y)\iny X$ is an affine morphism and $$(b)\qquad H^i(f^{-1}(U_y),\M)=H^i(X,\M_{f^{-1}(U_y)}),$$ for any $y\in Y$. Then,
%$$H^i(X,\M)\otimes_{\OO(Y)}\OO_y\overset{(a)}=H^i(X,\M_{f^{-1}(U_y)})\overset{(b)}=H^i(f^{-1}(U_y),\M).$$ Hence, $R^if_*\M=\wt{H^i(X,\M)}$.
%
\edemo

This proposition, Theorem \ref{Tcohsch} and \cite{KG} Theorem 5.6 show that
the definitions of schematic morphism given in this paper and in \cite{KG} are equivalent.

\lema  \label{K1} Let  $X$ be an $fr$-space and $\delta\colon X\to X\times X$, $\delta(x):=(x,x)$ be the diagonal morphism. Let $\M$ be an $\OO_X$-module.
$R^i\delta_*\M$ is quasi-coherent iff $$H^i(U_{pq},\M)\otimes_{\OO_p}\OO_{p'}=H^i(U_{p'q},\M),$$ for any $p\leq p'$ and for any  $q$.\elema

\demo $(R^i\delta_*\M)_{(q,q')}=H^i(U_{qq'},\M)$ and
$$\aligned (R^i\delta_*\M)_{(p,p')}\otimes_{\OO_{(p,p')}} \OO_{(q,q')} & =H^i(U_{pp'},\M)\otimes_{\OO_p\otimes \OO_{p'}} (\OO_{q}\otimes \OO_{q'})\\ & =
 (H^i(U_{pp'},\M)\otimes_{\OO_p}\OO_q)\otimes_{\OO_{p'}}\OO_{q'},\endaligned$$
for any $(p,p')\leq (q,q')$. Now, the proof  is easily completed.

\edemo

%\lema \label{K2} Let $X$ be a ringed finite space and assume the morphism $\OO_p\to \OO_q$ is flat, for any $p\leq q$. Then, 
%$\OO_{pq}\otimes_{\OO_p}\OO_{p'}=\OO_{p'q}$, for any $p\leq p'$ and for any  $q$ iff
%$\OO_{p'}\otimes_{\OO_p}\OO_{p''}=\OO_{p'p''}$, for any $p\leq p' ,p''$.\elema
%
%\demo $\Rightarrow)$ $\OO_{p'}\otimes_{\OO_p}\OO_{p''}=\OO_{pp'}\otimes_{\OO_p}\OO_{p''}=\OO_{p'p''}$.
%
%$\Leftarrow)$ Let $U\subseteq U_p$ be an open set  and $p\leq p'$. Consider the exact sequence of morphisms
%$$\OO(U)\to \prod_{x\in U} \OO_x\dosflechas \prod_{U\ni x\leq z} \OO_z.$$
%Tensoring by $\otimes_{\OO_p} \OO_{p'}$ we obtain the exact sequence of morphisms
%$$\OO(U)\otimes_{\OO_p} \OO_{p'}\to \prod_{x\in U} \OO_{p'x}\dosflechas \prod_{U\ni x\leq z} \OO_{p'z},$$
%which shows that $\OO(U)\otimes_{\OO_p} \OO_{p'}=\OO(U\cap U_{p'})$. In particular,
% $$\OO_{pq}\otimes_{\OO_p}\OO_{p'}=\OO(U_{pq}\cap U_{p'})=\OO_{p'q}.$$\edemo
% 
 
 \prop Let $X$ be an $fr$-space and $\delta\colon X\to X\times X$  the diagonal morphism. Let $\M$ be an $\OO_X$-module. Then, $\delta_*\M$ is a quasi-coherent  $\OO_{X\times X}$-module iff   
 $\M_{p'}\otimes_{\OO_p}\OO_{p''}=\M_{p'p''}$, for any $p\leq p',p''$.\eprop
 
 \demo It is a consequence of Lemma \ref{K1} and Proposition \ref{K2?}. 
 \edemo

\vaci[Cohomological characterization of schematic finite spaces (\cite{KS} 4.7,\,4.4\,):]  \label{ccefe} Let  $X$ be an $fr$-space and   $\delta\colon X\to  X\times X$  the diagonal morphism. $X$ is a schematic finite space iff  $R^i\delta_*\OO_X$ is a quasi-coherent module, for any $i\geq 0$.\evaci

\demo

 $\Leftarrow)$ We have to prove that $U_p$ is affine. $U_p$ is acyclic, and satisfies the property 2' of Theorem \ref{afin'}, by the previous proposition. We only need to prove that $U_{qq'}$ is acyclic, for any $q,q'\in U_p$:
 $$0=H^i(U_{q},\OO)\otimes_{\OO_p}\OO_{q'}=H^i(U_{pq},\OO)\otimes_{\OO_p}\OO_{q'}\overset{\text{\ref{K1}}}=H^i(U_{q'q},\OO).$$

$\Rightarrow)$ The diagonal morphism $\delta$ is schematic by 
Theorem \ref{Cdelta} and Theorem \ref{Tcohsch}. By Proposition \ref{K4}, we are done.

\edemo

\coro[(\cite{KG} 4.5)\,] An $fr$-space $X$ is schematic iff  for any open set $j\colon U_q\iny X$, $R^ij_*\OO_{U_q}$ is a quasi-coherent  $\OO_X$-module, for any $i$.\ecoro

\demo Let $\delta\colon X\to X\times X$ be the diagonal morphism. Then,
$X$ is a schematic finite space iff  $R^i\delta_*\OO_X$ is a quasi-coherent module, for any  $i$, which is equivalent to say that  $H^i(U_{pq},\OO)\otimes_{\OO_p}\OO_{p'}=H^i(U_{p'q},\OO)$,  for any  $p\leq p'$, and any  $q$, that is to say, $R^ij_*\OO_{U_q}$ is a quasi-coherent $\OO_X$-module, for any  $i$ and any open set $j\colon U_q\iny X$.

\edemo

A scheme is said to be a semiseparated scheme if the intersection of two affine open sets is affine. For example, the line with a double point is a semiseparated scheme (but it is not separated). The plane
with a double point is not semiseparated, but it is quasi-separated.

\defi A ringed finite space $X$ is said to be semiseparated if the open sets $U_{pq}$ are acyclic, for any $p,q\in X$.\edefi

\prop Let $X$ 
be a ringed finite space and let $\delta\colon X\to  X\times X$ be the diagonal morphism. $X$ is seimiseparated iff $R^i\delta_*\OO_X=0$, for any $i>0$.\eprop

\demo 
$(R^i\delta_*\OO_X)_{(p,q)}=H^i(U_{pq},\OO)$, and $R^i\delta_*\OO_X=0$ iff $(R^i\delta_*\OO_X)_{(p,q)}=0$ for any $p,q\in X$. Hence, $X$ is semiseparated iff $R^i\delta_*\OO_X=0$, for any $i>0$.
\edemo 

\teor Let $X$ 
be an $fr$-space and $\delta\colon X\to  X\times X$ be the diagonal morphism. $X$ is a semiseparated schematic finite space iff  $R^i\delta_*\OO_X=0$, for any $i>0$ and $\delta_*\OO_X$ is a quasi-coherent module.\eteor

\demo $\Rightarrow)$
By \ref{ccefe}, $\delta_*\OO_X$ is quasi-coherent, and 
$R^i\delta_*\OO_X=0$ by the previous proposition.
 
$\Leftarrow)$ $X$ is a schematic finite space, by \ref{ccefe}, and it is semiseparated by the previous proposition.
 \edemo

\prop A schematic finite space is semiseparated iff it satisfies any of the following equivalent conditions: 

\enumera 
\item The intersection of any two affine open subspaces is affine.

\item There exists an affine open covering of $X$, $\U=\{U_1,\ldots, U_n\}$ such that $U_i\cap U_j$ is affine
for any $i,j$.
\eenumera
\eprop

\demo  Assume $X$ is a semiseparated schematic finite space. Let $\delta\colon X\to X\times X,$ $\delta(x)=(x,x)$ be the diagonal morphism and  $U$, $U'$ two affine open subspaces. Since $R^i\delta_*\OO_X=0$, for any $i>0$,   $H^i(U\cap U',\OO)=H^i(U\times U',\delta_*\OO)=0$. By  \ref{12} $U\cap U'$ is affine. 

Assume  that there exists an affine open covering of $X$, $\U=\{U_1,\ldots, U_n\}$ such that $U_i\cap U_j$ is affine
for any $i,j$. $R^i\delta_*\OO_X$ is quasi-coherent, by \ref{ccefe}. $R^i\delta_*\OO_X(U_i\times U_j)=H^i(U_i\cap U_j,\OO_X)=0$, then $R^i\delta_*\OO_X=0$ and $X$ is semiseparated.

%Given $p,q\in X$, there exist  $i,$ and $j$ such that $U_p\subseteq U_i$ and $U_q\subset U_j$.
%Then, $V=U_p\cap U_i\cap U_j$ is affine, because  $U_p$ and  $U_i\cap U_j$ are affine open subspces of  $U_i$. Therefore, $U_{pq}=V\cap U_q$ is affine, because $V$ and $U_q$ are affine open subspaces of $U_j$. Hence, 
\edemo

All the examples in Examples \ref{Eejemplos} are semiseparated finite spaces.

Finally, let us give some cohomological characterizations of schematic morphisms.

\prop \label{ri=0} Let $X$ be an affine finite space and $Y$ a schematic finite space.  A morphism of ringed spaces $f\colon X\to Y$ is affine iff  $f_*\OO_X$ is quasi-coherent and  $R^if_*\OO_X=0$,  for any $i>0$.

\eprop 

\demo $\Rightarrow)$ Affine finite spaces are acyclic,  then $(R^if_*\OO_X)_y=H^i(f^{-1}(U_y),\OO_X)=0$, for any $i>0$ and any $y\in Y$. Hence, $R^if_*\OO_X=0$, for any $i>0$.

$\Leftarrow)$ Let $U\subseteq Y$ be an affine open  subspace. $H^i(f^{-1}(U),\OO_X)=
H^i(U,f_*\OO_X)\overset{\text{\ref{2}}}=0$, for any $i>0$, then $f^{-1}(U)$ is acyclic, therefore it is  affine by Proposition \ref{12}.\edemo

\prop 
A morphism of ringed spaces $f\colon X\to Y$ between schematic finite spaces is affine iff 
$f_*\OO_X$ is quasi-coherent, $R^if_*\OO_X=0$ for any $i>0$ and there exists an open covering  $\{U_i\}$ of $Y$ such that $f^{-1}(U_i)$ is affine, for any $i$.
\eprop

\demo $\Rightarrow)$ $(R^if_*\OO_X)_y=H^i(f^{-1}(U_y),\OO_X)=0$, for any $i>0$ and any $y\in Y$. Hence, $R^if_*\OO_X=0$, for any $i>0$.

$\Leftarrow)$ The morphisms $f^{-1}(U_i)\to U_i$ are affine, by the previous proposition.
Then, $f$ is affine.

\edemo

\prop Let $f\colon X\to Y$ be a morphism of ringed spaces between schematic finite spaces. Then, $f$ is schematic iff  $\,\Gamma_f\colon X\to X\times_\ZZ Y$, $\Gamma_f(x)=(x,f(x))$ is schematic.\eprop

\demo $\Leftarrow)$ It is easy to check that  $\pi_2\colon X\times_\ZZ Y\to Y$, $\pi_2(x,y)=y$ is schematic. Then, $f$ is schematic because $f= \pi_2\circ \Gamma_f$ and $\pi_2$ and  $\Gamma_f$ are schematic.

$\Rightarrow)$ Let $x\in X$ and  $(x',y)\in X\times_\ZZ Y$ (where $(x,f(x))\leq (x',y)$). $U_{x(x',y)}=U_{x'}\cap U_{xy}$ is affine because it is the intersection of two affine open subsets of the affine finite space $U_x$.
Observe that 
$$\aligned \OO_{x}\otimes_{\OO_{(x,f(x))}} \OO_{(x',y)} & =\OO_{x}\otimes_{\OO_{x}\otimes_\ZZ \OO_{f(x)}} \OO_{x'}\otimes_\ZZ \OO_{y}=\OO_{x'}\otimes_{\OO_{f(x)}} \OO_{y}
\\ & \overset{\text{\ref{K6}}}=\OO_{x'y} =\OO_{x(x',y)}\endaligned$$
Then, $\Gamma_f$ is schematic, by Proposition \ref{K6}.

\edemo

\teor A morphism of ringed spaces $f\colon X\to Y$ between schematic finite spaces is  schematic iff  $R^i{\Gamma_f}_*\OO_X$ is a  quasi-coherent $\OO_{X\times Y}$-module, for any $i\geq 0$.\eteor

\demo  $\Rightarrow)$ By Proposition \ref{K4}, 
$R^i{\Gamma_f}_*\OO_X$ is a quasi-coherent $\OO_{X\times Y}$-module, for any $i\geq 0$.

$\Leftarrow)$ 
$R^i{\Gamma_f}_*\OO_X$ is a quasi-coherent $\OO_{X\times Y}$-module, for any $i\geq 0$. Then,  $$H^0(U_x,\OO_X)\otimes_{f(x)}\OO_y=H^0(U_{xf(x)},\OO_X)\otimes_{\OO_{f(x)}}\OO_{y}=H^0(U_{xf(x)},\OO_X)\otimes_{\OO_{(x,f(x))}}\OO_{(x,y)}=H^0(U_{xy},\OO_X),$$ for any $x$ and $y\geq f(x)$. Therefore,
$\OO_{x}\otimes_{\OO_{f(x)}}\OO_{y}=\OO_{xy}$. Besides,
$$0=H^i(U_{x},\OO_X)\otimes_{\OO_{f(x)}}\OO_{y}=H^i(U_{xf(x)},\OO_X)\otimes_{\OO_{f(x)}}\OO_{y}=H^i(U_{xf(x)},\OO_X)\otimes_{\OO_{(x,f(x))}}\OO_{(x,y)}=H^i(U_{xy},\OO_X),$$
for any  $i>0$. Then, the open subsets $U_{xy}$ are acyclic, hence $f$ is schematic by Proposition \ref{K6}. 

\edemo

\teor Let $f\colon X\to Y$ be a morphism 
of ringed spaces between schematic finite spaces.  Let $x\in X$ and let $f_{U_x}$ be the composite morphism $U_x\iny X\to Y$.
Then, $f$ is  schematic iff  $R^i{f_{U_x}}_*\OO_{U_x}$ is a  quasi-coherent $\OO_{Y}$-module,  for any  $i\geq 0$ and any  $x\in X$.\eteor

\demo  $\Rightarrow)$ If $f$ is schematic,  $f_{U_x}$ is schematic and  $R^i{f_{U_x}}_*\OO_{U_x}$ is a  quasi-coherent $\OO_{Y}$-module,   for any $i\geq 0$,  by Proposition \ref{K4}.

$\Leftarrow)$  $R^i{f_{U_x}}_*\OO_{U_x}$ is a  quasi-coherent $\OO_{Y}$-module. Then,
$$H^0(U_{xf(x)},\OO_X)\otimes_{\OO_{f(x)}}{\OO_{y}}=H^0(U_{xy},\OO_X),$$ for any  $x$ and $y\geq f(x)$. Therefore, $\OO_x\otimes_{\OO_{f(x)}}\OO_y=\OO_{xy}$.
 Besides,
$$0=H^i(U_{x},\OO_X)\otimes_{\OO_{f(x)}}\OO_{y}=H^i(U_{xf(x)},\OO_X)\otimes_{\OO_{f(x)}}\OO_{y}=H^i(U_{xy},\OO_X),\,\text{for any } i>0.$$
Then, the open sets $U_{xy}$ are acyclic. By Proposition \ref{K6}, $f$ is schematic.

\edemo

\section{Quasi-isomorphisms}

\defi A schematic morphism $f\colon X\to Y$ is said to be a quasi-isomorphism  if
\enumera \item f is affine.

\item 
 $f_*\OO_X=\OO_Y$.\eenumera 
 \edefi
 
If $f\colon X\to Y$ is a quasi-isomorphism we shall say that $X$ is quasi-isomorphic to $Y$.
 
\ejems \label{E8.2} \enumera \item  If $X$ is an affine finite space, the morphism
 $X\to (*,\OO(X))$ is a quasi-isomorphism.
 
\item Let $X$ be a schematic finite space and let $\tilde X$  be the Kolmogorov quotient of $X$.
The quotient morphism $\pi\colon X\to \tilde X$ is a quasi-isomorphism (see Example \ref{Kolmogorov}.3.).
 
\item If $X$ is a schematic finite $T_0$-topological space, 
 $X_M \hookrightarrow X$ is a quasi-isomorphism.

\item Let $f\colon X\to Y$ be a schematic morphism. $(Y,f_*\OO_X)$ is a schematic finite space
by Example \ref{E5.8}.
Let us prove that the obvious ringed morphism $f'\colon X\to (Y,f_*\OO_X)$, $f'(x)=f(x)$ is schematic. By Theorem \ref{K7}, the morphism
$$\OO_x\otimes_{(f_*\OO_X)_{f(x)}} (f_*\OO_X)_y=
\OO_x\otimes_{(f_*\OO_X)_{f(x)}} (f_*\OO_X)_{f(x)}\otimes_{\OO_{f(x)}}\OO_y=\OO_x\otimes_{\OO_{f(x)}}\OO_y
\to \prod_{z\in U_{xy}}\OO_z$$
is surjective on spectra. Again by Theorem \ref{K7}, $f'$ is schematic. We have the obvious commutative diagram
$$\xymatrix{ X \ar[rr]^-f \ar[rd]^-{f'} & & Y \\ & (Y,f_*\OO_X) \ar[ru]^-{\Id} &}$$
$\Id$ is affine. If $f$ is affine, then $f'$ is a quasi-isomorphism.

\eenumera

\eejems

\defi Let $f\colon X\to Y$ be a schematic morphism. We shall say that $f$ is  flat
if the morphism $\OO_{Y,f(x)}\to \OO_{X,x}$ is flat, for any $x\in X$.
We shall say that $f$ is  faithfully flat
if the morphism $\OO_{Y,y}\to \prod_{x\in f^{-1}(U_y)}\OO_{X,x}$ is faithfully flat, for any $y\in Y$.
\edefi 

If $\{U_i\}$ is an open covering of $X$, the natural morphism $\coprod_i U_i\to X$ is faithfully flat.
 
\obse \label{O12.3} Quasi-isomorphisms are faithfully flat morphisms: Given $y\in Y$, $f^{-1}(U_y)$ is affine, then the morphism 
$$\OO_y=(f_*\OO_X)_y=\OO_X(f^{-1}(U_y))\to \prod_{x\in f^{-1}(U_y)} \OO_x$$ is faithfully flat.
\eobse

\prop Let  $X$ be a schematic finite space, $\U=\{U_1,\ldots,U_n\}$ a minimal  affine open covering of $X$ and  $Y$ the ringed  finite space associated with $\U$. Then, $Y$ is a schematic finite space and the quotient morphism  $\pi\colon X\to Y$ is a quasi-isomorphism.\eprop

\demo  $Y=\{y_1,\ldots,y_n\}$, where $\pi^{-1}(U_{y_i})=U_i$. Recall that $\pi_*\OO_X=\OO_Y$. 
Let $y_1\leq y_2$, then
$$\OO_{y_1}=\OO_X(U_1)\to \OO_X(U_2)=\OO_{y_2}$$
is a flat morphism by Lemma \ref{flat}. Let $y_i,y_j\geq y_k$, then
$$\aligned \OO_{y_i}\otimes_{\OO_{y_k}}\OO_{y_j} & =\OO_X(U_i)\otimes_{\OO_X(U_k)}\OO_X(U_j)\overset{\text{\ref{11}}}=
\OO_X(U_i\cap U_j)=\OO_X(\pi^{-1}(U_{y_i}\cap U_{y_j}))\\ &=\OO_Y(U_{y_i}\cap U_{y_j}) =\OO_{y_iy_j}.\endaligned$$
$U_i\cap U_j$ is an affine finite space by
Proposition \ref{inter}. $U_i\cap U_j=\cupa{U_k\subset U_i\cap U_j} U_z$. The morphisms $\OO(U_k)\to \proda{x\in U_k}\OO_x$, $\OO(U_i\cap U_j)\to \proda{U_k\subset U_i\cap U_j,x\in U_k} \OO_x$  are faithfully flat. Then, the morphism
$\OO(U_i\cap U_j)\to \proda{U_k\subset U_i\cap U_j} \OO(U_k)$ is faithfully flat. Therefore, the morphism
$\OO_{y_iy_j}\to \proda{y_k\in U_{y_iy_j}} \OO_{y_k}$ is faithfully flat. Then, $Y$ is a schematic finite space.
%
%
%Observe that $R^i\pi_*\OO_X=0$, for any $i>0$, because $(R^i\pi_*\OO_X)_{y_i}=H^i(U_i,\OO_X)=0$.
%Then,
%$$H^i(U_{y_j},\OO_Y)=H^i(U_j,\OO_X)=0,\text{ for any } i>0 \text{ and } j.$$
%Then, $U_{y_j}$ are acyclic. Likewise, $U_{y_iy_j}$ are acyclic, for any $y_i,y_j\geq y_k$. By Theorem \ref{afin'}, 
%$U_{y_k}$ is an affine finite space. Then, $Y$ is a schematic finite space. 

Finally, $\pi$ is affine  by Proposition \ref{hmm}.

\edemo

\prop \label{P9.6} The composition of quasi-isomorphisms is a quasi-isomorphism. \eprop

\hskip-0.65cm\colorbox{white}{\,\begin{minipage}{15.15cm}
\teor \label{T12.7} Let $f\colon X\to Y$ be a schematic morphism. The functors $$f_*\colon\,{\bf Qc\text{-}Mod}_X\,\to\,{\bf Qc\text{-}Mod}_Y\,\text{ and }\, f^*\colon\,{\bf Qc\text{-}Mod}_Y\,\to\,{\bf Qc\text{-}Mod}_X\,$$ are mutually inverse (i.e., the natural morphisms  $\mathcal M\to f_*f^*\M$, $f^*f_*\Nc\to \Nc$ are isomorphisms) iff $f$ is a quasi-isomorphism.\eteor\end{minipage}}

\demo $\Leftarrow)$ The morphism $f^*f_*\M\to \M$ is an isomorphism: It is a local property on $Y$. We can suppose that $Y$ is an affine finite space (then $X$ is affine).
Consider a free presentation of $\M$,  $\oplus_I \OO_X\to
\oplus_J\OO_X\to \M\to 0$. Taking $f_*$, which is an exact functor because $R^1f_*=0$, one has the exact sequence of morphisms
$\oplus_I \OO_Y\to
\oplus_J\OO_Y\to f_*\M\to 0$. Taking $f^*$, one has the exact sequence of morphisms  $\oplus_I \OO_X\to
\oplus_J\OO_X\to f^*f_*\M\to 0$, then $f^*f_*\M=\M$.

Likewise,  $f_*f^*\Nc=\Nc$.

$\Rightarrow)$ $\OO_Y=f_*f^*\OO_Y=f_*\OO_X$.  
Obviously, $f_*$ is an exact functor. By Theorem \ref{Tafex},
$f$ is affine.

\edemo

\coro \label{superafin} Let  $f\colon X\to Y$ be a  quasi-isomorphism.
$Y$ is affine iff  $X$ is affine.\ecoro

\demo $\Leftarrow)$ 
For any quasi-coherent $\OO_Y$-module $\Nc$, 
$\Nc=f_*f^*\Nc$, then
$$H^1(Y,\Nc)=H^1(X,f^*\Nc)=0.$$ By the Serre Theorem
$Y$ is affine.\edemo

\coro \label{qc-iso} Let  $X\overset f\to Y\overset g\to Z$ be schematic morphisms and assume $g\circ f$ is a  quasi-isomorphism. Then,

\enumera
\item If $g$ is a quasi-isomorphism, then  $f$ is a quasi-isomorphism.

\item If $f$  is a quasi-isomorphism, then  $g$ is a quasi-isomorphism. \eenumera\ecoro

\demo 1. Considering the diagram
 $$\xymatrix{\,{\bf Qc\text{-}Mod}_X\, \ar[r]^-{f_*} \ar@/^{5mm}/[rr]^-{(g\circ f)_*} &\,{\bf Qc\text{-}Mod}_Y\, \ar[r]^-{g_*} \ar@<1ex>[l]^-{f^*} &\,{\bf Qc\text{-}Mod}_Z\, \ar@<1ex>[l]^-{g^*} \ar@/^{9mm}/[ll]^-{(g\circ f)^*}}$$
it is easy to prove that $f_*$ and $f^*$ are mutually inverse functors.

2. Proceed likewise.
\edemo

\coro Let $X$ and $Y$ be affine finite spaces. Then, a schematic morphism $f\colon X\to Y$ is a quasi-isomorphism iff $\OO_Y(Y)=\OO_X(X)$.\ecoro

\demo $\Leftarrow)$  Observe that the diagram
$$\xymatrix{X \ar[r]^-f \ar[d] & Y\ar[d]\\ (*,\OO_X(X))
\ar[r] & (*,\OO_Y(Y))}$$
is commutative, $(X,\OO_X)$ is quasi-isomorphic to $(*,\OO_X(X))$ and $(Y,\OO_Y)$ is quasi-isomorphic to  $(*,\OO_Y(Y))$.

\edemo

\coro \label{proposiciontonta} Let  $f\colon X\to Y$ be a schematic morphism and  let $f_{M}\colon X_M\to Y_M$ be the induced morphism. Then, $f$ is a quasi-isomorphism iff $f_{M}$ is a quasi-isomorphism.

\ecoro

\demo It is an immediate consequence of Corollary \ref{qc-iso}.\edemo

\coro \label{C11.12} Let $f\colon X\to X'$ be a quasi-isomorphism , $X''$ a schematic finite space and $g\colon X'\to X''$ a morphism of ringed spaces. Then, $g$ is schematic (resp. affine)  iff $g\circ f$ is schematic (resp. affine) \ecoro

\demo Recall Theorem \ref{Tcohsch} (resp. Theorem \ref{Tafex}).

\edemo

\prop  \label{P11.13} Let $f\colon X\to X'$ be an affine  morphism of schematic spaces. Assume that $X'$ is $T_0$. Let $\U:=\{f^{-1}(U_{x'})\}_{x'\in X'}$, let $X/\!\sim$ be the schematic space associated with the open covering $\U$, and $\pi\colon X\to X/\!\sim$  the quotient morphism.  The morphism
$f'\colon X/\!\sim\,\to X'$, $f'([x])=f(x)$, induced by $f$, is  affine
and $Y$ is homeomorphic to $\Ima f$.

\eprop

\demo  By Corollary \ref{C11.12}, we only have to prove that $Y$ is homeomorphic to $\Ima f$.
The morphism $f'\colon X/\!\sim\,\to \Ima f$ is clearly bijective  and continuous. Given
$[x],[x']\in X/\!\sim$, if $f'([x])\leq f'([x'])$, then $f(x)\leq f(x')$
and $U_{f(x')}\subseteq U_{f(x)}$. Hence,
$f^{-1}(U_{f(x')})\subseteq f^{-1}(U_{f(x)})$ and $U_{[x']}\subseteq U_{[x]}$. Therefore, $[x]\leq [x']$. That is, $f'$ is a homeomorphism.

\edemo

\lema \label{lematonto} Let $h\colon X\to Y$ be a quasi-isomorphism.  Then, $Y\backslash h(X)$ is a set of removable points of $Y$.\elema

\demo  Let $y\in Y\backslash h(X)$. Since $h^{-1}(U_y)$ is affine, the morphism
$$\OO_Y(U_y)=\OO_{X}(h^{-1}(U_y))\to \prod_{x'\in h^{-1}(U_y)} \OO_{X,x'}$$ is faithfully flat. This morphism factors through the morphism $$\OO_Y(U_y)\to \prod_{h(x')\in U_y} \OO_{Y,h(x')},$$ then this last morphism is faithfully flat. Hence, $y$ is a removable point of $Y$.
\edemo

\hskip-0.65cm\colorbox{white}{\,\begin{minipage}{15.15cm}
\teor \label{comosonqc} Let  $f\colon X\to Y$ be a quasi-isomorphism. Assume that $Y$ is $T_0$. Consider the affine open covering of $X$, $\{f^{-1}(U_y)\}_{y\in Y}$ and let $X/\!\sim$ be the associated schematic finite space. Then, $f$ is the composition of the quotient morphism  $\pi\colon X\to X/\!\sim$ and an isomorphism $f'\colon X/\!\sim\,\to Y\backslash P$, $f'([x])=f(x)$,  where $P$ is a set of removable points of $Y$. 

Therefore, if $f\colon X\to Y$ is a quasi-isomorphism and   $Y$ is minimal,
$f$  is the composition of the quotient morphism  $X\to X/\!\sim$ and an isomorphism $X/\!\sim\,\simeq Y$.
\eteor\end{minipage}}

\demo By Lemma \ref{lematonto}, we can suppose that $f$ is surjective. By Proposition \ref{P11.13}, $f'\colon Y'\to Y$ is a homeomorphism and it is affine.
Finally, $\OO_{Y,f'([x])}=(f_*\OO_{X})_{f(x)}=\OO_{Y',[x]}$, for any $[x]$.
\edemo

\section{Change of base and flat schematic morphisms}

\prop[(\cite{KG} 5.27)\,] \label{XxY} Let $X$, $X'$ and  $Y$ schematic finite spaces and  $f\colon X\to Y$ and $f'\colon X'\to Y$ schematic morphisms. Then,   

\enumera 
\item $X\times_Y X'$ is a schematic finite space.

\item If $X$, $X'$ and $Y$ are affine, then   $\OO(X\times_YX')=\OO(X)\otimes_{\OO(Y)}\OO(X')$ and $X\times_YX'$ is affine.

\item Given a commutative diagram of schematic morphisms
$$\xymatrix@R=8pt{ U \ar[rd]^-g \ar[rr]^-h & & X \ar[rd]^-f & \\ & V \ar[rr] & & Y\\ U' \ar[ru]_-{g'} \ar[rr]_-{h'}& & X' \ar[ru]_-{f'}} $$
the morphism  $h\times h'\colon U\times_VU'\to X\times_Y X'$, $h\times h'(u,u'):=(h(u),h'(u'))$ is schematic.

\item $\pi\colon X\times_Y X'\to X$, $\pi(x,x')=x$, is schematic.

\item If $h\colon X\to X'$ is a schematic $Y$-morphism, then $\Gamma_h\colon X\to X\times_Y X'$, $\Gamma_h(x):=(x,h(x))$ is schematic.

\item The diagonal morphism  $\delta\colon X\to X\times_Y X$, $\delta(x)=(x,x)$ is schematic.

\eenumera
\eprop

\demo 1. We only need to prove 2.

2. $X\times_{\OO(Y)} X'$ is an affine schematic space  and
$\OO(X\times_{\OO(Y)} X')=\OO(X)\otimes_{\OO(Y)} \OO(X')$, by Proposition \ref{ProdAf}. Let $(x,x')\in X\times_{\OO(Y)} X'$. 
If $f(x)=f'(x')$, then $$\OO_x\otimes_{\OO(Y)}\OO_{x'}=
\OO_{x}\otimes_{\OO_{f(x)}}(\OO_{f(x)}\otimes_{\OO(Y)}\OO_{f(x)})\otimes_{\OO_{f(x)}}\OO_{x'}\overset{\text{\ref{11}}}=
\OO_{x}\otimes_{\OO_{f(x)}}\OO_{f(x)}\otimes_{\OO_{f(x)}}\OO_{x'}=\OO_x\otimes_{\OO_{f(x)}}\OO_{x'}$$
If $f(x)\neq f'(x')$, then $(x,x')\in X\times_{\OO(Y)}X'$ is a removable point: 
Consider the morphism
$f_x\colon U_x \to Y,$  $f_x(z):=f(z)$. Then, $\OO_{xy}=(f_{x*}\OO_{U_x})_{y}=(f_{x*}\OO_{U_x})(Y)\otimes_{\OO(Y)} \OO_{y}=\OO_x\otimes_{\OO(Y)} \OO_{y}$.
Observe that
$U_{x'f(x)}\underset\neq\subset X'$ and  it is affine since $f'_{x'}\colon U_{x'} \to Y,$  $f'_{x'}(z):=f'(z)$ is affine. The morphism
$$\OO_x\otimes_{\OO(Y)} \OO_{x'} =\OO_{xf(x)} \otimes_{\OO(Y)} \OO_{x'}=\OO_x\otimes_{\OO(Y)} \OO_{f(x)}
\otimes_{\OO(Y)} \OO_{x'}=\OO_x\otimes_{\OO(Y)} \OO_{x'f(x)}\to \proda{z\in U_{x'f(x)} } \OO_x \otimes_{\OO(Y)} \OO_z $$
is faithfully flat. Therefore, $(x,x')$ is removable.
In conclusion, $X\times_Y X'=(X\times_{\OO(Y)} X')\backslash \{$A set of removable points$\}$, then $X\times_Y X'$ is affine and
$\OO(X\times_Y X')=\OO(X\times_{\OO(Y)} X')=\OO(X)\otimes_{\OO(Y)}\OO(X')$.

3. By Proposition \ref{K6},  given  $(u,u')\in U\times_VU'$,
$(x,x')\in X\times_Y X'$ (where $(x,x')\geq  (h(u), h'(u'))$),  we have to prove that $U_{(u,u')(x,x')}$ is acyclic and
$\OO_{(u,u')(h(u),h(u'))}\otimes_{\OO_{(h(u),h(u'))}} \OO_{(x,x')}=\OO_{(u,u')( x,x')}$:
$U_{(u,u')(x,x')} =U_{ux}\times_{U_{g(u)f(x)}} U_{u'x'}$, which is affine  (then acyclic) and 
$$\aligned  & \OO_{(u,u')} \otimes_{\OO_{(h(u),h(u'))}} \OO_{(x,x')} 
=(\OO_{u}\otimes_{\OO_{g(u)}} \OO_{u'})\otimes_{\OO_{h(u)}\otimes_{\OO_{f(h(u))}}\OO_{h'(u')}} (\OO_{x}\otimes_{\OO_{f(x)}}\OO_{x'}) 
\\ & = (\OO_u\otimes_{\OO_{h(u)}}\OO_x)\otimes_{\OO_{g(u)}\otimes_{\OO_{f(h(u))}}\OO_{f(x)}}
 (\OO_{u'}\otimes_{\OO_{h'(u')}}\OO_{x'})
 \overset{\text{\ref{K6}}}=\OO_{u x}\otimes_{\OO_{g(u)f(x)}} \OO_{u' x'}
 =\OO_{(u,u')(x,x')}.
\endaligned$$

4.,5. and  6. are particular cases of 3.

\edemo

Obviously, $$\Hom_Y(Z,X\times_Y X')=\Hom_Y(Z,X)\times \Hom_Y(Z,X')$$
for any schematic finite $Y$-spaces $Z,X,X'$.

\prop \label{qc1} Affine morphisms and  quasi-isomorphisms are stable by base change.
\eprop

\demo Let $f\colon X\to Y$ be an affine morphism and  $Y'\to Y$ a schematic morphism. In order to prove that the schematic morphism $X\times_YY'\to Y'$ is affine, it is sufficient to prove that  $X\times_YU_{y'}$ is affine, for any $y'\in Y'$.  Observe that $X\times_YU_{y'}=f^{-1}(U_{f(y')})\times_{U_{f(y')} }U_{y'}$, which is affine because  $f^{-1}(U_{f(y')})$ is affine.

Let $f$ be a quasi-isomorphism. We only have to prove that $\OO(X\times_YU_{y'})=\OO_{y'}$:

$$\OO(X\times_YU_{y'})=\OO(f^{-1}(U_{f(y')}))\times_{U_{f(y')}}U_{y'})=
\OO(f^{-1}(U_{f(y')}))\otimes_{\OO_{f(y')}} \OO_{y'} ={\OO_{f(y')}}\otimes_{\OO_{f(y')}} \OO_{y'}=\OO_{y'}.$$

\edemo

\lema \label{L11.3} Let $f\colon X\to Y$ and 
$g\colon Y'\to Y$ be schematic morphisms and let 
$g'\colon
X\times_Y Y'\to X$ be defined by $g'(x,y'):=x$. Let $\M$ be a quasi-coherent $\OO_X$-module. If $X,Y$ and $Y'$ are affine, then
$$\Gamma(X\times_Y Y',g'^*\M)=\Gamma(X,\M)\otimes_{\OO(Y)}\OO(Y').$$
\elema

\demo Consider an exact sequence of $\OO_X$-modules
$\oplusa{I}\, \OO_X\to \oplusa{J}\,\OO_X\to \M\to 0$.

1. Taking $g'^*$, $\oplusa{I}\, \OO_{X\times_Y Y'}\to \oplusa{J}\,\OO_{X\times_Y Y'}\to g'^*\M\to 0$ is exact. Taking sections, the sequence 
$\oplusa{I} \,\OO(X\times_Y Y')\to \oplusa{J}\,\OO(X\times_Y Y')\to g'^*\M(X\times_Y Y')\to 0$ is exact. By Proposition \ref{XxY}, $\OO(X\times_Y Y')=\OO(X)\otimes_{\OO(Y)}\OO(Y')$. Hence, the sequence of morphisms
$$\oplusa{I}\, \OO(X)\otimes_{\OO(Y)}\OO(Y')\to \oplusa{J}\,\OO(X)\otimes_{\OO(Y)}\OO(Y')\to g'^*\M(X\times_Y Y')\to 0$$ is exact.

2. The sequence of $\OO(X)$-modules
$\oplusa{I} \,\OO(X)\to \oplusa{J}\,\OO(X)\to \M(X)\to 0$ is exact.
Hence, the sequence $\oplusa{I} \,\OO(X)\otimes_{\OO(Y)}\OO(Y')\to \oplusa{J} \,\OO(X)\otimes_{\OO(Y)}\OO(Y')\to \M(X)\otimes_{\OO(Y)} \OO(Y')\to 0$ is exact.

3. Therefore, the natural morphism $\Gamma(X,\M)\otimes_{\OO(Y)}\OO(Y')\to \Gamma(X\times_Y Y',g'^*\M)$ is an isomorphism.

\edemo

\teor \label{Tcbp} Let $f\colon X\to Y$ be a schematic morphism and
$g\colon Y'\to Y$ a flat schematic morphism. Denote 
$f'\colon X\times_Y Y'\to Y'$, $f'(x,y')=y'$, $g'\colon
X\times_Y Y'\to X$, $g'(x,y')=x$ the induced morphisms.
Then, the natural morphism
$$g^*R^if_*\mathcal M\to R^if'_*(g'^*\mathcal M)$$
is an isomorphism.\eteor 

\demo We have to prove that the morphism is an isomorphism on stalks at $z$, for any $z\in Y'$. That is to say, we have to prove that the morphism
$$H^i(f^{-1}(U_{g(z)}), \M)\otimes_{\OO_{g(z)}} \OO_z\to 
H^i(f^{-1}(U_{g(z)})\times_{U_{g(z)}} U_z, g'^*\M)$$
is an isomorphism. We can suppose that $Y'=U_z$, $Y=U_{g(z)}$ and $X=f^{-1}(U_{g(z)})$. Then, $g'$ is an affine morphism and
$H^i(X\times_Y Y', g'^*\M)=H^i(X,g'_*g'^*\M)$. By Lemma \ref{L11.3}, 
$$(g'_*g'^*\M)_x=\Gamma(U_x\times_{Y}Y',g'^*\M)=
\M_x\otimes_{\OO(Y)}\OO(Y').$$
Hence, $C^\punto(g'_*g'^*\M)=(C^\punto\M)\otimes_{\OO(Y)}\OO(Y')$ and $H^i(X,g'_*g'^*\M)=H^i(X,\M)\otimes_{\OO(Y)}\OO(Y')$. That is, 
$$H^i(X\times_Y Y', g'^*\M)=H^i(X,g'_*g'^*\M)=H^i(X,\M)\otimes_{\OO(Y)}\OO(Y').$$
\edemo

\prop \label{Pff} Let $X$ and $Y$ be  schematic finite spaces and 
$f\colon X\to Y$ a schematic morphism.  Then, $f$ is flat (resp. faithfully flat) iff the functor
$$f^*\colon\,{\bf Qc\text{-}Mod}_Y\,\to\,{\bf Qc\text{-}Mod}_X, \,\,\mathcal M \functor f^*\mathcal M$$  is exact (resp. faithfully exact). 
\eprop

\demo Obviously, if $f$ is flat then $f^*$ is exact. If $f$ is faithfully flat then $f^*$ is faithfully exact: We only have  to check that $\M=0$ if $f^*\M=0$. Given $y\in Y$, the morphism $\OO_y\iny \prod_{x\in f^{-1}(U_y)} \OO_x$ is faithfully flat. Tensoring by $\M_y\otimes_{\OO_y}$ one has the injective morphism
$$\aligned \M_y\iny \prod_{x\in f^{-1}(U_y)} \M_y\otimes_{\OO_y}\OO_x
& =\prod_{x\in f^{-1}(U_y)} \M_y\otimes_{\OO_y}\OO_{f(x)} \otimes_{\OO_{f(x)} }\OO_x=\prod_{x\in f^{-1}(U_y)} \M_{f(x)} \otimes_{\OO_{f(x)} }\OO_x\\ & =\prod_{x\in f^{-1}(U_y)} (f^*\M)_x=0.\endaligned$$
Hence, $\M=0$.

 If  $f^*$ is exact $f$ is flat: Let $f(x)\in Y$. Consider an ideal $I\subseteq \OO_{f(x)}$ and let $\tilde I\subseteq \OO_{U_{f(x)}}$ be the quasi-coherent $\OO_{U_{f(x)}}$-module associated with $I$. 
Consider the inclusion morphism $i\colon U_{f(x)}\iny Y$. The morphism
$i_*\tilde I\iny i_*\OO_{U_{f(x)}}$ is injective, then $f^*i_*\tilde I\iny f^* i_*\OO_{U_{f(x)}}$ is injective. Hence, $I\otimes_{\OO_{f(x)}} \OO_x=(f^*i_*\tilde I)_x\iny (f^* i_*\OO_{U_{f(x)}})_x=\OO_x$ is injective. Therefore,
the morphism $\OO_{f(x)}\to \OO_x$ is flat.
Finally, 
if  $f^*$ is faithfully exact $f$ is faithfully  flat: Let $y\in Y$ be a maximal point of $Y$, if it exists,  such that the flat morphism $\OO_y\to \prod_{x\in f^{-1}(U_y)} \OO_x$ is not faifhfully flat. Then, there exists an ideal $I\underset\neq\subset \OO_y$
such that $I\cdot \prod_{x\in f^{-1}(U_y)} \OO_x=\prod_{x\in f^{-1}(U_y)} \OO_x$. Let $\M$ be the quasi-coherent $\OO_{U_y}$-module associated with  $\OO_y/I$ and let  $f'\colon f^{-1}(U_y)\to U_y$ be the morphism defined by $f'(x):=f(x)$.
Obviously, $f'^*\M=0$. Let $i\colon U_y\iny Y$ and $i'\colon
f^{-1}(U_y)\iny X$ be the inclusion morphisms. By Theorem \ref{Tcbp},
$0=i'_*f'^*\M=f^*i_*\M$, then $i_*\M=0$. Hence, $0=i^*i_*\M=\M$ which is contradictory.

\edemo

\prop Let $f\colon X\to Y$ be a  schematic morphism. Then, $f$ is a quasi-isomorphism iff it is faithfully flat and the natural morphism $f^*f_*\M\to \M$ is an isomorphism, for any quasi-coherent $\OO_X$-module $\M$.
\eprop

\demo $\Rightarrow)$ It is Remark \ref{O12.3} and Theorem \ref{T12.7}.

$\Leftarrow)$ $f_*$ is an exact functor since $f^*$ is faithfully exact  and $\Id=f^*f_*$. By Theorem \ref{Tafex}, $f$ is affine. Finally, the morphism $\OO_Y\to f_*\OO_X$ is an isomorphism, since $f^*\OO_Y=\OO_X\to f^*f_*\OO_X$ is an isomorphism.

\edemo

\prop Let $f\colon X\to Y$ be a schematic morphism. Then, $f$ is faithfully flat iff
$f_M\colon X_M\to Y_M$ is surjective and flat.
\eprop

\demo Let $g\colon X'\to X$ be a quasi-isomorphism. Then, $f$ is faithfully flat iff $f\circ g$ is faithfully flat, since $f^*$ is faithfully exact iff $g^*\circ f^*$ is faithfully exact.
Let $g\colon Y\to Y'$ 
be a quasi-isomorphism. Likewise, $f$ is faithfully flat iff $g\circ f$ is faithfully flat.

Therefore, we have to prove that $f_M$ is a faithfully flat iff is surjective and flat. The converse implication is obvious. Let us prove the direct implication. If $f$ is not surjective, let $y\in Y_M$ be maximal satisfying
$f^{-1}(y)=\emptyset$. Consider the commutative diagram of obvious morphisms
$$\xymatrix{\OO_y \ar[r]^-{i_1} \ar[d]^-{i_2}  & \prod_{f_M(x)\geq y} \OO_x \ar@{=}[r] &
\prod_{f_M(x)> y} \OO_x \ar[dl]^-{i_4}\\ \prod_{y'\geq y} \OO_{y'} \ar[r]^-{i_3} & 
\prod_{y'\geq y} \prod_{f_M(x)\geq y'} \OO_x & }$$
The morphisms $i_1$ and $i_3$ are faithfully flat since $f_M$ is faithfully flat, $i_4$ is obviously faithfully flat, hence $i_2$ is faithfully flat and $y$ is a removable point, which is contradictory.

\edemo 

\prop Let $f\colon X\to Y$ be a schematic morphism and
$g\colon Y'\to Y$ a  faithfully flat  schematic morphism. Let
$f'\colon X\times_Y Y'\to Y'$ be the morphism defined by $f'(x,y')=y'$. 
Then, \enumera \item $f$ is affine   iff $f'$ is affine.

\item $f$ is a quasi-isomorphism iff $f'$ is a quasi-isomorphism.\eenumera

\eprop

\demo We can suppose that $X$,$Y$ and $Y'$ are minimal schematic spaces. 
The morphism $g'\colon  X\times_Y Y'\to X$, $g'(x,y'):=x$ is faithfully flat since it is flat and surjective.

1. $\Leftarrow)$ The functor $g^*f_*=f'_*g'^*$ is exact since $f'_*$ and $g'^*$ are exact.
Hence, $f_*$ is exact since $g^*$ is faithfully exact and $f$ is affine.

2.  $\Leftarrow)$ We only have to prove that the morphism $\OO_Y\to f_*\OO_X$ is an isomorphism.
Taking $g^*$, we obtain the isomorphism $\OO_X\to g^*f_*\OO_X=f'_*g'^*\OO_X=f'_*\OO_{X\times_Y Y'}=\OO_X$. Hence, $\OO_Y\to f_*\OO_X$ is an isomorphism.

\edemo

\section{Quasi-open immersions}

 \defi We shall  say that a schematic morphism $f\colon X\to Y$
 is a quasi-open immersion if it is flat and the diagonal morphism $X\to X\times_YX$ is a quasi-isomorphism.\edefi

\ejem If $X$ is a schematic finite space and $U\subseteq X$ an open subset, then the inclusion morphism $U\hookrightarrow X$ is a quasi-open immersion.
\eejem

\prop \label{P3.13} If $f\colon X\to Y$ is a quasi-isomorphism, then it is a quasi-open immersion\eprop

\demo Quasi-isomorphisms are faithfully flat morphisms, by Observation \ref{O12.3}. 
The morphism $X\times_Y X\to X$ is a quasi-isomorphism by Proposition \ref{qc1}. The composite morphism
$X\to X\times_Y X\to X$ is the identity morphism, then  $X\to X\times_Y X$ is a quasi-isomorphism, by Corollary \ref{qc-iso}.

\edemo

\prop \label{P13.4} If $f\colon X\to Y$ is a quasi-open immersion and $Y'\to Y$ a schematic morphism, then $X\times_YY'\to Y'$ is a quasi-open immersion.\eprop

\demo The morphism $X\to Y$ is flat. Taking $\times_YY'$, the morphism $X\times_YY'\to Y'$ is flat.
The morphism $X\to X\times_Y X$ is a quasi-isomorphism. Taking $\times_YY'$, the morphism $X\times_YY'\to (X\times_Y Y')\times_{Y'} (X\times_YY')$ is a quasi-isomorphism, by Proposition \ref{qc1}. Hence,  $X\times_YY'\to Y'$ is a quasi-open immersion.

\edemo

\prop Let $f\colon X\to Y$ be a schematic morphism. Let $Y'\to Y$ be  a faithfully flat schematic morphism
and $f'\colon X\times_Y Y'\to Y'$, $f'(x,y'):=f(x)$ the induced morphism. Then $f$ is a quasi-open immersion iff $f'$ is a quasi-open immersion.\eprop

\prop The composition of two quasi-open immersions is a quasi-open immersion.\eprop

 \demo The composition of two flat morphisms is flat. 
 Let $f\colon X\to Y$, $g\colon Y\to Z$ be quasi-open immersions.
 Consider the commutative diagram
$$\xymatrix{X \ar[r]^-{\delta_X} & X\times_YX \ar[r]^-{\Id\times \Id} \ar[d] &
X\times_ZX\ar[d]^-{f\times f}\\ & Y \ar[r]_-{\delta_Y}  
& Y\times_ZY}$$ 
Observe that $X\times_Y X=Y\times_{Y\times_ZY} (X\times_Z X)$ and $\delta_Y$ is a quasi-isomorphism. Then, $\Id\times \Id$ is a quasi-isomorphism  by Proposition
 \ref{P13.4}. The morphism $(\Id\times \Id)\circ \delta_X$ is a quasi-isomorphism by Proposition \ref{P9.6}, hence $g\circ f$ is a quasi-open immersion.

\edemo 

\prop  Let $f\colon X\to Y$, $g\colon Y\to Z$ be schematic morphisms and suppose $g\circ f$ is a quasi-open immersion.
\begin{enumerate}
\item If $g$ is a quasi-open immersion, then $f$ is a quasi-open immersion.

\item If $f$ is a quasi-isomorphism, then $g$ is a quasi-open immersion.

\end{enumerate}
\eprop

\demo 1.  Consider the commutative diagram
$$\xymatrix{X \ar[r]^-{\delta_X} & X\times_YX \ar[r]^-{\Id\times \Id} \ar[d] &
X\times_ZX\ar[d]^-{f\times f}\\ & Y \ar[r]_-{\delta_Y}  
& Y\times_ZY}$$ 
$\Id\times \Id$ is a quasi-isomorphism since $\delta_Y$ is a quasi-isomorphism. $(\Id\times \Id)\circ \delta_X$ is a quasi-isomorphism, since $g\circ f$ is a quasi-open immersion. Hence, $\delta_X$ is a quasi-isomorphism, by Corollary \ref{qc-iso}, that is, $f$ is a quasi-open immersion.

2. The obvious morphism $X\times_ZX\to Y\times_ZY$ is a quasi-isomorphism, since is the composition of the quasi-isomorphisms
$X\times_ZX\to X\times_ZY$, $X\times_ZY\to Y\times_ZY$.
Consider the commutative diagram
$$\xymatrix{X \ar[r]^-{\delta_X} \ar[d]_-f &
X\times_ZX\ar[d]^-{f\times f}\\ Y \ar[r]_-{\delta_Y}  
& Y\times_ZY}$$
Then, $\delta_Y$ is a quasi-isomorphism since $f,\delta_X$ and $f\times f$ are quasi-isomorphisms. That is, $g$ is a quasi-open immersion.

\edemo

\defi Let $X$ and $Y$ be ringed finite spaces and $f\colon X\to Y$ a morphism of ringed spaces.
$C(f):=X\coprod Y$ is a finite ringed space as follows: the order relation on $X$ and on $Y$ is the pre-stablished order relation, and given $x\in X$ and $y\in Y$ we shall say that $x>y$ if $f(x)\geq  y$; $\OO_{C(f),x}:=\OO_{X,x}$ for any $x\in X$,  $\OO_{C(f),y}:=\OO_{Y,y}$ for any $y\in Y$;  the morphisms
between the stalks of $\OO_{C(f)}$  are defined in the obvious way.
\edefi 

Observe that $X$ is an open subset of $C(f)$ and $F\colon C(f)\to Y$, $F(x):= f(x)$, for any $x\in X$ and $F(y):= y$, for any $y\in Y$ is a morphism of ringed spaces. $F_*\OO_{C(f)}=\OO_Y$ because 
$$(F_*\OO_{C(f)})_y=\OO_{C(f),y}=\OO_{Y,y}$$
Besides, $R^iF_*\OO_{C(f)}=0$ for any $i>0$, because 
$$(R^iF_*\OO_{C(f)})_y=H^i(U_y,\OO_{C(f)})=0.$$

\teor Let $f\colon X\to Y$ be a schematic morphism. 
Then, $f$ is a quasi-open immersion iff $C(f)=X\coprod Y$ is a schematic finite space. If $f$ is a quasi-open immersion,
then it is the composition of the open inclusion $X\hookrightarrow C(f)$ and the quasi-isomorphism $F\colon C(f)\to Y$. 

%\begin{enumerate}
%
%
%\item $F\colon C(f)\to Y$ is a quasi-isomorphism.
%
%\item $\Ima f\subset C(f)$ is a subset of removable points. Then, the composite morphism
%$$X\coprod (Y-\Ima f) \iny C(f)\overset F\to Y$$
%is a quasi-isomorphism.
%
%\end{enumerate} 
\eteor

\demo  $\Rightarrow)$ Given, $x\geq x'\in X\subset C(f)$, the morphism
$$\OO_{C(f), x'}=\OO_{X,x'}\to \OO_{X,x}=\OO_{C(f),x}$$
is flat. Given, $y\leq y'\in Y\subset C(f)$,  the morphism
$$\OO_{C(f), y}=\OO_{Y,y}\to \OO_{Y,y'}=\OO_{C(f),y'}$$
is flat. Given $x\in X\subset C(f)$ and $f(x)\geq y \in Y\subset C(f)$, the morphism
$$\OO_{C(f),y}\to \OO_{C(f),f(x)}=\OO_{Y,f(x)}\to \OO_{X,x}=\OO_{C(f),x}$$
is flat.

Given $c\in C(f)$, we shall denote $\tilde U_c:=\{z\in C(f)\colon z\geq c\}$.
We have to prove that $\tilde U_c$ is affine. Recall Theorem \ref{afin'}.

a.  If $c=x\in X$, then
$\tilde U_c=U_x\subseteq X$ is affine.  If $c=y\in Y$, $\tilde U_y$ is acyclic. 

b. Given $x,x'\in \tilde U_y\cap X$, $\tilde U_x\cap \tilde U_{x'}=U_x\cap U_{x'}$ which is quasi-isomorphic to  $U_x\times_{U_y} U_{x'}$,  then $\tilde U_x\cap \tilde U_{x'}$
is acyclic and 
$$\OO_{C(f),xx'}=\OO_{xx'}=\OO_{x}\otimes_{\OO_y}\OO_{x'}=\OO_{C(f),x}\otimes_{\OO_{C(f),y}} \OO_{C(f),x'}$$

c. Given, $y',y''\in \tilde U_{y}\cap Y$, $\tilde U_{y'}\cap \tilde U_{y''}=F^{-1}(U_{y'}\cap U_{y''})$, which is acyclic
because $U_{y'}\cap U_{y''}$ is acyclic, and
$$\OO_{C(f),y'y''}=\OO_{y'y''}=\OO_{y'}\otimes_{\OO_y}\OO_{y''}=\OO_{C(f),y'}\otimes_{\OO_{C(f),y}} \OO_{C(f),y''}$$

d. Given  $x,y'\in \tilde U_y$, where $x\in X$ and $y'\in Y$.  Observe that $\tilde U_{x}\cap \tilde U_{y'}=U_{xy'}$ and $U_{xy'}=f_x^{-1}(U_{f(x)}\cap U_{y'})$, which  is affine since $f_x\colon U_x\to U_{f(x)}$ is affine.  and $U_{f(x)}\cap  U_{y'}\subset U_y$ is affine. 
Finally,
$$\OO_{C(f),xy'}=\OO_{xy'}\overset*=\OO_{x}\otimes_{\OO_y}\OO_{y''}=\OO_{C(f),x}\otimes_{\OO_{C(f),y}} \OO_{C(f),y'}$$
($*$ observe that $U_x\times_{U_y}U_{y'}=U_{xy'})$.

Therefore $C(f)$ is schematic.

If $\M$ is a $\OO_{C(f)}$-quasi-coherent module, $F_*\M$ is a quasi-coherent $\OO_Y$-module since
$$(F_*\M)_y\otimes_{\OO_{Y,y}}\OO_{Y,y'}=
\M_y\otimes_{\OO_{C(f),y}}\OO_{C(f),y'}=\M_{y'}=(F_*\M)_{y'}.$$
By Theorem \ref{Tcohsch}, $F$ is schematic.
$F$ is a quasi-isomorphism since $F_*\OO_{C(f)}=\OO_Y$ and  $F^{-1}(U_y)=\tilde U_y$, for any $y\in Y$.

%3. 
%Let $y\in W\subset C(f)$.
%The morphism $$\OO_{C(f),y}=\OO_{Y,y}= \OO_X(f'^{-1}(U_{y}))\to \prod_{x'  \in f'^{-1}(U_{y})} \OO_{X,x'}=\prod_{x'\in X\cap \tilde U_y} \OO_{C(f),x'}$$
%is faithfully flat, hence $y$ is removable.

$\Leftarrow)$ The morphism $f$ is the composition of the open immersion $X\hookrightarrow C(f)$ and the quasi-isomorphism $F\colon C(f)\to Y$, hence $f$ is a quasi-open immersion.

\edemo

\prop Let $f\colon X\to Y$ be a schematic morphism. Then, $f$ is a quasi-isomorphism iff it is a faithfully flat quasi-open immersion.\eprop

\demo $\Rightarrow)$ It is Remark \ref{O12.3} and Proposition \ref{P3.13}.

$\Leftarrow$) If $y\in Y$, then $y$ is a removable point of $C(f)$, since the morphism
$$\OO_{C(f),y}=\OO_{Y,y}\to \prod_{x\in f^{-1}(U_y)} \OO_{X,x}=\prod_{x>y,x\in X} \OO_{C(f),x}$$
is faithfully flat. The morphism $X\to C(f)$ is a quasi-isomorphism, since $X=C(f)-Y$ and $C(f)-Y$ is quasi-isomorphic to  $C(f)$. Finally, $X$ is quasi-isomorphic to  $Y$, since $C(f)$ is quasi-isomorphic to  $Y$.

\edemo

%\coro Let $X$ be a schematic finite space. Then, $X$ is affine iff
%$\OO(X)\to \proda{x\in X} \OO_x$ is faithfully flat, $U_{xy}$ is acyclic and $\OO_{xy}=\OO_x\otimes_{\OO(X)} \OO_y$, for any $x,y\in X$.\ecoro
%
%\demo 
%$X$ is affine iff the morphism $\pi\colon X\to (*,\OO(X))$ is a quasi-isomorphism, that is, it is faithfully flat and it is a quasi-open immersion, which is equivalent to saying that $\pi$ is faithfully 
%flat, and the diagonal morphism $\delta\colon X \to X\times_{\OO(X)} X$ is a quasi-isomorphism. Observe that $U_{xy}$ is acyclic iff it is affine, since it is an open subset of the affine finite space $U_x$. Then, $\delta$ is a quasi-isomorphism iff $\delta_*\OO_X=\OO_{X\times_{\OO(X)} X}$ and $U_{xy}=\delta^{-1}(U_x\times U_y)$ is affine , for any $x,y\in X$, that is, 
%$\OO_{xy}=\OO_x\otimes_{\OO(X)} \OO_y$ and $U_{xy}$ is acyclic, for any $x,y\in X$.
%\edemo
%

\teor Let $f\colon X\to Y$ be a schematic morphism. Then, $f$ is a quasi-open immersion iff $f$ is flat and the morphism $f^*f_*\mathcal M\to \mathcal M$ is an isomorphism
for any quasi-coherent $\OO_X$-module.\eteor

\demo $\Rightarrow)$ The diagonal morphism $\delta\colon X\to X\times_YX$ is a quasi-isomorphism.
Then, $\delta$ is affine. Consider the projections, $\pi_1,\pi_2\colon X\times_Y X\to X$. 
The morphism $f$ is flat, then $f^*f_*\mathcal M=\pi_{2*}\pi_1^*\mathcal M$, by Theorem \ref{Tcbp}. Observe that
$$\aligned (\pi_{2*} \pi_1^*\mathcal M)_x & =\Gamma(X\times_Y U_x, \pi_1^*\mathcal M)=
\Gamma(X\times_Y U_x, \delta_*\delta^*\pi_1^*\mathcal M)
=\Gamma(X\times_Y U_x, \delta_*\mathcal M)\\ & =\Gamma(U_x, \mathcal M) =\mathcal M_x,\endaligned$$
for any $x\in X$. Therefore, the morphism $f^*f_*\mathcal M\to \mathcal M$ is an isomorphism.

$\Leftarrow)$ 
Let $i\colon U_x\to X$ be the obvious  inclusion and denote $i_*\M_{|U_x}=\M_{U_x}$.
The natural morphism $f^*f_*\M_{U_x}\to \M_{U_x}$ is an isomorphism. Then, 
$$\aligned (\pi_1^*\M)_{(x,x')} & =\M_x\otimes_{\OO_{f(x)}}\OO_{x'}=(f_*\M_{U_x})_{f(x)}\otimes_{\OO_{f(x)}}\OO_{x'}
 =(f^*f_*\M_{U_x})_{x'}=(\M_{U_x})_{x'}\\ &=\M(U_{xx'})=(\delta_*\M)_{(x,x')}.\endaligned$$
Hence, $\delta_*$ is an exact functor, since $\pi_1$ is flat.
By Theorem \ref{Tafex}, $\delta$ is an affine morphism. 
Besides, $\OO_{X\times_Y X}=\pi_1^*\OO_X=\delta_*\OO_X$. Hence, $\delta$ is a quasi-isomorphism and $f$ is a quasi-open immersion.

\edemo

\lema  Let $f\colon X\to Y$ be a schematic morphism and suppose that $X$ is affine and $f_*\OO_X=\OO_Y$. Then, $f$ is a quasi-open immersion.\elema

\demo We have to prove that $C(f)$ is a schematic finite space.
By Proposition \ref{4pg20}, we have to prove that
$\OO_{z_1}\otimes_{\OO_z}\OO_{z_2} \to \prod_{w\in U_{z_1z_2}} \OO_w$ is faithfully flat, for any $z\leq z_1,z_2\in C(f)$.

Suppose that  $z_1,z_2\in X$. The epimorphism
$\OO_{z_1}\otimes_{\OO(X)}\OO_{z_2}\to \OO_{z_1}\otimes_{\OO_z}\OO_{z_2}$ is an isomorphism, since
the composite morphism $\OO_{z_1}\otimes_{\OO(X)}\OO_{z_2}\to \OO_{z_1}\otimes_{\OO_z}\OO_{z_2}\to \OO_{z_1z_2}$ is an isomorphism. Besides, the morphism
$\OO_{z_1z_2}\to \prod_{w\in U_{z_1z_2}} \OO_w$ is faithfully flat, since $X$ is affine.

Suppose that $z_1\in X$ and that $z_2 \in Y$ (then $z\in Y$). Observe that
$U_{z_1}\cap f^{-1}({U_{z_2}})=U_{z_1z_2}=\{c\in C(f)\colon c\geq z_1,z_2\}$. Then,
$$\OO_{z_1}\otimes_{\OO_z}\OO_{z_2}=\OO_{z_1}\otimes_{\OO_z}\OO_z\otimes_{\OO_z}\OO_{z_2}
=\OO_{z_1}\otimes_{\OO_z}(\OO_z\otimes_{\OO_{f(z)} }\OO_z)\otimes_{\OO_z}\OO_{z_2}=\OO_{z_1}\otimes_{\OO_{f(z)}}\OO_{z_2}=\OO_{z_1z_2}\to \prod_{w\in U_{z_1z_2}} \OO_w$$
is faithfully flat, since $U_{z_1z_2}\subset X$ is affine, by Proposition \ref{K6}.

Suppose that $z_1,z_2\in Y$. Observe that $U_{C(f),z_1z_2}=
U_{Y,z_1z_2}\coprod f^{-1}(U_{Y,z_1z_2})$ and
$\OO_Y(U_{Y,z_1z_2})=\OO_X(f^{-1}(U_{Y,z_1z_2}))$. 
The morphism
$\OO_{z_1}\otimes_{\OO_z}\OO_{z_2}=\OO_{Y,z_1z_2}\to 
\prod_{w\in U_{Y,z_1z_2}} \OO_w$ is faithfully flat, since
$U_{Y,z_1z_2}$ is affine. 
For any open subset $V\subset X$ and any $x\in V$, the morphism
$\OO_X(V)=\OO(X)\otimes_{\OO(X)} \OO(V)\to \OO_x\otimes_{\OO(X)}\OO(V)\overset{\text{\ref{11}}}=\OO_x$ is flat.  Hence, the morphism 
$\OO_{z_1}\otimes_{\OO_z}\OO_{z_2}=\OO_{Y,z_1z_2}=
\OO_X(f^{-1}(U_{Y,z_1z_2}))\to 
\prod_{x\in f^{-1}(U_{Y,z_1z_2})} \OO_x$ is flat. Therefore,
the morphism
$$\OO_{z_1}\otimes_{\OO_z}\OO_{z_2}\to 
\prod_{w\in U_{C(f),z_1z_2}} \OO_w=\prod_{w\in U_{Y,z_1z_2}} \OO_w\times \prod_{x\in f^{-1}(U_{Y,z_1z_2})} \OO_x$$
is faithfully flat.

\edemo

\prop Let $f\colon X\to Y$ be a schematic morphism and suppose that $X$ is affine. Then, there exist an open inmersion $i\colon X\to Z$ such that $i_*X=\OO_Z$ and an affine morphism $g\colon
Z\to Y$ such that $f=g\circ i$. \eprop

\demo The obvious morphism $f'\colon X\to (Y,f_*\OO_X)$, $f'(x)=f(x)$ is a quasi-open inmersion by the lemma above and Example \ref{E8.2}. Let $i\colon X\to C(f')$, $\pi\colon 
C(f')\to (Y,f_*\OO_X)$ and let $\Id\colon (Y,f_*\OO_X)\to Y$
be the obvious morphism. Observe that  $i$ is an open immersion, $g:=\Id\circ \pi$ is affine, since $\pi$ is a quasi-isomorphism and $\Id$ is affine, and 
$$f=\Id \circ f'=\Id\circ \pi\circ i=g\circ i.$$

\edemo 

\section{Quasi-closed immersions}

%\prop Let $(X,\OO)$ be a schematic finite space, $\mathcal O\to \mathcal B$ a morphism of sheaves of rings such that $\mathcal B$ is a quasi-coherent module. Then, $(X,\mathcal B)$ is schematic finite space and $\Id:(X,\mathcal B)\to (X,\OO)$ is an affine morphism.\eprop
%
%\demo The morphism $\OO_y\otimes_{\OO_x}\OO_{y'} \to \prod_{z\in U_{yy'}} \OO_z$ is faithfully flat, for any $y,y'\geq x$, by Proposition \ref{4pg20}. Tensoring by $\mathcal B_x\otimes_{\OO_x}$,
%the morphism
%$$\B_y\otimes_{\B_x}\B_{y'} \to \prod_{z\in U_{yy'}} \B_z$$
%is faitfully flat. By Proposition \ref{4pg20}, $(X,\B)$ is a schematic finite space. If $\M$ is a quasi-coherent $\B$-module,
%then $\M_x\otimes_{\OO_x}\OO_{x'}= \M_x\otimes_{\B_x}\B_x\otimes_{\OO_x}\OO_{x'}=\M_x\otimes_{\B_x}\B_{x'}=\M_{x'}$.
%Hence, $\Id_*\M$ is quasi-coherent and $\Id$ is schematic by Theorem \ref{Tcohsch}. $\Id$ is affine since $\Id^{-1}(U_x)=U_x$, for any $x\in X$.
%
%\edemo 
%
%
%\prop Let $f\colon X\to Y$ be a schematic morphism. The obvious morphism
%$$f'\colon (X,\OO_X)\to (Y,f_*\OO_Y),\quad f'(x)=f(x)$$ is a 
%schematic morphism, and the diagram
%$$\xymatrix{(X,\OO_X)  \ar[rr]^-f \ar[rd]_-{f'} & & (Y,\OO_Y)\\ & (Y,f_*\OO_X) \ar[ru]_-{\Id} &}$$
%is commutative. In adittion, if $f$ is affine $f'$ is a quasi-isomorphism.\eprop
%

Let $\mathcal I\subset \OO_X$ be a quasi-coherent ideal and 
$(\mathcal I)_0:=\{x\in X\colon (\OO_X/\mathcal I)_x\neq 0\}$, which is a closed subspace of 
$X$. Consider the schematic space $(X,\OO_X/\mathcal I)$ and observe that $x\in X\backslash Y$ iff $(\OO_X/\mathcal I)_x=0$.
Hence, $X\backslash Y$ is a set of removable points of 
$(X,\OO_X/\mathcal I)$ and the obvious morphism 
$((\mathcal I)_0,\OO_X/\mathcal I)\to (X,\OO_X/\mathcal I)$ is a quasi-isomorphism.
We shall say that the composition of the affine morphisms 
$$((\mathcal I)_0,\OO_X/\mathcal I)\to (X,\OO_X/\mathcal I)\to (X,\OO_X)$$
is a closed immersion.

Let $f\colon X'\to X$ be a schematic morphism. Let $\mathcal I=\Ker[\OO_{X}\to f_*\OO_{X'}]$.
The obvious morphism $f'\colon X'\to (X,\OO_X/I)$, $f'(x)=f(x)$ is schematic since
$f'_*\M=f_*\M$ is a quasi-coherent $\OO_X/I$-module, for any quasi-coherent $\OO_{X'}$-module $\M$, because  it is a quasi-coherent $\OO_X$-module.  Obviously, $f$ is the composition of the morphisms $X'\to (X,\OO_X/I) \to (X,\OO_X)$. 
Assume that $\OO_{X',x'}\neq 0$ for any $x'\in X'$ (recall that if $\OO_{X'.x'}=0$ then $x'$ is a removable point).
 The closure of $f(X')$ in $X$ is $ (\mathcal I)_0$:  $x\in X\backslash  (\mathcal I)_0$
iff $(\OO_X/\mathcal I)_x= 0$, which is equivalent to saying that
$(f_*\OO_{X'})_x= 0$ since the morphism of sheaves (of rings) $\OO_{X}/\mathcal I\to f_*\OO_{X'}$ is inyective. 
$(f_*\OO_{X'})_x=\OO_{X'}(f^{-1}(U_x))= 0$ iff  $f^{-1}(U_x)=\emptyset$. 
The morphism  $X'\to ((\mathcal I)_0,\OO_X/\mathcal I)$, $x\mapsto f(x)$  is schematic and
$f'$ is the composition of the morphisms $X'\to ((I)_0,\OO_X/I) \to (X,\OO_X/I)$
(see Proposition \ref{P6.11}).

\defi Let $f\colon X'\to X$ be a schematic morphism. We shall say that $f$ is a quasi-closed immersion if it is affine and the morphism $\OO_X\to f_*\OO_{X'}$ is an epimorphism.\edefi

Suppose  that $\OO_{X',x'}\neq 0$, for any $x'\in X'$, let
$f\colon X'\to X$ be a quasi-closed immersión and $\mathcal I=\Ker[\OO_X\to f_*\OO_{X'}]$.  Then $
f$ is the composition of a quasi-isomorphism $X'\to ((I)_0,\OO_X/I)$ and a closed immersion
$ ((I)_0,\OO_X/I)\to (X,\OO_X)$.

\section{{\rm Spec}\hspace{0.1cm}{\sl O}$_X$}

Let $\{X_i, f_{ij}\}_{i,j\in I}$ (where $\#I<\infty$)
be a direct system of morphisms of ringed spaces.
Let $\ilim{i} X_i$ be the direct limit of the topological spaces $X_i$: $\ilim{i} X_i=\coproda{i} X_i/\!\sim$, where $\sim$ is the equivalence relation generated by the relation $x_i\sim f_{ij}(x_i)$, and
 $U\subseteq \ilim{i} X_i$ is an open subset iff $f_j^{-1}(U)$ is an open subset for any $j\in I$, where $f_j\colon X_j\to \ilim{i} X_i$ is the natural map. 
Define $\OO_{\ilim{i} X_i}(V):=\plim{i\in I} \OO_{X_i}(f_i^{-1}(V))$, for any open set $V\subseteq \ilim{i} X_i$.
It is well known that $(\ilim{i} X_i, \OO_{\ilim{i} X_i})$ is the direct limit of the direct system of morphisms $\{X_i,f_{ij}\}$ in the category of ringed spaces.

\defi Given a schematic finite space $X$  we shall denote
 $$\Spec \OO_X:=\ilim{x\in X} \Spec \OO_x.$$
(the sheaf of  rings considered on $\Spec O_x$ is the sheaf of localizations of $\OO_x$, $\tilde \OO_x$.)
\edefi

Observe that $\OO_{\Spec \OO_X}(\Spec \OO_X)=\plim{x\in X} \OO_x=\OO_X(X)$.

\ejem Obviously, if  $X=U_x$, then $\Spec \OO_X=\Spec \OO_x$ and $\OO_{\Spec\OO_X}=\tilde \OO_x$.
\eejem

Consider the following relation on $\coprod \Spec \OO_{x_i}$: Given $\pp\in \Spec\OO_x$ and  $\qq\in\Spec \OO_y$
we shall say that $\pp\equiv \qq$ if there exist $u\geq x,y$ and  $\rr \in \Spec \OO_u$ such that  
the given morphisms $\Spec \OO_u\iny \Spec \OO_x$ and $\Spec \OO_u\iny \Spec \OO_y$ map $\rr$ to $\pp$ and $\rr $ to $\qq$, respectively (recall Proposition \ref{last}). Let us prove that $\equiv$ is an equivalence relation: Let  $\pp\equiv \qq$,  ($\rr\mapsto \pp,\qq$) and   $\qq\equiv \qq'$ ($\qq'\in \Spec \OO_z$, there exist $u'\geq y,z$ and $\rr'\in\Spec\OO_{u'}$ such that $\rr'\mapsto \qq,\qq'$). Recall that 
$\OO_{uu'}=\OO_{u}\otimes_{\OO_y}\OO_{u'}$ and that $\OO_{uu'}\to \prod_{w\in U_{uu'}} \OO_{w}$ is faithfully flat. Then, $$\Spec \OO_u\cap \Spec \OO_{u'}=\Spec \OO_{uu'}=\cup _{w\in U_{uu'}} \Spec \OO_{w}$$
Since, $\qq=\rr=\rr' \in \Spec \OO_u\cap \Spec \OO_{u'}$ there exists  $v\in U_{uu'}$ and  $\rr''\in 
\Spec \OO_{v}$ such that $\rr''\mapsto \rr,\rr'$. Then, $\rr''\mapsto \pp,\qq'$ and $\pp\equiv \qq'$.

Observe that $\Spec \OO_X=\coprod_{x\in X}  \Spec \OO_{x}/\equiv$ as topological spaces. Besides, the morphisms
$\Spec \OO_{x}\to \Spec \OO_X$ are injective, $\Spec \OO_X=\cup_{x\in X} \Spec\OO_x$ as topological spaces ($U\subseteq \Spec \OO_X$ is an open set iff $U\cap \Spec\OO_x$ is an open set, for any $x\in X$).

\lema \label{radical} Let $A\to B$ be a flat morphism and assume $B\otimes_A B=B$. If $I\subseteq A$ is a radical ideal, then $I\cdot B$ is a radical ideal of $B$. \elema

\demo Let $\pp\in\Spec B\subset\Spec A$ and recall Notation \ref{Notation}. Then,
$$(\rad (I\cdot B))_\pp=\rad (I\cdot B_\pp)\overset{\text{\ref{last}}}=\rad (I\cdot A_\pp)=\rad(I)\cdot A_\pp=I\cdot A_\pp\overset{\text{\ref{last}}}=I\cdot B_\pp=(I\cdot B)_\pp$$
Therefore, $\rad (I\cdot B)=I\cdot B$.
\edemo

%\prop Sea $X$ un espacio finito esquem\'{a}tico y $x\leq y,y'\in X$. En $\Spec \OO_x$ se cumple que
%$$\Spec \OO_{yy'}=(\Spec \OO_y)\cap (\Spec \OO_{y'})=\cup_{z\in U_{yy'}} \Spec \OO_z$$\eprop
%
%\demo El morfismo $\OO_{U_{yy'}}\to \prod_{z\in {U_{yy'}}} \OO_z$ es fielmente plano, luego $\Spec \prod_{z\in {U_{yy'}}} \OO_z\to \Spec \OO_{U_{yy'}}$ es epiyectivo. Adem\'{a}s, $\OO_{U_{yy'}}=\OO_y\otimes_{\OO_x}\OO_{y'}$ y es localmente una localizaci\'{o}n de $\OO_x$.
%Del diagrama
%$$\xymatrix{  & & \Spec\OO_y \ar@{^{(}->}[rd] & \\ \Spec \prod_{z\in {U_{yy'}}} \OO_z \ar[r]^-{epi} 
%& \Spec \OO_{U_{yy'}}=\Spec \OO_y\times_{\Spec \OO_x} \Spec\OO_{y'} \ar@{^{(}->}[ru] \ar@{^{(}->}[rd] & & \Spec\OO_x\\  & & \Spec\OO_{y'} \ar@{^{(}->}[ru]& } $$
%se concluye.
%
%\edemo

\prop Let $X$ be a schematic finite space. Let $\mathcal I\subseteq \OO_X$ be  a quasi-coherent ideal.  The ideal $\rad\mathcal I\subset \OO_X$, defined by $(\rad\mathcal I)_x:=\rad\mathcal I_{x}$, for any $x\in X$, is a quasi-coherent ideal of $\OO_X$.\eprop

\demo We only have to prove that given a flat morphism $A\to B$ such that $B\otimes_AB=B$ and an ideal $I\subseteq A$,  then $(\rad I)\cdot B=\rad (I\cdot B)$. This is a consequence of  Lemma \ref{radical}.
\edemo

\obse $(\rad\mathcal I)(U)=\rad (\mathcal I(U))$, for any open subset 
$U\subset X$: $$(\rad\mathcal I)(U)=\plim{x\in U} (\rad\mathcal I)_x=
\plim{x\in U} (\rad\mathcal I_x)=\rad\plim{x\in U} \mathcal I_x=\rad\mathcal I(U).$$

\eobse

\defi Let $X$ be a schematic finite space. We shall say that a  quasi-coherent ideal  $\mathcal I\subseteq \OO_X$ is radical if $\mathcal I=\rad(\mathcal I)$. 
\edefi

\nota Let $X$ be a schematic finite space. Given a quasi-coherent ideal $\mathcal I\subset \OO_X$,  we shall denote
$$(\mathcal I)_0:=\cupa{x\in X} \{\pp\in \Spec \OO_x\colon \mathcal I_x\subseteq \pp\}\subseteq\Spec\OO_X$$
Given a closet subset $C\subset \Spec \OO_X$,  let $\mathcal I_C\subset \OO_X$ be the radical quasi-coherent ideal  defined
by $\mathcal I_{C,x}:=\capa{\pp'\in C\cap \Spec \OO_x} \pp'\subset \OO_x$,\footnote{If  $C\cap \Spec \OO_x=\emptyset$, then $\capa{\pp'\in C\cap \Spec \OO_x} \pp':=\OO_X$.}
for any $x\in X$.
\enota

\prop The maps
$$\{\text{Closed subspaces of $\Spec \OO_X$}\} \longleftrightarrow
\{\text{Radical quasi-coherent ideals of $\OO_X$}\}$$
$$\xymatrix @R8pt{ C \ar@{|->}[r] \quad & \quad \qquad \mathcal I_C \\ (\mathcal I)_0 \quad  & \quad \quad   \qquad   \ar@{|->}[l] \mathcal I \quad }$$
are mutually inverse.
\eprop

\nota Given a ring $B$ and $b\in B$
we denote $B_b=\{1,b^{-1},b^{-2},\cdots\}\cdot B$.\enota

\prop \label{P16.2} If $X$ is an affine finite space, then $\Spec \OO_X=\Spec \OO(X)$. \eprop

\demo The morphism $\OO(X)\to \OO_x$ is flat 
and $\OO_x\otimes_{\OO(X)}\OO_x\overset{\text{\ref{11}}}=\OO_x$, then $\Spec \OO_x\iny \Spec \OO(X)$ is a subspace. The morphism $\OO(X)\to \prod_{x\in X} \OO_x$
is faithfully flat, then the induced morphism $\coprod_{x\in X}  \Spec \OO_x\to \Spec \OO(X)$ is surjective.
The sequence of morphisms
$$\OO(X)\to \prod_{x\in X} \OO_x\dosflechas \prod_ {x\leq x'\in X}  \OO_{x'}$$ 
is exact.  Then, the natural morphism $f\colon\Spec \OO_X\to \Spec \OO(X)$ is continuous and bijective. Given a closed set $C\subset\Spec \OO_X$, let $\mathcal I_C$ be the  radical quasi-coherent ideal of $\OO_X$ associated with $C$. $\mathcal I_C=\tilde{\mathcal I_C(X)}$ since $X$ is affine. Recall Notation \ref{N3.9}. Then,
$C\cap \Spec \OO_x=(\mathcal I_{C,x})_0=(\mathcal I_C(X)\cdot \OO_x)_0$. Hence
$f(C)=(\mathcal I_C(X))_0$ and $f$ is a homeomorphism.

Also observe that $(\plim{x} \OO_x)_{(a_x)}=
\plim{x\in X} \OO_{x,a_x}$, for any $(a_x)\in \plim{x\in X} \OO_x\subseteq \prod_{x\in X} \OO_x$, hence $\OO_{\ilim{x\in X} \Spec \OO_x}=\widetilde{\OO(X)}$.
\edemo

\defi Let  $f\colon X\to Y$ be a schematic morphism. Consider the morphisms $\OO_{f(x)}\to \OO_x$, which induce the scheme morphisms $\Spec \OO_x\to \Spec \OO_{f(x)}$, which induce a morphism of ringed spaces $$\tilde f\colon \Spec \OO_X\to \Spec \OO_Y.$$
We shall say that $\tilde f$ is the morphism induced by $f$.
\edefi

\prop \label{P16.4} Let $f\colon X\to Y$ be a quasi-isomorphism. Then, the morphism  induced by $f$, $\tilde f\colon \Spec \OO_X\to \Spec \OO_Y$,  is an isomorphism.
\eprop

\demo Observe that
$$\aligned \Spec \OO_X & =\ilim{x\in X} \Spec \OO_x=
\ilim{y\in Y}\ilim{x\in f^{-1}(U_y)} \Spec \OO_x\overset{\text{\ref{P16.2}}}
=\ilim{y\in Y}
\Spec \OO_X(f^{-1}(U_y))\\ &= 
\ilim{y\in Y} \Spec \OO_Y(U_y) = \Spec \OO_Y.\endaligned$$

\edemo

\prop Let $X$ be a schematic finite space, $U\overset i\subset X$ an open subset   and $\mathcal I\subseteq \OO_U$  a quasi-coherent ideal. Then, there exists a quasi-coherent ideal $\mathcal J\subseteq \OO_X$ such that $\mathcal J_{|U}=\mathcal I$.
\eprop

\demo $\mathcal J:=\Ker[\OO_X\to i_*(\OO_U/\mathcal I)]$
holds $\mathcal J_{|U}=\mathcal I$. \edemo

\nota Given a schematic finite space $X$ we shall denote
$\tilde X=\Spec\OO_X$.

\enota

\prop \label{P16.12} Let $X$ be a schematic finite space and $U\subset X$ an open subset. Then,  \enumera

\item $\tilde U$ is a topological subspace of $\tilde X.$

\item $\tilde U=\capa{\tilde U\subseteq \text{ open subset  }\bar V\subseteq \tilde X } \bar V.$

\eenumera
 \eprop

\demo 1. Given a closed set $C\subset \tilde U$, let $\mathcal I_C\subseteq \OO_U$ be the radical quasi-coherent ideal associated. Let $\mathcal J\subseteq \OO_X$ be the quasi-coherent ideal such that $\mathcal J_{|U}=\mathcal I$.  Then, the closed subset  $D=(\mathcal J)_0=(\rad\mathcal J)_0$ of $\tilde X$ holds that $D\cap \tilde U=C$.

2. Let   $\pp\in \tilde X-\tilde U$. Let $\mathcal P\subset \OO_X$ be the sheaf of ideals defined by $\mathcal P_x=\pp\subset \OO_x$ if $\pp\in \Spec \OO_x$ and $\mathcal P_x=\OO_x$ if $\pp\notin \Spec \OO_x$. By Proposition \ref{sudor}, $\mathcal P$ is quasi-coherent.
$(\mathcal P)_0\subset \Spec \OO_X$ is the closure of $\pp$ and 
$(\mathcal P)_0\cap \tilde U_x=\emptyset$, for any 
$\tilde U_x\subset \tilde U$, hence  $(\mathcal P)_0\cap \tilde U=\emptyset$.
Then, $\tilde U$ is equal to the intersection of the open subsets $\bar V\subseteq \tilde X$,   such that $\tilde U\subseteq \bar V$.

\edemo

\defi Let $X$ be a schematic finite space. We shall say that a quasi-coherent $\OO_X$-module  $\mathcal M$ is finitely generated if $\M_x$ is a finitely generated $\OO_{x}$-module, for any $x\in X$.\edefi

\prop Let $X$ be an affine finite space and $\M$ a a quasi-coherent $\OO_X$-module. Then, $\M$ is finitely generated iff $\M(X)$ is a finitely generated $\OO(X)$-module.\eprop

\demo $\Rightarrow)$ Given $x\in X$, $\M_x=\M(X)\otimes_{\OO(X)}\OO_x$.
Let $N^x\subset \M(X)$ be a finitely generated $\OO(X)$-submodule such that $N^x\otimes_{\OO(X)}\OO_x=\M_x$ and  $N:=\sum_{x\in X} N^x$. Then $N=\M(X)$, since $N\otimes_{\OO(X)}\OO_x=M_x$ for any $x\in X$.

$\Leftarrow)$ $\M_x=\M(X)\otimes_{\OO(X)}\OO_x$ is a finitely generated $\OO_x$-module, for any $x\in X$.
\edemo

\prop Let $X$ be a schematic finite space. Any quasi-coherent $\OO_X$-module is the direct limit of its finitely generated submodules.\eprop

\demo Let $\M$ be a quasi-coherent $\OO_X$-module.
Let us fix $x_1\in X$ and  a finitely generated submodule  $N_1\subset \M_{x_1}$.
Consider the inclusion morphism $i_1\colon U_{x_1}\iny X$ and let
$\M_1:=\Ker[\M\to i_{1*} (\M_{|U_{x_1}}/\tilde N_1)]$. Observe that 
$\M_1\subset\M$ and $\M_{1|U_1}=\tilde N_1$. Given $x_2\in X$, let $N_{2}\subset \M_{1,x_2}$ be a finitely generated submodule such that $N_{2}\otimes_{\OO_{x_2}}\OO_y =\M_{1,y}$, for any
$y\in U_ {x_1}\cap U_ {x_2}$. Let $U_2=U_{x_1}\cup U_{x_2}$ and let $\Nc_2\subset \M_{1|U_2}$ be  
the finitely generated $\OO_{U_2}$-module such that $\Nc_{2|U_{x_1}}=\tilde N_1$ and $\Nc_{2|U_{x_2}}=\tilde N_2$.
Consider the inclusion morphism  $i_2\colon U_2\iny X$
and let $\M_2:=\Ker[\M_1\to i_{2*} (\M_{1|U_2}/\Nc_{2})]$.
Observe that $\M_{2|U_2}=\Nc_2$.
Given $x_3\in X$, let $N_{3}\subset \M_{2,x_3}$ be a finitely generated submodule such that $N_{3}\otimes_{\OO_{x_3}}\OO_y=\M_{2,y}$ for any $y\in U_{x_3}\cap U_2$. Let $U_3:=U_{2}\cup U_{x_3}$ and let $\Nc_3\subset \M_{2|U_3}$ be  
the finitely generated $\OO_{U_3}$-module such that $\Nc_{3|U_{2}}=\Nc_2$ and $\Nc_{3|U_{x_3}}=\tilde N_3$. Consider the inclusion morphism  $i_{3}\colon U_{3}\iny X$
and let $\M_3:=\Ker[\M_2\to i_{3*} (\M_{2|U_{3}}/\Nc_{3})]$.
Observe that $\M_{3|U_3}=\Nc_3$. So on we shall get a finitely quasi-coherent $\OO_X$-submodule $\M_n\subset\M$ such that $\M_{n,x_1}=N_1$. Now it is easy to prove this proposition.
\edemo

\coro \label{C15.23} Let $X$ be a schematic finite space. Any quasi-coherent ideal  $\mathcal I \subset \OO_X$ is the direct limit of its finitely generated ideals $\mathcal I_i\subset \mathcal I$.
\ecoro

\lema \label{L15.20} Let $X$ be a schematic finite space, $\bar U\subset \tilde X$ an open subset and $C=\tilde X-\bar U $. Then, $\bar U$ is quasi-compact
iff there exists a finitely generated ideal $\mathcal I\subset \OO_X$ such that $({\mathcal I})_0=C$. 
\elema

\demo $\Rightarrow)$ Consider the quasi-coherent ideal $\mathcal I_C\subset \OO_X$. Let $J=\{\mathcal I_j\}_{j\in J}$ the set of finitely generated ideals of $\OO_X$ contained in $\mathcal I_C$.
By Corollary \ref{C15.23}, $\mathcal I_C=\ilim{j\in J} \mathcal I_j$. Then, $C=(\mathcal I_C)_0=(\ilim{j\in J} \mathcal I_j)_0=\cap_{j\in J} (\mathcal I_j)_0$ and $\bar U=\cup_{j\in J} (\tilde X-(\mathcal I_j)_0)$. There exists $j\in J$ such that $\bar U=\tilde X-(\mathcal I_j)_0$, since $\bar U$ is quasi-compact. Hence, $C=(\mathcal I_j)_0$.

$\Leftarrow$) Let $x\in X$, then $\mathcal I_x=(a_1,\ldots,a_n)\subset \OO_{x}$ is finitely generated. $C\cap \tilde U_{x}=(\mathcal I_x)_0=
\cap_{i} (a_i)_0$, then $\bar U\cap \tilde U_{x}=\cup_i \Spec \OO_{x,a_i}$ is quasi-compact. Therefore, $\bar U=\cup_{x} (\bar U\cap \tilde U_{x})$ is quasi-compact.

\edemo

\prop \label{P16.16} Let $X$ be a schematic finite space. Then, 
\enumera 
\item The intersection of two quasi-compact open subsets of $\tilde X$ is quasi-compact.

\item The family  of quasi-compact open subsets of $\tilde X$ is a basis for the topology of  $\tilde X$.

\item If $\bar V\subseteq \tilde X$ is a quasi-compact open subset then $\bar V\cap \tilde U$ is quasi-compact, for any open subset $U\subset X$.

\eenumera
\eprop

\demo  1. Let $\bar U_1,\bar U_2\subset \tilde X$ be two quasi-compact open subsets, $C_1:=\tilde X-\bar U_1$, $C_2:=\tilde X-\bar U_2$,  and $\mathcal I_1,\mathcal I_2\subset \OO_X$ two finitely generated ideals such that $C_1=(\mathcal I_1)_0$ and
$C_2=(\mathcal I_1)_0$.  Then, $C_1\cup C_2=(\mathcal I_1)_0\cup 
(\mathcal I_2)_0=(\mathcal I_1\cdot\mathcal I_2)_0$ and $\bar U_1\cap\bar U_2=\tilde X-(\mathcal I_1\cdot \mathcal I_2)_0$. By Lemma \ref{L15.20}, $\bar U_1\cap\bar U_2$ is quasi-compact.

2. Let $\bar U\subset \tilde X$ be an open subset and 
$C=\tilde X-\bar U$. $\mathcal I_C=\ilim{j\in J} \mathcal I_j$, where
$\{\mathcal I_j\}_{j\in J}$ is the set of the finitely generated of $\OO_X$ contained in $\mathcal I_C$. Then, $C=(\mathcal I_C)_0=
(\ilim{j\in J} \mathcal I_j)_0=\cap_{j\in J} ( \mathcal I_j)_0$
and $\bar U=\cup_{j\in J} (\tilde X-(\mathcal I_j)_0)$, where the open subsets $\tilde X-(\mathcal I_j)_0$ are quasi-compact by Lemma 
\ref{L15.20}.

3. Let $C=\tilde X-\bar V$· and let $\mathcal I\subset \OO_X$ be a finitely generated ideal such that
$C=(\mathcal I)_0$. Then, $C\cap \tilde U=(\mathcal I_{|U})_0$
and $\bar V\cap\tilde U=\tilde U-(\mathcal I_{|U})_0$. By Lemma 
\ref{L15.20}, $\bar V\cap\tilde U$ is quasi-compact.

\edemo

\coro Let $X$ be a schematic finite space, $U\subseteq X$ an open subset and $\bar U\subset \tilde U$  a quasi-compact open subset. Then,

\enumera \item There exists a quasi-compact open subset $\bar W\subset \tilde X$, such that, $\bar W\cap \tilde U=\bar U$.

\item $\bar U$ is equal to the intersection of the quasi-compact open subsets of $\tilde X$ which contain it.\eenumera

\ecoro

\demo 1. By Proposition \ref{P16.12}, there exists an open subset $\bar W'\subseteq \tilde X$ such that $\bar W'\cap \tilde U=\bar U$.
Given $\pp\in \bar U$, there exists a quasi-compact open subset $\bar W_{\pp}\subset \bar W'$ such that $\pp\in \bar W_{\pp}$. There exist $\pp_1,\ldots,\pp_n\in \bar U$ such that
$\bar U\subset \cup_{i=1}^n \bar W_{\pp_i}\subset \bar W'$. Hence, $\bar W:= \cup_{i=1}^n \bar W_{\pp_i}$ holds $\bar W\cap \tilde U=\bar U$.

2. Given an open subset $\bar V\subset \tilde X$ such that $\bar U\subset \bar V$, there exists a quasi-compact open subset $\bar V'\subset \tilde X$ such that $\bar U\subset \bar V'\subset \bar V$. By Proposition \ref{P16.12}, we are done.
%2. Let   $\pp\in \tilde X \backslash \tilde U$ Let $U_x\subset X$ be such that  $\pp\in \tilde U_x$. Let $C$ be the closure of $\pp$ en $\tilde U_x$. Observe  that  $C\cap  \tilde U=\emptyset$, by \ref{sudor} 3.).  $\bar U:=(\tilde U_x\backslash C)\cup \tilde U$ is an open subset of ${\tilde U_x}\cup \tilde U=\tilde U'$ ($U'=U_x\cup U$) and $\pp\notin \bar U$. Let $\bar V\subset \tilde X$ be an open subset such that $\bar V\cap \tilde U'=\bar U$. Observe that $\pp\notin \bar V$ and $\tilde U\subset \bar V$.  As $\tilde U$ is quasi-compact, there exists an quasi-compact open subset $\tilde W$ of $\tilde X$ such that $\pp\notin \bar W$ and $\tilde U\subset \bar W$.Therefore, $\tilde U$ is equal to the intersection of the quasi-compact open subsets of $\tilde X$ which contain it.

\edemo

\lema \label{L16.18} Let $X$ be a schematic finite space, $U_1,U_2\subset X$ open subsets, $\bar V_1\subset \tilde U_1$ and $\bar V_2\subset \tilde U_2$ quasi-compact open subsets
and $\bar W\subset \tilde X$ an open subset such that $\bar V_1\cap \bar V_2\subset \bar W$. Then, there exist open subsets $\bar W_1,\bar W_2\subset \tilde X$ such that $\bar V_1\subset \bar W_1$, $\bar V_2\subset \bar W_2$ and $\bar W_1\cap \bar W_2\subset \bar W$.
\elema

\demo By the quasi-compactness of $\bar V_1$ and $\bar V_2$, to prove this theorem we can easily reduce ourselves to the case in which $\bar V_1=\Spec \OO_{x_1,a_1}\subset \tilde U_{x_1}$ ($a_1\in \OO_{x_1}$) and $\bar V_2=\Spec \OO_{x_2,a_2}\subset \tilde U_{x_2}$ ($a_2\in \OO_{x_2}$).

1. Suppose that  $\bar V_1\cap \bar V_2=\emptyset$.
Let $\OO_{U_{x_1,a_1}}$ be the quasi-coherent $\OO_{U_{x_1}}$-module defined by $\OO_{U_{x_1,a_1}}(U_z)=\OO_{z,a_1}$, for any $z\in U_{x_1}$. Let $i_{x_1}\colon U_{x_1}\subset X$ be the inclusion morphism.
Let $\mathcal I_1$ be the kernel of the natural morphism
$\OO_X\to i_{x_1*}\OO_{U_{x_1,a_1}}$. 
Likewise, define $\OO_{U_{x_2,a_2}}$, $i_{x_2}$ and $\mathcal I_2$.

Observe that $i_{x_1*}\OO_{U_{x_1,a_1}}(U_z)=\OO_{x_1z, a_1}$ for any 
$z\in X$, and the natural morphism 
$$i_{x_1*}\OO_{U_{x_1,a_1}}(U_z)\to
\prod_{y\in U_{x_1z}} \OO_{y,a_1}$$ is injective.
 Then, the sequence of morphisms
$$(**)\qquad 0\to \mathcal I_{1,z}\to \OO_z\to \prod_{y\in U_{x_1z}} \OO_{y,a_1}$$
is exact. Observe that $\bar V_1\cap \tilde U_z=\cup_{y\in U_{x_1z}} \Spec\OO_{y,a_1}$. 
Then, $(\mathcal I_1)_0$ is equal to the closure $Cl(\bar V_1)$ of
$\bar V_1$ in $\tilde X$.

Let $z=x_2$, tensoring $(**)$ by $\otimes_{\OO_{x_2}}  \OO_{x_2,a_2}$,  we obtain the exact sequence
$$0\to \mathcal I_{1,x_2}\otimes_{\OO_{x_2}} \OO_{x_2,a_2}\to \OO_{x_2,a_2}\to \prod_{y\in U_{x_1x_2}} \OO_{y,a_1}\otimes_{\OO_{x_2}}  \OO_{x_2,a_2}.$$
 $\OO_{y,a_1}\otimes_{\OO_{x_2}}  \OO_{x_2,a_2}=0$ since
$\Spec \OO_{y,a_1}\cap \Spec \OO_{x_2,a_2}\subset \bar V_1\cap\bar V_2=\emptyset$. Then, $\mathcal I_{1,x_2}\otimes_{\OO_{x_2}}  \OO_{x_2,a_2}= \OO_{x_2,a_2}$, hence $\mathcal I_{1,x_2}\cdot \OO_{x_2,a_2}= \OO_{x_2,a_2}$. Therefore, $(\mathcal I_1)_0\cap \bar V_2=\emptyset$, that is, 
$Cl(\bar V_1)\cap \bar V_2\overset{*}=\emptyset$.
Let $\mathcal J_1\subset\mathcal I_1$ be a finitely generated ideal such that $\mathcal J_{1,x_2} \cdot \OO_{x_2,a_2}= \OO_{x_2,a_2}$. Again, 
$(\mathcal J_1)_0\cap \bar V_2=\emptyset$ and $\bar V_1\subset (\mathcal J_1)_0$. Likewise, define a finitely generated ideal $\mathcal J_2$ such that $(\mathcal J_2)_0\cap \bar V_1=\emptyset$ and $\bar V_2\subset (\mathcal J_2)_0$. 

Given a subset $Y\subset \tilde X$  denote
$Y^c:=\tilde X-Y$.
Let  $\bar W_1:=(Cl((\mathcal J_1\cdot \mathcal J_2)_0^c)\cup (\mathcal J_2)_0)^c$ 
and $\bar W_2:=(Cl((\mathcal J_1\cdot \mathcal J_2)_0^c)\cup (\mathcal J_1)_0)^c$. Obviously, $\bar W_1\subset ((\mathcal J_1\cdot \mathcal J_2)_0^c\cup (\mathcal J_2)_0)^c=(\mathcal J_1)_0-(\mathcal J_2)_0$ 
and $\bar W_2\subset (\mathcal J_2)_0-(\mathcal J_1)_0$.
Then, $\bar W_1\cap \bar W_2=\emptyset$. We only have to prove
that $\bar V_1\subset \bar W_1$. We know that $\bar V_1\cap 
(\mathcal J_2)_0=\emptyset$, it remains to prove that $\bar V_1\cap Cl((\mathcal J_1\cdot \mathcal J_2)_0^c)=\emptyset$. $\bar V_1\cap (\mathcal J_1\cdot \mathcal J_2)_0^c=\emptyset$ and 
$(\mathcal J_1\cdot \mathcal J_2)_0^c$ is the union
of a finite set of subsets $\Spec \OO_{y,b}\subset \tilde U_y\subset \tilde X$ (with $b\in \OO_y$ and $y\in X$).   As we have proved above ($\overset*=$), $\bar V_1\cap 
Cl(\Spec \OO_{y,b})=\emptyset$ since  $\bar V_1\cap 
\Spec \OO_{y,b}=\emptyset$.
Then,  $\bar V_1\cap Cl((\mathcal J_1\cdot \mathcal J_2)_0^c)=\emptyset$.

2. Suppose that $V_1\cap V_2\neq \emptyset$. Let $\mathcal I:=\mathcal I_{\tilde X-\bar W}\subset \OO_X$. $\tilde X-\bar W\subset \tilde X$ is equal to $\tilde Y:=\Spec\, \OO_X/\mathcal I$.   $\bar V_1\cap \tilde Y=\Spec\, (\OO_X/\mathcal I)_{x_1,[a_1]}$ and $\bar V_2\cap \tilde Y=\Spec (\OO_X/\mathcal I)_{x_2,[a_2]}$.
By 1., there exist 
open subsets $\bar W'_1,\bar W'_2\subset \tilde Y$ such that
$\bar W'_1\cap \bar W'_2=\emptyset$, $\bar V_1\cap \tilde Y\subset \bar W'_1$ and $\bar V_2\cap \tilde Y\subset \bar W'_2$.
Then, $\bar W_1=\bar W\cup \bar W'_1$ and $\bar W_2=\bar W\cup \bar W'_2$ are the searched open subsets.

\edemo

\lema \label{L15.25}
 Let $X$ be a schematic finite space and  $B$ the family of quasi-compact open subsets of $\tilde X$.  Let $\mathcal F'$ be a presheaf on $\tilde X$ and $\mathcal F$ the sheafification of $\mathcal F'$. If for any $\bar V\in B$ and any finite open covering $\{\bar V_i\in B\}$ of $\bar V$ the sequence of morphisms
$$\mathcal F'(\bar V)\to \prod_{i} \mathcal F'(\bar V_i) \dosflechas  \prod_{ij}  \mathcal F'(\bar V_i\cap \bar V_j)$$
is exact, then $\mathcal F'(\bar V)=\mathcal F(\bar V)$.
\elema 

\demo It is well known.
\edemo

\coro Let $X$ be a schematic finite space and let $\{F_i\}_i$ be a direct system of sheaves of abelian groups on $\tilde X$. Then,
$$H^n(\tilde X,\ilim{i\in I} F_i)=\ilim{i\in I} H^n(\tilde X,F_i),$$
for any $n\geq 0$.\ecoro

\coro \label{coro25} Let $X$ be a schematic finite space, $U\subset X$ an open subset, $\bar V\subset \tilde U$ a quasi-compact open subset and 
$\mathcal F$ a sheaf of abelian groups on $\tilde X$. Then,
$$\mathcal F_{|\tilde U}(\bar V)=\ilim{ \bar V\subset\bar W} \mathcal F(\bar W).$$
Therefore, $H^n(\bar V,\mathcal F_{|\tilde U})=\ilim{ \bar V\subset\bar W} H^n(\bar W,\mathcal F)$, for any $n\geq 0$.

\ecoro

\demo 
Let $\mathcal G$ be the presheaf on $\tilde U$ defined by
$\mathcal G(\bar V):=\ilim{ \bar V\subset \bar W } \mathcal F(\bar W)$.  $\mathcal F_{|\tilde U}$ is the sheafification of $\mathcal G$. Let $\{\bar V_i\}$ a finite quasi-compact  open covering of $\bar V$.
Let $I_i$ be the family of open subsets $\bar W_i\subset \tilde X$ such that $\bar V_i\subset\bar W_i$ and $I=\prod_i I_i$. The sequence of morphisms
$$\mathcal F(\cup_i\bar W_i)\to \prod_{i} \mathcal F(\bar W_i) \dosflechas \prod_{i,j} \mathcal F(\bar W_i\cap \bar W_j)$$
is exact. Taking direct limits, we obtain the exact sequence
$$\ilim{(\bar W_i)\in I} \mathcal F(\cup_i\bar W_i)\to \ilim{(\bar W_i)\in I} \prod_{i} \mathcal F(\bar W_i) \dosflechas \ilim{(\bar W_i)\in I} \prod_{ij} \mathcal F(\bar W_i\cap \bar W_j)$$
Observe that $\ilim{(\bar W_i)\in I} \mathcal F(\cup_i\bar W_i)=\mathcal G(\cup_i \bar V_i)$, 
$\ilim{(\bar W_i)\in I} \prod_{i} \mathcal F(\bar W_i)=
\prod_{i} \mathcal G(\bar V_i)$, and by Lemma \ref{L16.18}, $\ilim{(\bar W_i)\in I} \prod_{ij} \mathcal F(\bar W_i\cap \bar W_j)=\prod_{ij} \mathcal G(\bar V_i\cap \bar V_j)$. Hence,
$$\mathcal G(\cup_i\bar V_i)\to \prod_{i} \mathcal G(\bar V_i) \dosflechas \prod_{i,j} \mathcal G(\bar V_i\cap \bar V_j)$$
is exact. By Lemma \ref{L15.25}, $\mathcal G(\bar V)=\mathcal F_{|\tilde U}(\bar V)$.

Finally, let $F\to C^\punto\mathcal F$ be the Godement resolution, then
$$H^n(\bar V,\mathcal F_{|\tilde U})=H^n\Gamma(\bar V, (C^\punto\mathcal F)_{|\tilde U})= \ilim{ \bar V\subset\bar W}H^n\Gamma(\bar W, C^\punto\mathcal F)=\ilim{ \bar V\subset\bar W} 
H^n(\bar W,\mathcal F).$$
\edemo

\teor \label{T16.22} Let $X$ be a schematic finite space and $U\subset X$ an open subset. Then, $${\OO_{\tilde X}}_{|\tilde U}=\OO_{\tilde U}.$$

Let $x\in X$ and $\pp\in \tilde U_x\subseteq \tilde X$. Then,
$\OO_{\tilde X,\pp}=\OO_{x,\pp}.$
\eteor

\demo Let $i\colon U\iny X$ be the inclusion morphism and $\tilde i\colon \tilde U\iny \tilde X$ the induced morphism. The natural morphism
$\OO_{\tilde X}\to \tilde i_*\OO_{\tilde U}$ defines by adjunction
the morphism ${\OO_{\tilde X}}_{|\tilde U}\to \OO_{\tilde U}$ and we have to prove that the morphism $\OO_{\tilde X,\pp}={\OO_{\tilde X}}_{|\tilde U,\pp} \to \OO_{\tilde U,\pp}$ is an isomorphism, for any $\pp\in\tilde U$.
Let $I:=\{(\bar W,\bar V)$, where $\bar V$ is any quasicompact open subset of $\tilde U$ such that $\pp\in \bar V$ and $\bar W$ is any open
subset of $\tilde X$ such that $\bar V\subset \bar W\}.$ Then,
$$\OO_{\tilde X,\pp}=\ilim{(\bar W,\bar V)\in I} \OO_{\tilde X}(\bar W)= \ilim{\pp\in \bar V}\ilim{\bar V\subset \bar W} \OO_{\tilde X}(\bar W) =\ilim{\pp\in \bar V}\OO_{\tilde U}(\bar V)=\OO_{\tilde U,\pp}.$$
Finally, $\OO_{\tilde X,\pp}=\OO_{\tilde U_x,\pp}=\OO_{x,\pp}$.

%
%
%
%it is an isomorphism. Let us proceed by induction on $\# X$.
%
%Suppose that $X$ is affine. Recall that $\tilde X=\Spec \OO(X)$
%and obviously ${\OO_{\tilde X}}_{|\tilde U_x} =\OO_{\tilde U_x}$.
%  Observe that 
%${{\OO_{\tilde X}}_{|\tilde U}}_{|\tilde U_x}=\OO_{\tilde U_x}$ and 
%${\OO_{\tilde U}}_{|\tilde U_x}=\OO_{\tilde U_x}$ by the induction hypothesis, then ${\OO_{\tilde X}}_{|\tilde U}= \OO_{\tilde U}$.
%
%Now in general. Let $\bar V\subset \tilde U$ be a quasi-compact open set.
%Then,
%
%$$\aligned {\OO_{\tilde X}}_{|\tilde U}&(\bar V)  =\ilim{\bar V\subset \bar W}
%\OO_{\tilde X}(\bar W)=\ilim{\bar W\cap \tilde U=\bar V} \plim{x\in X}
%\OO_{\tilde U_x}(\bar W\cap \tilde U_x)=
% \plim{x\in X} \ilim{\bar W\cap \tilde U=\bar V}
%\OO_{\tilde U_x}(\bar W\cap \tilde U_x)\\ &
%\underset*{\overset{\text{\ref{coro25}}}=}\plim{x\in X}{\OO_{\tilde U_x}}_{|\widetilde{U_x\cap  U}}(\bar V\cap \widetilde{U_x\cap U})
%=\plim{x\in X}
%\OO_{\widetilde{U_x\cap  U}}(\bar V\cap \widetilde{U_x\cap U})=
%\plim{x\in X} \plim{y\in U_x\cap  U}\OO_{\tilde U_y}(\bar V\cap \tilde U_y)\\ & = \plim{y\in U}\OO_{\tilde U_y}(\bar V\cap \tilde U_y)=
%\OO_{\tilde U}(\bar V).\endaligned$$
%($*$: Obviously $(\bar W\cap \tilde U_x)\cap \widetilde{U_x\cap U}=
%\bar V\cap \widetilde{U_x\cap U}$. If an open subset $\bar V_x\subset \tilde U_x$ satisfies $\bar V_x\cap \widetilde{U_x\cap U}=
%\bar V\cap \widetilde{U_x\cap U}$, let $\bar W'\subset \tilde X$
%be an open subset such that $\bar W'\cap (\tilde U\cup \tilde U_x)=\bar V\cup \bar V_x$, then $\bar W'\cap \tilde U=\bar V$ and $\bar W'\cap \tilde U_x=\bar V_x$).
\edemo

\section{$H^n(X,\M)=H^n(\tilde X,\tilde{\M})$}

\nota Given an affine scheme $\Spec R$ and an $R$-module $M$ we shall denote $\tilde M$ the sheaf of
localizations of the $R$-module $M$.

\enota

\defi Let $X$ be a schematic finite space and $\tilde X=\Spec\OO_X$. We shall say that an $\OO_{\tilde X}$-module $\bar \M$
is quasi-coherent if $\bar \M_{|\tilde U_x}$ is a quasi-coherent $\OO_{\tilde U_x}$-module for any $x\in X$.\edefi

I warn the reader that this definition is not the usual definition of quasi-coherent module.

Let $X$ be a schematic finite space and $\bar \M$ a quasi-coherent $\OO_{\tilde X}$-module. Let $\M$ be the $\OO_X$-module defined by $\M_x=\bar \M_{|\tilde U_x}(\tilde U_x)$, then it easy to check that $\M$ is a  quasi-coherent $\OO_X$-module (see \cite{Hartshorne} II 5.1 (d) and 5.2 c.) 

Let  $\M$ be a quasi-coherent $\OO_X$-module. Define
$\tilde\M:=\plim{x\in X} \tilde i_{x*}\widetilde{\M_x}$, where $\tilde i_x\colon \tilde U_x\to \tilde X$ is the morphism induced by the inclusion morphism $i_x\colon U_x\iny X$.
Observe that $\tilde \M(\tilde X)=\plim{x\in X} \M_x=\M(X)$.

\prop Let $X$ be an  affine finite space and $\M$ a quasi-coherent $\OO_X$-module. then $\tilde\M=\widetilde{\M(X)}$.\eprop

\demo Observe that
$\tilde i_{x*}\widetilde{\M_x}=\widetilde{\M_x}$, then  $$\tilde\M=\plim{x\in X} \tilde i_{x*}\widetilde{\M_x}=\plim{x\in X} \widetilde{\M_x}=\widetilde{\plim{x\in X}  \M_x}=\widetilde{\M(X)}.$$

\edemo

\prop \label{P16.17} Let $X$ be a schematic finite space, $U\subset X$ an open subset and $\M$ a quasi-coherent $\OO_X$-module. Then, 
$\tilde\M_{|\tilde U}=\widetilde{\M_{|U}}.$ In particular, $\tilde M$ is a quasi-coherent $\OO_{\tilde X}$-module.
\eprop

\demo Proceed like in the proof of Theorem \ref{T16.22}.
\edemo

Let $\M$ and $\M'$ be quasi-coherent $\OO_X$-modules. Any morphism of $\OO_X$-modules $ \M\to \M'$ induces a natural morphism
$\tilde \M = \plim{x\in X} \tilde i_{x*}\widetilde{\M_x}\to \plim{x\in X} \tilde i_{x*}\widetilde{\M'_x}=\tilde \M'$.  

Let $\bar \M$ and $\bar \Nc$ be quasi-coherent $\OO_{\tilde X}$-modules. Any morphism of $\OO_{\tilde X}$-modules $\bar \M\to \bar \Nc$ induces natural morphisms $\M_x:=\bar \M_{|\tilde U_x}(\tilde U_x)\to \bar \Nc_{|\tilde U_x}(\tilde U_x)=:\Nc_x$ and then a morphism
$\M\to \Nc$.

\medskip

\teor Let $X$ be a schematic finite space. The category of
quasi-coherent $\OO_X$-modules is equivalent to the category of  quasi-coherent $\OO_{\tilde X}$-modules.

\eteor

\demo The functors $\bar {\mathcal M} \functor \{\bar {\mathcal M}_{|\tilde U_x}(\tilde U_x)\}_{x\in X}$ and  $\mathcal M \functor 
\tilde{\mathcal M} $ are mutually inverse.\edemo

\prop \label{P15.6} Let $f\colon X \to Y$ be a schematic morphism and 
 $\tilde f\colon \tilde X\to \tilde Y$ the induced morphism. Let $\M$ be a quasi-coherent $\OO_X$-module and $\Nc$  a quasi-coherent $\OO_Y$-module.
Then, 
\enumera

\item $\tilde f_*\tilde \M=\widetilde{f_*\M}$.

\item $\tilde f^*\tilde \Nc=\widetilde{f^*\Nc}$.

\eenumera
\eprop

\demo Consider the obvious commutative diagram
$$\xymatrix{\tilde X \ar[r]^-{\tilde f} & \tilde Y & & \\
\tilde U_x \ar@{^{(}->}[u]^-{\tilde i_{x}} \ar[r]_-{\tilde f_{xy}} &
\tilde U_y \ar@{^{(}->}[u]_-{\tilde i_{y}} & &  (f(x)\geq y)}$$

1. Observe that $(f_*\M)(U_y)=\M(f^{-1}(U_y))=\plim{x\in f^{-1}(U_y)} \M_x$, then 
$$\aligned \widetilde{f_*\M} & =\plim{y\in Y}
\tilde i_{y*}(\plim{x\in f^{-1}(U_y)} \widetilde{ \M_x})=
\plim{y\in Y}
\tilde i_{y*}(\plim{x\in f^{-1}(U_y)} \tilde f_{xy*}\widetilde{\M_x})=
\plim{y\in Y}
\plim{x\in f^{-1}(U_y)} \tilde i_{y*}\tilde f_{xy*}\widetilde{\M_x}
\\ &=\plim{y\in Y}\plim{x\in f^{-1}(U_y)} \tilde f_*\tilde i_{x*}\widetilde{\M_x}=\tilde f_*\plim{y\in Y}\plim{x\in f^{-1}(U_y)} \tilde i_{x*}\widetilde{\M_x}=\tilde f_*\plim{x\in X} \tilde i_{x*}\widetilde{\M_x}=
\tilde f_*\tilde\M.\endaligned
$$

2. Observe that $(f^*\Nc)_x=\Nc_{f(x)}\otimes_{\OO_{f(x)}}\OO_x$, then $\widetilde{(f^*\Nc)}_{|\tilde U_x}\overset{\text{\ref{P16.17}}}=\widetilde{({f^*\Nc})_{|U_x}}=\widetilde{\Nc_{f(x)}\otimes_{\OO_{f(x)}}\OO_x}$.

On the other hand,  $(\tilde f^*\tilde\Nc)_{|\tilde U_x}=\tilde f_{xf(x)}^*(\tilde\Nc_{|\tilde U_{f(x)}})=\tilde f_{xf(x)}^*(\widetilde{\Nc_{f(x)}})=
\widetilde{\Nc_{f(x)}\otimes_{\OO_{f(x)}}\OO_x}$.

\edemo

\lema \label{L17.7} Let $X$ be a schematic finite space and $\mathcal F$ a sheaf of abelian groups on $X$. Let $U,V\subset X$ be two open subsets. Consider the obvious commutative diagram 
$$\xymatrix{\tilde U \ar@{^{(}->}[r]^-i & \tilde X\\ \widetilde{U\cap V}=\tilde U\cap \tilde V \ar@{^{(}->}[r]_-{\bar i} \ar@{^{(}->}[u]_-{\bar j}& \tilde V \ar@{^{(}->}[u]^-{j}} $$
Then, $j^*(R^ni_*\mathcal F)=R^n\bar i_*(\bar j^*\mathcal F)$, for any $n\geq 0$.
\elema

\demo Let $\pp\in \tilde V$. Then, $$\aligned (j^*(R^ni_*\mathcal F))_\pp & =
(R^ni_*\mathcal F)_\pp=\ilim{\pp\in \bar W\subset \tilde X} H^n(\bar W\cap \tilde U, \mathcal F_{|\tilde U})\overset{\text{\ref{coro25}}}=\ilim{\pp\in \bar W'\subset \tilde V} H^n(\bar W'\cap \tilde U\cap \tilde V, \mathcal F_{|\tilde U\cap \tilde W})\\  & =(R^n\bar i_*(\bar j^*\mathcal F))_\pp.\endaligned$$

\edemo

\teor Let $X$ be a semiseparated schematic finite space and 
$\mathcal M$ a quasi-coherent $\OO_X$-module. Then,
$$H^n(X,\mathcal M)=H^n(\tilde X,\tilde{\M}),$$ 
for any $n\geq 0$.

\eteor

\demo Let $\tilde i\colon \tilde U_x\iny \tilde X$ be the inclusion morphism. The morphism $\bar i\colon \tilde U_y\cap\tilde U_x\iny \tilde U_y,$ $\pp\mapsto i(\pp)$, is an affine morphism of schemes since  $\tilde U_x\cap \tilde U_y$ is an affine scheme because
$X$ is semiseparated. 
Let $\tilde\Nc$ be  a quasi-coherent $\OO_{\tilde U_x}$-module and $\pp\in \tilde U_y$. By Lemma \ref{L17.7},
$(R^n\tilde i_{*}\tilde\Nc)_\pp=(R^n\bar i_{*}\tilde\Nc_{|\tilde U_y\cap \tilde U_x})_\pp=0$, for any $n>0$.  
Hence, $R^n\tilde i_{*}\Nc=0$ and $H^n(
\tilde X,\tilde i_*\Nc)=H^n(\tilde U_x,\Nc)=0$, for any $n>0$. That is, 
$\tilde i_*\Nc$ is acyclic. 

Given a quasi-coherent $\OO_X$-module $\M$
denote $\tilde{\M}_{\tilde U_x}=\tilde i_*\tilde i\,^*\tilde\M$. Observe that $\tilde{\M}_{\tilde U_x}$ is acyclic and $\tilde{\M}_{\tilde U_x}(\tilde X)=\tilde \M_{|\tilde U_x}(\tilde U_x)=\widetilde{\M_{|U_x}} (\tilde U_x)=\M(U_x)=\M_x$. The obvious sequence of morphisms
$$\tilde\M \to \prod_{x\in X} \tilde \M_{\tilde U_x} \to  \prod_{x_1<x_2} \tilde \M_{\tilde U_{x_2}} \to  \prod_{x_1<x_2<x_3} \tilde \M_{\tilde U_{x_3}}\to \cdots$$
is exact. Denote this resolution  $\tilde\M\to \tilde C^\punto\tilde \M$ and let $\M\to C^\punto \M$ be the standard resolution of $\M$. Then,
$$H^n(\tilde X,\tilde\M)=H^n(\Gamma(\tilde X,\tilde C^\punto \tilde \M))=
H^n(\Gamma(X,C^\punto \M))=H^n(X,\M).$$

\edemo

\section{\rm Hom$_{sch}(\tilde X,\tilde Y)$ = Hom$_{[sch]}(X,Y)$}

\prop 
Let $f\colon X\to Y$ be a schematic morphism. Then, the induced morphism $\tilde f\colon \tilde X\to \tilde Y$ is quasi-compact, that is, $\tilde f^{-1}(\bar V)$ is quasi-compact for any quasi-compact open subset $\bar V\subset \tilde Y$
\eprop

\demo 
Any affine scheme morphism is quasi-compact. Given $x\in X$, denote
$\tilde f_x \colon \tilde U_x\to \tilde U_{f(x)}$ the morphism induced by $f_x\colon U_x\to U_{f(x)}$, $f_x(x'):=f(x')$. 
By  Proposition \ref{P16.16}.2., $\bar V\cap \tilde U_{f(x)}$ is quasi-compact. Then,
$$\tilde f^{-1}(\bar V)=\cup_{x\in X} \tilde f^{-1}(\bar V)\cap \tilde U_x=
\cup_{x\in X} \tilde f_x^{-1}(\bar V\cap \tilde U_{f(x)})$$
is quasi-compact.
\edemo

\prop \label{L18.8}

Let 
$f\colon X\to Y$  be a schematic morphism and
$\tilde f\colon \tilde X\to\tilde Y$ the induced morphism.
Then, $\tilde f\,^{-1}(\tilde U)=\widetilde{f^{-1}(U)}$, for any open subset $U\subset X$.
\eprop

\demo 
Given two open subsets $V,V'$ of a  schematic finite space, observe that $\tilde V\cap \tilde V'=
\widetilde{V\cap V'}$ and $\tilde V\cup \tilde V'=\widetilde{V\cup V'}$.

Given $y\geq f(x)$, 
$\Spec \prod_{z\in U_{xy}} \OO_z\to \Spec (\OO_x\otimes_{\OO_{f(x)}} \OO_y)=\Spec \OO_x\times_{\Spec \OO_{f(x)}} \Spec \OO_y$  is surjective, by Theorem \ref{K7}. 
Hence, $\tilde U_x\cap \tilde f^{-1}(\tilde U_y)=\cup_{z\in U_{xy}} \tilde U_z=\tilde{U_{xy}}=\!\!\!\!\widetilde{\,\,\,\,U_x\cap f^{-1}(U_y)}$.

Obviously, 
$ \widetilde{f^{-1}(U)}\subseteq \tilde f\,^{-1}(\tilde U)$.
 Let $\pp\in \tilde f\,^{-1}(\tilde U)$ and $x\in X$ such that $\pp\in \tilde U_x$. Then, $\tilde f(\pp)\in \tilde U_{f(x)}\cap \tilde U$. Let
$y\in U_{f(x)}\cap U$ such that $\tilde f(\pp)\in \tilde U_y$. Then, $\pp\in \tilde U_x\cap \tilde f\,^{-1}(\tilde U_y)=\widetilde{U_x\cap f\,^{-1}(U_y)}\subset \widetilde{f^{-1}(U)}$. Therefore, $\tilde f\,^{-1}(\tilde U)\subseteq \widetilde{f^{-1}(U)}$. 

\edemo

\defi 
Let $X$ and $Y$ be schematic finite spaces. We shall say
that a morphism of ringed spaces $f'\colon \tilde X\to \tilde Y$ is a schematic morphism if $f'_*\tilde \M$ is a quasi-coherent 
$\OO_{\tilde Y}$-module for any quasi-coherent $\OO_{\tilde X}$-module $\tilde\M$.
\edefi

\ejem
 If $f\colon X\to Y$ is a schematic morphism, then
$\tilde f\colon \tilde X\to \tilde Y$ is a schematic morphism, by Proposition \ref{P15.6}.
\eejem

\prop \label{P16.21} Let $f\colon \Spec B\to \Spec A$, $f'\colon \tilde A\to f_*\tilde B$  be a morphism of
ringed spaces. If $f_*\tilde M$ is a quasi-coherent $\tilde A$-module for any $B$-module $M$, then $(f,f')$ is a morphism of schemes.\eprop

\demo Let $f''\colon \Spec B\to \Spec A$ be the morphism defined on spectra by $f'_{\Spec A}\colon A\to B$. We only have to prove that $f=f''$. 

Let $\pp\in \Spec B$, $a\in A$ and $U_a=\Spec A\backslash (a)_0$. By the hypothesis, $f_*\widetilde{B/\pp}$ is a quasi-coherent $\tilde A$-module. Then, $$(B/\pp)_a=((f_*\widetilde{B/\pp})(\Spec A))_a=(f_*\widetilde{B/\pp})(U_a)=\widetilde{B/\pp}(f^{-1}(U_a)).$$
Then,
$$\aligned f'_{\Spec A}(a)\in\pp & \iff (B/\pp)_a=0  \iff \widetilde{B/\pp}(f^{-1}(U_a))=0\iff f^{-1}(U_a)\cap (\pp)_0=\emptyset\\ & \iff \pp\notin f^{-1}(U_a) \iff
f(\pp)\notin U_a \iff a\in f(\pp).\endaligned$$
Therefore, ${f'_{\Spec A}}^{-1}(\pp)=f(\pp)$, that is to say,
$f''=f$.

\edemo

\prop Any schematic morphism $f'\colon \tilde X\to \tilde Y$ is
a morphism of locally ringed spaces.\eprop

\demo Let $\pp\in\tilde U_x\subset \tilde X$. Let $g$ be the composition of the schematic  morphisms  
$$\Spec \OO_{x,\pp}  \iny \Spec \OO_x\iny \tilde X\to \tilde Y.$$
 Let $y\in Y$  be a point  such that $f(\pp)\in \tilde U_y$, then $g^{-1}(\tilde U_y)=\Spec \OO_{x,\pp}$. Consider the continuous morphism
$h\colon \Spec \OO_{x,\pp}\to \tilde U_y$, $\qq\mapsto g(\qq)$ and let $\tilde i_y\colon \tilde U_y \iny \tilde Y$ be the inclusion morphism. Consider the morphism $\OO_{\tilde Y}\to g_*\tilde \OO_{x,\pp}$. Taking $i^*$, we obtain a morphism $\phi\colon \OO_{\tilde U_y}\to h_*\tilde\OO_{x,\pp}$. The morphism of ringed spaces $(h,\phi)$
is schematic, since $h_*\M=\tilde i_y^*\tilde i_{y*}h_*\M=\tilde i_y^*g_*\M$ is quasi-coherent, for any quasi-coherent $\OO_{\tilde U_y}$-module $\M$. By Proposition \ref{P16.21}, $h$ is a morphism of locally ringed spaces. We are done.

\edemo 

\lema \label{L17.11} Let $\M$ be a finitely generated $\OO_{X}$-module.
For any $\pp\in \tilde  X$, there exist an open neighbourhood $\bar U$ of $\pp$ and an epimorphism of sheaves $\OO^n_{\bar U} \to \tilde\M_{|\bar U}$.
\elema

\demo $\bar V:=\{\qq\in \tilde X\colon \tilde \M_{\qq}=0\}$ is an open subset of $\tilde X$, since $\bar V\cap  \tilde U_x=\{\qq\in \tilde U_x\colon {\M_x}_{\qq}=0\}$ is an open subset of $\tilde U_x$. Hence, given 
another quasi-coherent module  $\tilde \M'$, and a  morphism of $\OO_{\tilde X}$-modules $\tilde\M'\to \tilde\M$, if
the morphism on stalks at $\pp$, $ \tilde \M'_\pp\to \tilde\M_\pp$ is an epimorphism then there exists an open neighbourhood $\bar V$ of $\pp$
such that the morphism of sheaves $\tilde\M'_{|\bar V}\to \tilde\M_{|\bar V}$ is an epimorphism.

Let $m_1,\ldots,m_r$ be a generator system of the $\OO_\pp$-module  $\tilde \M_\pp$.
Let $\bar W\subset \tilde X$ be an open neighbourhood of $\pp$, such that there exist $n_1,\ldots,n_r\in \tilde\M(\bar W)$ satisfying $n_{i,\pp}=m_i$. The morphism of $\OO_{\bar W}$-modules $\OO_{\bar W}^r \to \tilde\M_{|\bar W}$, $(a_i)\mapsto \sum_i a_i\cdot n_i$ is an epimorphism
on stalks at $\pp$. Hence, it is an epimorphism in an open neighbourhood $\bar U$ of $\pp$.
\edemo

\prop Any schematic morphism $f'\colon \tilde X\to \tilde Y$ is
quasi-compact.
\eprop

\demo Let $\bar V\subset \tilde Y$ be a quasi-compact  open subset, we have to prove that $f'^{-1}(\bar V)$ is quasi-compact.
We only have to prove that $\tilde U_x\cap f'^{-1}(\bar V)$ is quasi-compact, for any $x\in X$. Hence, we can suppose that $\tilde X=\tilde U_x$.

Let $C:=\tilde Y-\bar V$. By Lemma \ref{L15.20}, there exists a finitely generated ideal $\mathcal I\subset \OO_Y$ such that $(\mathcal I)_0=C$. Consider the exact sequence of morphisms 
$0\to \mathcal I \to \OO_{\tilde Y}\to \OO_{\tilde Y}/\mathcal I\to 0$.
By Lemma \ref{L17.11}, there exist an open covering 
$\{\bar U_i\}$ of $\tilde Y$ and epimorphisms $\OO_{\bar U_i}^{n_i}\to \mathcal I_{|\bar U_i}$.
Taking $f^*$, one has an exact sequence of morphisms 
$$f'^*\mathcal I\to \tilde\OO_x \to f'^*(\OO_{\tilde Y}/\mathcal I)\to 0$$
and $\mathcal J:=\Ima[f'^*\mathcal I\to \tilde\OO_x ]$ is a finitely generated  quasi-coherent ideal of $\tilde \OO_x$, since it is so locally (over $f'^{-1}(\bar U_i)$). 
$(\mathcal J)_0=f'^{-1}(C)$, therefore $f'^{-1}(\bar V)=\tilde X-(\mathcal J)_0$ is a quasi-compact open subset.
\edemo

Let $\C_{sch}$ be the category of schematic finite spaces and $W$ the family of quasi-isomorphisms. Let us construct the localization of  $\C_{sch}$ by $W$, $\C_{sch}[W^{-1}]$.

\defi A schematic pair of morphisms from $X$ to $Y$, $X\,\overset{f}{-->} Y$  is a pair of schematic morphisms $(\phi',f')$
$$\xymatrix{ & X' \ar[ld]_-{\phi'} \ar[rd]^-{f'} & \\ X & & Y}$$
where $\phi'$ is a quasi-isomorphism.
\edefi

\ejem A schematic morphism $f\colon X\to Y$ can be considered as a schematic pair of morphisms: Consider the pair of morphisms $(Id_X,f)$
$$\xymatrix{ & X \ar[ld]_-{Id_X} \ar[rd]^-{f} & \\ X & & Y}$$
\eejem 

\defi Let  $f=(\phi',f')\colon X --> Y$ y
$g=(\varphi',g')\colon Y --> Z$ be two schematic pairs of morphisms, where $\phi'\colon X'\to X$ is a quasi-isomorphism and $f'\colon X'\to Y$ is a schematic morphism and $\varphi'\colon Y'\to Y$ is a quasi-isomorphism and $g'\colon Y'\to Z$ is a schematic morphism.
Let $\pi_1,\pi_2\colon X'\times_Y Y'\to X',Y'$ be the two obvious projection maps (observe that $\pi_1$ is a quasi-isomorphism). We define
$g\circ f:=(\phi'\circ\pi_1,g'\circ\pi_2)\colon X --> Z$
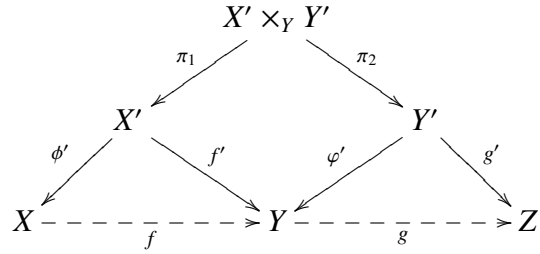
$$\xymatrix{ && X'\times_Y Y' \ar[ld]_-{\pi_1}  \ar[rd]^-{\pi_2}&& \\& X' \ar[ld]_-{\phi'} \ar[rd]^-{f'} &  & Y' \ar[ld]_-{\varphi'} \ar[rd]^-{g'} &\\ X \ar@{-->}[rr]_-{f} & & Y\ar@{-->}[rr]_-{g}   & & Z}$$
\edefi

Let $f\colon X\to Y$ and $g\colon Y\to Z$ be schematic morphisms. Then, 

$$(Id_X,f)\circ (Id_Y,g)=(Id_X,f\circ g).$$

\defi \label{d5} Two schematic pairs of morphisms $(\phi',f'),(\phi'',f'')\colon X-->Y$
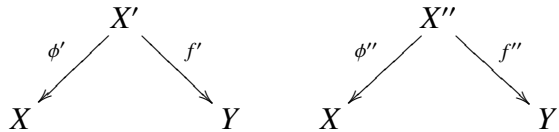
$$\xymatrix{ & X' \ar[ld]_-{\phi'} \ar[rd]^-{f'} & & & X'' \ar[ld]_-{\phi''} \ar[rd]^-{f''} &\\ X & & Y &  X & & Y}$$
are said to be equivalent, $(\phi',f')\equiv (\phi'',f'')$, if there exist a schematic space $T$ and two quasi-isomorphisms $\pi'\colon T\to X'$, $\pi''\colon T\to X''$
such that the diagram  
$$
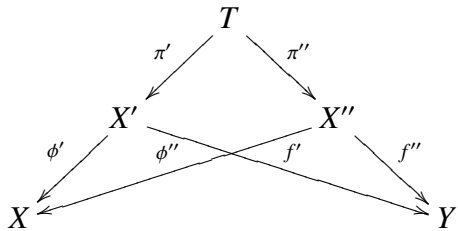
\xymatrix{  & & T \ar[dl]_-{\pi'} \ar[dr]^-{\pi''} &&\\
& X' \ar[dl]_-{\phi'} \ar[drrr]^-{f'} & & X'' \ar[dlll]_-{\phi''} \ar[rd]^-{f''} & \\ X & & & & Y
}$$
is commutative. \edefi

In order to prove the associative property the reader should consider the following commutative diagram (where the double arrows are quasi-isomorphisms)
$$
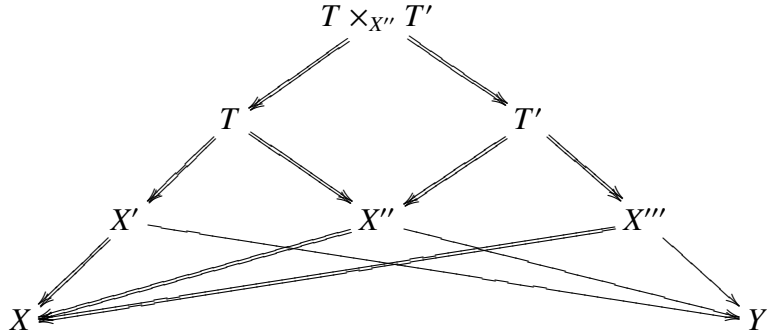
\xymatrix{ & & & T\times_{X''} T'  \ar@<0.2ex>[rd] \ar[rd]   \ar@<0.2ex>[ld]\ar[ld] &&& \\ && T   \ar@<0.2ex>[rd]\ar[rd]  \ar@<0.2ex>[ld]\ar[ld] & & T'  \ar@<0.2ex>[rd]\ar[rd]  \ar@<0.2ex>[ld]\ar[ld] & & \\  & X'  \ar@<0.2ex>[ld]\ar[ld]  \ar[rrrrrd]& & X''  \ar[rrrd]  \ar@<0.2ex>[llld]\ar[llld] & & X'''  \ar[rd]  \ar@<0.2ex>[llllld]\ar[llllld] & \\ X &  &  & &  & & Y }$$

\nada[Definitions and notations] Let $f=(\phi,f')\colon X --> Y$ be a schematic pair of morphisms.
The equivalence class of $f$ (resp. $(\phi,f')$) will be denoted $[f]$ (resp. $[\phi,f']$). We shall say that $[f]$ (or  $[f]\colon X\to Y$) is a [schematic] morphism from $X$ to $Y$.

Let $f\colon X\to Y$ be a schematic morphism. The equivalence class of $(Id,f)$
will be denoted $[f]$. \enada

\prop \label{Prop6} Let $(\phi',f')\colon X --> Y$ be  a schematic pair of morphisms, where $\phi'\colon X'\to X$ is a quasi-isomorphism and $f'\colon X'\to Y$ is a schematic morphism. Let $\varphi\colon X''\to X'$ be a quasi-isomorphism.
Then, $[\phi',f']=[\phi'\circ \varphi,f'\circ \varphi]$.\eprop

\demo Consider the commutative  diagram

$$
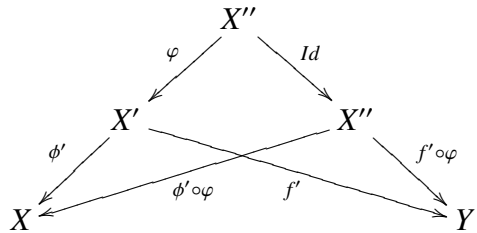
\xymatrix{  & & X'' \ar[dl]_-{\varphi} \ar[dr]^-{Id} &&\\
& X' \ar[dl]_-{\phi'} \ar[drrr]_-{f'} & & X'' \ar[dlll]^-{\phi'\circ \varphi} \ar[rd]^-{f'\circ \varphi} & \\ X & & & & Y
}$$
\edemo

%In Definition \ref{d5}, we can suppose $T$ is minimal. The morphism
%$T\to X'\times_X X'',$ $t\mapsto (\pi'(t),\pp''(t)$ is a qc-isomophism, because the morphisms $T\to X'$ and $X'\times_X X''\to X'$ are equivalences. Then, the morphism $T\to (X'\times_X X'')_M$
%is a quasi-isomorphism. The 

\teor Let $f,F\colon X -->  Y$ and $g,G\colon Y --> Z$ be schematic pairs of morphisms. If
$[f]=[F]$ and  $[g]=[G]$,  then  $[g\circ f]=[G\circ F]$.

\eteor

\demo 1. Let us prove that $[g\circ f]= [g\circ F]$.
Write $f=(\phi,f')$, $\phi\colon X'\to X$, $f'\colon X'\to Y$.
Let $\varphi\colon X''\to X'$ be a quasi-isomorphism and let $F':=(\phi\circ \varphi, f'\circ \varphi)$. By Proposition \ref{Prop6}, $[f]=[F']$. Consider the commutative diagram
$$\xymatrix{ X'' \ar@<0.2ex>[d]\ar[d]_\varphi & X''\times_Y Y' \ar@<0.2ex>[d]\ar[d] \ar@<0.2ex>[l]\ar[l]&  & \\ X' \ar[rrd]_{f'} \ar@<0.2ex>[d]\ar[d]_\phi & X'\times_Y Y' \ar@<0.2ex>[l]\ar[l] \ar[r]  & Y' \ar@<0.2ex>[d]\ar[d] \ar[rd] & \\
X \ar@{-->}[rr]_{f,F'} & & Y\ar@{-->}[r]_g & Z}$$
By Proposition \ref{Prop6}, $[g\circ f]=[ g\circ F']$. Finally, since $[f]=[ F]$, there exists $F'$ such that $[g\circ f]=[g\circ F']=[ g\circ F]$.

2. Likewise, $[g\circ f]=[ G\circ f]$.

3. The theorem is a consequence of 1. and 2.
\edemo

Let $f\colon X--> Y$ and
$g\colon Y--> Z$ be two schematic pairs of morphisms. We define
$[g]\circ [f]:=[g\circ f].$

\prop \label{inv} Let $(\phi,f')\colon X --> Y$ be a schematic pair of morphisms. Then, 
\enumera 
\item If $f'$ is a quasi-isomorphism, then $[\phi,f']^{-1}=[f',\phi]$.

\item  $[\phi,f']=[f']\circ [\phi]^{-1}$.
\eenumera

\eprop

\demo We have the quasi-isomorphism $\phi\colon X'\to X$ and the morphism $f'\colon X' \to Y$.
Let $$\pi_1,\pi_2\colon X'\times_XX'\to X', \, \pi_1(x_1',x'_2):=x'_1,\, \pi_2(x_1',x'_2):=x'_2,$$ which are quasi-isomorphisms. 

1. Let $\delta\colon X'\to X'\times_X X'$ be the diagonal morphism, which is a quasi-isomorphism because $\pi_1$ and $\pi_1\circ \delta=Id$ are quasi-isomorphisms. Then,
$$[\phi,f']\circ [f',\phi]=[f'\circ \pi_1,f'\circ \pi_2]\overset{\text{\ref{Prop6}}}=[f'\circ \pi_1\circ\delta,f'\circ \pi_1\circ\delta]=[f',f']\overset{\text{\ref{Prop6}}}=[Id_Y,Id_Y]=[Id_Y]$$
Likewise, $[f',\phi]\circ [\phi,f']=[Id_X]$.

2. It is easy to check that
$[f']\circ [\phi]^{-1}=[Id_{X'}, f']\circ [\phi,Id_{X'}]=[\phi,f'].$

\edemo

\prop \label{paque} Let  $X$ be a minimal schematic finite space. Let $f,g\colon X\to Y$ be two schematic morphisms.
Then, $[f]=[g]$ iff $f=g$.\eprop

\demo $\Rightarrow)$ There exists a (surjective) quasi-isomorphism 
 $\pi\colon T\to X$ such that $f\circ \pi=g\circ \pi$. Then $f$ and $g$ are equal as continuous maps. Finally, the morphism  $\OO_Y\to f_*\pi_*\OO_T=f_*\OO_X$
coincides with the morphism $\OO_Y\to g_*\pi_*\OO_T=g_*\OO_X$.

\edemo

\prop \label{paque+} Let $X$ be a minimal schematic finite space. Let  $f,g\colon X\to Y$ be two schematic morphisms and $\pi\colon Y\to Y'$ a quasi-isomorphism. Then,  $\pi\circ f=\pi\circ g$ iff $f=g$.\eprop

\demo $\Rightarrow)$ Observe that $[f]=[g]$ since $[\pi]\circ [f]=[\pi\circ f]=[\pi\circ g]=[\pi]\circ [g]$. By Proposition \ref{paque}, $f=g$.\edemo

\prop \label{P14.11} A [schematic] morphism $[\phi,f]\colon X \to Y$ is an isomorphism iff $f$ is a quasi-isomorphism.\eprop

\demo $\Leftarrow)$ $[\phi,f]^{-1}=[f,\phi]$, by Proposition \ref{inv}.

$\Rightarrow)$ $[\phi,f]=[f]\circ [\phi]^{-1}$ is invertible, then $[f]$ is invertible. 
Put $f\colon Z\to Y$ and 
let $[\varphi,g]\colon Y \to Z$ be the inverse morphism of $[f]$, where
$\varphi \colon T\to Y$ is a quasi-isomorphism (and we can assume that $T$ is minimal) and 
$g\colon T\to Z$ is a schematic morphism. Then, $[\Id_Y]=[f]\circ [\varphi,g]=[f]\circ [g]\circ [\varphi]^{-1}$ and $[\varphi]=[f\circ g]$. Hence, $\varphi=f\circ g$, by Proposition \ref{paque}.
Besides, $[\Id_Z]=[\varphi,g]\circ [f]$. If we consider the commutative diagram
$$\xymatrix{Z\times_YT \ar[r]^-{\pi_2} \ar[d]_-{\pi_1} & T \ar[d]^-{\varphi} \ar[rd]^-g & \\
Z \ar[r]_-f & Y & Z}$$
then $[\Id_Z]=[\varphi,g]\circ [f]=[\pi_1,g\circ \pi_2]$. Let $i\colon (Z\times_YT)_M\subseteq Z\times_Y T$ be the natural inclusion, $\pi_1';=\pi_1\circ i$ and $\pi_2'=\pi_2\circ i$. Then,
$[\pi'_1,g\circ \pi'_2]=[\Id_Z]$ and $[\pi'_1]=[g\circ \pi'_2]$.
By Proposition \ref{paque},
$\pi'_1=g\circ \pi'_2$.
Let $\tilde g\colon T \to  (Z\times_YT)_M$, $\tilde g(t)=(g(t),t)$. Observe that
$$\pi'_1\circ \tilde g \circ \pi'_2=g\circ \pi'_2=\pi'_1.$$
By Proposition \ref{paque+}, $\tilde g\circ \pi_2'=\Id_{(Z\times_YT)_M}$. Obviously, $\pi_2'\circ \tilde g=\Id_T$. Hence, $\pi_2'$ is an isomorphism and  $\pi_2$ a quasi-isomorphism. Finally, $f$ is a quasi-isomorphism since $\pi_1$, $\phi$ and $\pi_2$ are  quasi-isomorphisms and $f\circ \pi_1=\phi\circ \pi_2$.

\edemo

Given a [schematic] morphism $g=[\phi,f]\colon X \to Y$, consider the functors
$$\alinea g_*\colon\, {\bf Qc\text{-}Mod}_X \,\to \,{\bf Qc\text{-}Mod}_Y,\,\, g_*\mathcal M:=f_*\phi^*\mathcal M\\ g*\colon \,{\bf Qc\text{-}Mod}_Y\, \to  \,{\bf Qc\text{-}Mod}_X,\,\, g^*\mathcal M:=\phi_*f^*\mathcal M\ealinea$$ Recall $[\phi,f]=[\phi\circ\varphi,f\circ \varphi]$, where $\varphi$ is a quasi-isomorphism. Then we have canonical isomorphisms $f_*\phi^*\mathcal M=f_*\varphi_*\varphi^*\phi^*\mathcal M=(f\circ \varphi)_*(\phi\circ \varphi)^*\mathcal M$ and 
$\phi_*f^*\mathcal M=\phi_*\varphi_*\varphi^*f^*\mathcal M=(\phi\circ \varphi)_*(f\circ \varphi)^*\mathcal M$. 

\prop Let $g=[\phi,f]\colon X\to Y$ be a [schematic] morphism. The functors $g_*$ and $g^*$ are mutually inverse iff $g$ is a [schematic] isomorphism.\eprop

\demo 
$\Rightarrow)$ If $\Id=g_*g^*$ and $\Id=g^*g_*$, then $\Id=f_*\phi^*\phi_*f^*=f_*f^*$
and $\Id=\phi_*f^*f_*\phi^*$, hence $f^*f_*=\phi^*\phi_*=\Id$. By Theorem \ref{T12.7}, $f$ is a quasi-isomorphism and $g$ is invertible.

$\Leftarrow)$ By Proposition \ref{P14.11}, $f$ is a quasi-isomorphism, then $g_*g^*=f_*\phi^*\phi_*f^*=f_*f^*=\Id$ and $g^*g_*=\phi_*f^*f_*\phi^*=\phi_*\phi^*=\Id$.
\edemo

\nota Let $X$ and $Y$ be two schematic finite spaces. 
$\Hom_{[sch]}(X,Y)$ will denote the family  of [schematic] morphisms from $X$ to $Y$. 
$\Hom_{sch}(\tilde X,\tilde Y)$ will denote the set  of schematic morphisms from $\tilde X$ to $\tilde Y$. \enota

\lema \label{L17.2} Let $g\colon Y'\to Y$ be a [schematic] isomorphism and $X$ a schematic finite space. Then, the maps
$$\alinea \Hom_{[sch]}(X,Y')\to \Hom_{[sch]}(X,Y), \,\,[f]\mapsto [g]\circ [f]\\
\Hom_{[sch]}(Y,X)\to \Hom_{[sch]}(Y',X), \,\,[f]\mapsto [f]\circ [g]\ealinea$$
are biyective.
\elema

\prop Let $X$ be a schematic finite space and $Y$ an affine finite space. Then,
$$\Hom_{[sch]}(X,Y)=\Hom_{rings}(\OO(Y), \OO(X)).$$
\eprop

\demo  For any schematic finite space $T$,
$\Hom_{sch}(T,(*,A))=\Hom_{rings}(A,\OO(T)).$ 

Consider the natural morphism $\pi\colon Y\to (*,\OO(Y))$, which is a quasi-isomorphism. Then, 
$$\Hom_{[sch]}(X,Y)=\Hom_{[sch]}(X,(*,\OO(Y))=\Hom_{rings}(\OO(Y),\OO(X)).$$\edemo

Let  $[\phi,f]\colon X\to Y$ be a [schematic] morphism, where
$\phi\colon X'\to X$ is a quasi-isomorphism and $f\colon X'\to Y$ a schematic morphism. Consider the morphisms
$$\xymatrix{\tilde{X'} \ar@{=}[d]^-{\tilde \phi} \ar[rd]^-{\tilde f}& \\ \tilde X & \tilde Y}$$
where $\tilde\phi$ is an isomorphism, by Proposition \ref{P16.4}. The morphism
$$\Hom_{[sch]}(X,Y) \to \Hom_{sch}(\tilde X,\tilde Y),\,\, [\phi,f]\mapsto \tilde f\circ \tilde\phi^{-1}$$ is well defined.

\lema \label{L17.4} Let $X$ be a minimal schematic space, $Y$ a schematic $T_0$-space, 
$f\colon X\to Y$  a schematic morphism and
$\tilde f\colon \tilde X\to\tilde Y$ the induced morphism. Given $x\in X$, $y=f(x)$  iff $y$ is the greatest element of $Y$ such that $\tilde f(\tilde U_x)\subseteq  \tilde U_{y}$.
\elema

\demo Obviiously  $f(U_x)\subset U_{f(x)}$ and  $\tilde f(\tilde U_x)\subset \tilde U_{f(x)}$.
If $\tilde f(\tilde U_x)\subseteq \tilde U_{y'}$, then
$\tilde U_x\subset \tilde f^{-1}(\tilde U_{y'})\overset{\text{\ref{L18.8}}}=\widetilde{f^{-1}(U_{y'})}$
Hence,  $\tilde U_x=\tilde U_x\cap \widetilde{f^{-1}(U_{y'})}=\widetilde{U_x\cap f^{-1}(U_{y'})}$. Therefore, $x\in U_x\cap f^{-1}(U_{y'})$, since $x$ is not a removable point. That is, $f(x)\in U_{y'}$ and $f(x)\geq y'$.
\edemo

\prop \label{P17.5} Let $X$ and $Y$ be schematic finite spaces. The natural morphism 
$$\Hom_{[sch]}(X,Y)\to \Hom_{sch}(\tilde X,\tilde Y),\,[\phi,f]\mapsto \tilde f\circ \tilde\phi^{-1}$$
is injective.
\eprop

\demo Let $[\phi,f],[\phi',f']\colon X\to Y$ be [schematic] morphisms such that $\tilde f\circ \tilde\phi^{-1}=\tilde f'\circ \tilde\phi'^{-1}$. We can suppose that $\phi=\phi'$, then $\tilde f
=\tilde f'$. Let us say that $f$ and $f'$ are morphisms from $X'$ to $Y$. By Proposition \ref{P16.4} and  Lemma \ref{L17.2}  we can  suppose that $X'$ and $Y$ are minimal schematic spaces.

By Lemma \ref{L17.4}, the map $f$ is determined by $\tilde f$. The morphism of rings $\OO_{f(x')}\to \OO_{x'}$ is determined by the morphism of schemes $\Spec \OO_{x'}\to \Spec \OO_{f(x')}$. Therefore, $f=f'$.

\edemo 

\defi Let $U\overset{i_1}\to U_1$ and $U\overset{i_2} \to U_2$ be quasi-open immersions. We denote $U_1\cupa{U} U_2:= C(i_1)\coproda{U} C(i_2)$.
\edefi 

\vaci \label{ll} Observe that $C(i_1)$ and $C(i_2)$ are open subsets of $U_1\cupa{U} U_2$, $C(i_1)\cup C(i_2)= U_1\cupa{U} U_2$, $C(i_1)\cap C(i_2)= U$
and the natural morphisms $C(i_j)\to U_j$ are quasi-isomorphisms, for $j=1,2$. \evaci

Let $U_1,U_2\subset X$ be open subsets. Then, the natural morphism
$U_1\cupa{U_1\cap U_2} U_2\to U_1\cup U_2$ is a quasi-isomorphism.

Let 
$$\xymatrix @R=8pt { & V_1 \ar[rr]^-{\tilde{}}& & U_1\\ V \ar[ur] \ar[dr] \ar[rr]^-{\tilde{}} &  & U \ar[ru] \ar[rd] & \\ & V_2 \ar[rr]^-{\tilde{}}& & U_2}
$$
be a commutative diagram of quasi-open immersions, where the arrows $\tilde\longrightarrow$ are quasi-isomorphisms. Then, the natural morphism
$V_1\cupa{V} V_2\to U_1\cupa{U} U_2$ is a quasi-isomorphism.

\teor \label{T18.14} Let $U\overset{i_1}\to U_1$ and $U\overset{i_2}\to U_2$ be quasi-open immersions. Then,
  $$\Hom_{[sch]}(U_1\cupa{U} U_2,Y)=\Hom_{[sch]}(U_1,Y)\times_{\Hom_{[sch]}(U,Y)}
 \Hom_{[sch]}(U_2,Y).$$
 In other words (by \ref{ll}), let $U_1,U_2\subset X$ be open subsets. Then,
  $$\Hom_{[sch]}(U_1\cup U_2,Y)=\Hom_{[sch]}(U_1,Y)\times_{\Hom_{[sch]}(U_1\cap U_2,Y)}
 \Hom_{[sch]}(U_2,Y).$$
 \eteor
  
\demo Let $U_1\overset{j_1}\iny U_1\cup U_2$ and $U_2\overset{j_2}\iny U_1\cup U_2$ be the obvious inclusion morphisms . We have to prove that
  
$$\alineas{rcl} \Hom_{[sch]}(U_1\cup U_2,Y) & \to & \Hom_{[sch]}(U_1,Y)\times_{\Hom_{[sch]}(U_1\cap U_2,Y)}\Hom_{[sch]}(U_2,Y)
\\  \left[f\right] & \mapsto & ([f]\circ [j_1],[f]\circ [j_2])\ealineas$$
is bijective.  

Let $[f],[g]$ such that  $([f]\circ [j_1],[f]\circ [j_2])=([g]\circ [j_1],[g]\circ [j_2])$. There exist a minimal schematic finite space $W$,
a quasi-isomorphism $\phi \colon W\to U_1\cup U_2$ and morphisms $f',g'\colon W\to Y$ such that $[f]=[\phi ,f']$ and $[g]=[\phi ,g']$. Then,
$$[f]\circ [j_1]=[\phi _{|\phi ^{-1}(U_1)},f'_{|\phi ^{-1}(U_1)}] \text{ and }
[g]\circ [j_1]=[\phi _{|\phi ^{-1}(U_1)},g'_{|\phi ^{-1}(U_1)}].$$
By Proposition \ref{paque}, $f'_{|\phi ^{-1}(U_1)}=g'_{|\phi ^{-1}(U_1)}$. Likewise, $f'_{|\phi ^{-1}(U_2)}=g'_{|\phi ^{-1}(U_2)}$. That is, $f'=g'$ and $[f]=[g]$.

Let $([f_1],[f_2])\in \Hom_{[sch]}(U_1,Y)\times_{\Hom_{[sch]}(U,Y)}\Hom_{[sch]}(U_2,Y)$. Write $[f_1]=[\phi_1,g_1]$ and $[f_2]=[\phi_2,g_2]$. Since $[f_1]\circ [i_1]=[f_2]\circ [i_1]$, there exist a schematic finite space $V$ and a commutative diagram
$$\xymatrix @R=10pt @C=32pt{ & V_1\times_{U_1} U_1\cap U_2 \ar@{.>}[dd] \ar[r] & V_1  \ar[rdd]^-{g_1} \ar@{.>}[d]_-{\phi_1} & 
\\  & & U_1 & 
\\V  \ar@{.>}[ruu] \ar@{.>}[rdd] &   U_1\cap U_2 \ar[ru]^-{i_1}  \ar[rd]_-{i_2} & & Y
\\ &  & U_2 & 
\\&  V_2\times_{U_2} U_1\cap U_2 \ar@{.>}[uu] \ar[r] & V_2 \ar@{.>}[u]^-{\phi_2}   \ar[ruu]_-{g_2} &  }$$
(where the arrows $\cdots\!\!>$ are quasi-isomorphisms). Then, we have a schematic morphism $g\colon V_1\cupa{V} V_2\to Y$ and the composition $\phi$ of the quasi-isomorphisms $V_1\cupa{V} V_2\to  U_1\cupa{U_1\cap U_2} U_2\to U_1\cup U_2$. The reader can check that $[\phi,g]\mapsto ([f_1],[f_2])$.

\edemo

%\prop Let $\{U_1,U_2\}$ be an open covering of  $X$ and 
%let $f_1\colon U_1--> Y$, $f_2\colon U_2--> Y$ be two local schematic  morphisms such that ${f_1}_{|U_1\cap U_2}\equiv {f_2}_{|U_1\cap U_2}$.  Then, there exists a local schematic  morphism
% $f\colon X-->Y$ such that 
%$f_{|U_1}\equiv f_1$ and $f_{|U_2}\equiv f_2$.\eprop
%
%
%\demo Let $f_1=(\phi_1,f'_1)$, where $\phi'_1\colon T_1\to U_1$ is a quasi-isomorphism and $f'_1\colon T_1\to Y$ is a schematic morphism.  
%Let $f_2=(\phi_2,f'_2)$, where $\phi'_2\colon T_2\to U_2$ is a quasi-isomorphism and $f'_2\colon T_2\to Y$ is a schematic morphism. Let
%$T_{12}=f_1^{-1}(U_1\cap U_2)$ and $T_{21}=f_2^{-1}(U_1\cap U_2)$. By the hypothesis there exist quasi-isomorphisms $g_1\colon T\to T_{12}$ and $g_2\colon T\to T_{21}$, such that ${f'_1}_{|T_{12}}\circ g_1=
%{f'_2}_{|T_{21}}\circ g_2$ and ${\phi'_1}_{|T_{12}}\circ g_1=
%{\phi'_2}_{|T_{21}}\circ g_2$. Let $$S:=[T\coprod (T_1-T_{12})]\coprod_T 
%[T\coprod (T_2-T_{21})]$$
%Let $\phi\colon S\to X$ be defined by $\phi(t):=\phi_1(g_1(t))=\phi_2(g_2(t))$, for any $t\in T$, $\phi(t_1):=\phi_1(t_1)$ for any $t_1\in T_1-T_{12}$, 
%$\phi(t_2):=\phi_2(t_2)$ for any $t_2\in T_2-T_{21}$.  
%Let $f'\colon S\to Y$ be defined by $f'(t):=f'_1(g_1(t))=f'_2(g_2(t))$, for any $t\in T$, $f'(t_1):=g_1(t_1)$ for any $t_1\in T_1-T_{12}$, 
%$f'(t_2):=g_2(t_2)$ for any $t_2\in T_2-T_{21}$.  
%Then, $f:=(\phi,f')$.
%
%
%\edemo

\teor 
Let $X$ and $Y$ be schematic finite spaces  and suppose that $Y$ is semiseparated. Then, the morphism 
$$\Hom_{[sch]}(X,Y)\to \Hom_{sch}(\tilde X,\tilde Y),\,[\phi,f]\mapsto \tilde f\circ \tilde\phi^{-1}$$
is bijective.
\eteor

\demo By Proposition \ref{P17.5}, it is injective. Let $f'\in \Hom_{sch}(\tilde X,\tilde Y)$.

1. Assume that $\tilde X=\Spec A$ is affine. We can suppose that $Y$ is a $T_0$-space. 
Observe that  $f'_*\tilde A$ is a quasi-coherent $\OO_{\tilde Y}$-module, then $f'_*\tilde A=\tilde {\mathcal A}$, where $\mathcal  A$ is a quasi-coherent $\OO_Y$-module and an $\OO_Y$-algebra.  Let us prove that $(Y,\A)$ is affine. We only have to prove that $(Y,\A_\pp)$ is affine for any $\pp\in\Spec A$. Given an open subset  $U\subset Y$, let  $I_U:=\{$quasi-compact open subsets
$\bar V\subset \tilde Y\colon \tilde U\subset \bar V\}$. Recall
$$\mathcal A(U)=\mathcal A_{|U}(U)=\widetilde{\A_{|U}}(\tilde U)\overset{\text{\ref{P16.17}}}=\tilde {\mathcal A}_{|\tilde U}(\tilde U)\overset{\text{\ref{coro25}}}=\ilim{\bar V\in I_U} \tilde{\mathcal A} (\bar V)=\ilim{\bar V\in I_U} \tilde A(f'^{-1}(\bar V)).
$$
Denote by $f'_\pp$ the composition of the morphisms $\Spec A_\pp\iny \Spec A\overset{f'}\to \tilde Y$.
Observe that
$$\A_{\pp}(U):=\mathcal A(U)_\pp=\ilim{\bar V\in I_U} \tilde A(f'^{-1}(\bar V))_\pp=\ilim{\bar V\in I_U} \tilde A_\pp({f'_\pp}^{-1}(\bar V)).$$
Then, we can suppose that $A=A_\pp$. The obvious morphism $\Id\colon (Y,\A)\to (Y,\OO_Y)$ is affine and $(U_{yy'},{\OO_Y}_{|U_{yy'}})$ is affine, for any $y,y'\in Y$, since $Y$ is semiseparated.
Hence, $(U_{yy'},\A_{|U_{yy'}})$ is affine. Then, the morphism
$\A_{yy'}\to \prod_{k\in U_{yy'}} \A_k$ is faithfully flat.
Let $y$ be the greatest point of $Y$ such that $f'(\pp)\in \tilde U_y$ and let $y'\in Y$ be another point. Observe that
$$\A_{y'}=\ilim{\bar V\in I_{U_{y'}}} \tilde{\mathcal A} (\bar V)
\overset*=\ilim{\bar V\in I_{U_{y'}},\bar W\in I_{U_{y}} } \tilde{\mathcal A} (\bar V\cap \bar W)\overset{\text{\ref{L16.18}}}=\ilim{\bar W\in I_{U_{yy'}}} \tilde{\mathcal A} (\bar W)\overset{\text{\ref{coro25}}}=\tilde\A_{|\tilde U_{yy'}}(\tilde U_{yy'})=
\A_{yy'}$$
($\overset*=$ since $f'^{-1}(\bar V)=f'^{-1}(\bar V\cap \bar W)$). 
Then, $\A_{y'}=\A_{yy'} \to \prod_{k\in U_{yy'}} \A_{k}$ is faithfully flat.
Hence, $y'$ is a removable point of $(Y,\A)$ if $y\not\leq y'$. 
Therefore, $(Y,\A)$ is quasi-isomorphic to $(U_y,\A_{|U_y})$, then it is affine.

Finally, $\Spec \A=\Spec \A(Y)=\Spec A=\tilde X$ and 
 the morphism $\tilde X=\Spec \A \to \Spec \OO_Y=\tilde Y$
 induced by the obvious morphism $\Id\colon (Y,\A)\to (Y,\OO_Y)$ is
 $f'$. Therefore, $\Hom_{[sch]}(X,Y)=\Hom_{sch}(\tilde X,\tilde Y)$.

2. Now, in general. 
$$\aligned \Hom_{[sch]}(X,Y) & =\Hom_{[sch]}(\ilim{x\in X} U_x,Y)
\overset{\text{\ref{T18.14}}}=
\plim{x\in X}\Hom_{[sch]}(U_x,Y)=\plim{x\in X}\Hom_{sch}(\tilde U_x,\tilde Y)\\ & \overset +=\Hom_{sch}(\ilim{x\in X} \tilde U_x,\tilde Y)=
\Hom_{sch}(\tilde X,\tilde Y)\endaligned$$
($\overset +=$ given $(f_x)\in \plim{x\in X}\Hom_{sch}(\tilde U_x,\tilde Y)$ the  induced morphism of ringed spaces $f\colon \ilim{x\in X} \tilde U_x\to \tilde Y$  is schematic, since for any quasi-coherent
$\OO_{\tilde X}$-module $\tilde \M=\plim{x\in X}\tilde i_{x*}\widetilde\M_x$, the $\OO_{\tilde Y}$-module $f_*\tilde \M=\plim{x\in X}f_*\tilde i_{x*}\widetilde\M_x=
\plim{x\in X}f_{x*}\widetilde\M_x$ is quasi-coherent).
\edemo

%\teor Let $S$ and $T$ be  quasicompact and quasi-separated 
%schemes. Let $$\U=\{U_1,\ldots,U_n\}\text{ and }\V=\{V_1,\ldots,V_n\}$$ be locally affine coverings of $S$ and $T$, and let $X$ and $Y$ be the schematic finite spaces associated with $\U$ and $\V$ respectively. Then, $$\Hom_{sch}(S,T)=\Hom_{[sch]}(X,Y).$$
%\eteor

\end{document}